\documentclass[english]{article}
\usepackage[T1]{fontenc}
\usepackage[utf8]{inputenc}
\usepackage{babel}
\usepackage{float}
\usepackage{enumitem}
\usepackage{dsfont}
\usepackage{amsmath}
\usepackage{amsthm}
\usepackage{amssymb}
\usepackage[authoryear,round]{natbib}
\usepackage[pdfusetitle,
 bookmarks=true,bookmarksnumbered=false,bookmarksopen=false,
 breaklinks=false,pdfborder={0 0 1},backref=false,colorlinks=false]
 {hyperref}

\makeatletter

\floatstyle{ruled}
\newfloat{algorithm}{tbp}{loa}
\providecommand{\algorithmname}{Algorithm}
\floatname{algorithm}{\protect\algorithmname}

\theoremstyle{remark}
\newtheorem{rem}{\protect\remarkname}
\theoremstyle{plain}
\newtheorem{assumption}{\protect\assumptionname}
\theoremstyle{plain}
\newtheorem{lem}{\protect\lemmaname}
\theoremstyle{definition}
 \newtheorem{example}{\protect\examplename}
\theoremstyle{plain}
\newtheorem{thm}{\protect\theoremname}
\theoremstyle{plain}
\newtheorem{cor}{\protect\corollaryname}
\newlist{casenv}{enumerate}{4}
\setlist[casenv]{leftmargin=*,align=left,widest={iiii}}
\setlist[casenv,1]{label={{\itshape\ \casename} \arabic*.},ref=\arabic*}
\setlist[casenv,2]{label={{\itshape\ \casename} \roman*.},ref=\roman*}
\setlist[casenv,3]{label={{\itshape\ \casename\ \alph*.}},ref=\alph*}
\setlist[casenv,4]{label={{\itshape\ \casename} \arabic*.},ref=\arabic*}

\@ifundefined{date}{}{\date{}}
\usepackage[paperwidth=8.5in, paperheight=11in, margin=1in]{geometry}

\usepackage{dsfont}

\makeatother

\providecommand{\assumptionname}{Assumption}
\providecommand{\casename}{Case}
\providecommand{\corollaryname}{Corollary}
\providecommand{\examplename}{Example}
\providecommand{\lemmaname}{Lemma}
\providecommand{\remarkname}{Remark}
\providecommand{\theoremname}{Theorem}

\begin{document}
\global\long\def\E{\mathbb{E}}%
\global\long\def\F{\mathcal{F}}%
\global\long\def\N{\mathbb{N}}%
\global\long\def\R{\mathbb{R}}%
\global\long\def\T{\mathcal{T}}%
\global\long\def\X{\mathbb{X}}%
\global\long\def\bv{\mathbf{v}}%
\global\long\def\bw{\mathbf{w}}%
\global\long\def\bx{\mathbf{x}}%
\global\long\def\by{\mathbf{y}}%
\global\long\def\bz{\mathbf{z}}%
\global\long\def\bzero{\mathbf{0}}%
\global\long\def\1{\mathds{1}}%
\global\long\def\defeq{\triangleq}%
\global\long\def\RR{\textsf{RR}}%
\global\long\def\SS{\textsf{SS}}%
\global\long\def\IG{\textsf{IG}}%
\global\long\def\SGD{\textsf{SGD}}%
\global\long\def\ShufflingSGD{\textsf{Shuffling SGD}}%
\global\long\def\RShufflingSGD{\textsf{Restarted Shuffling SGD}}%
\global\long\def\Lavg{\bar{L}}%
\global\long\def\Lmax{\hat{L}}%
\global\long\def\mydots{\dots}%
\global\long\def\Breg{\mathrm{B}}%
\global\long\def\Rate{\textsf{Rate}}%
\global\long\def\Order{\textsf{Order}}%
\global\long\def\Old{\textsf{Old}}%
\global\long\def\New{\textsf{New}}%
\global\long\def\Avg{\textsf{Avg}}%
\global\long\def\Last{\textsf{Last}}%
\global\long\def\Alg{\textsf{Alg}}%
\global\long\def\dist{\mathrm{dist}}%

\title{Random Reshuffling for Smooth Convex Optimization:\\Dominates Stochastic
Gradient Descent\thanks{A preliminary version \citep{pmlr-v336-liu26d}, titled \textit{Random
Reshuffling Dominates Stochastic Gradient Descent}, has been accepted
at COLT 2026. This extended version contains additional results. See
the end of the paper for a timeline and the AI disclosure.}}
\author{Zijian Liu\thanks{Stern School of Business, New York University, zl3067@stern.nyu.edu.}}
\maketitle
\begin{abstract}
Stochastic Gradient Descent ($\textsf{SGD}$) is one of the most classical
optimization algorithms with favorable theoretical guarantees, yet
its practical implementation differs subtly from its well-known form
and is often referred to as Shuffling Stochastic Gradient Descent
($\textsf{Shuffling SGD}$). A particularly popular strategy in $\textsf{Shuffling SGD}$
is Random Reshuffling ($\textsf{RR}$), which has achieved great empirical
success. Despite its strong performance, $\textsf{RR}$ has long been
considered a heuristic due to a lack of theoretical support.

Over the last decade, people have finally established provable convergence
rates for $\textsf{RR}$, thus justifying its observed superiority.
However, for smooth convex optimization, two clouds over the convergence
theory of $\textsf{RR}$ remain. According to current theory, $\textsf{Shuffling SGD}$
under $\textsf{RR}$ converges only when the stepsize is smaller than
a threshold proportional to $1/n$, where $n$ is the number of summands
in the objective (or the number of data points). Consequently, the
optimally tuned theoretical rate of $\textsf{Shuffling SGD}$ under
$\textsf{RR}$ is strictly worse than that of $\textsf{SGD}$ when
the number of epochs is smaller than another threshold proportional
to $n$. These two restrictions heavily limit the applicability of
existing theories and leave a critical mismatch with practice.

In this work, for the first time, we prove that $\textsf{RR}$ dominates
$\textsf{SGD}$ in smooth convex optimization for any reasonable stepsize
after any finite number of epochs, thereby addressing a longstanding
open question. Moreover, we establish new algorithm-dependent lower
bounds for $\textsf{RR}$ and $\textsf{SGD}$, respectively. Altogether,
these results yield (nearly) sharp worst-case characterizations of
$\textsf{RR}$ and $\textsf{SGD}$. Furthermore, we develop a shuffling-based
method to optimize functions satisfying the quadratic growth condition,
a strict relaxation of strong convexity, and derive a new complexity
bound improving existing results proved under stronger assumptions.
\end{abstract}

\section{Introduction\label{sec:Introduction}}

One of the fundamental tasks in machine learning is to optimize functions
in a finite-sum form, i.e., $f(\bx)\defeq\frac{1}{n}\sum_{i=1}^{n}f_{i}(\bx)$.
Among different optimization algorithms, Stochastic Gradient Descent
($\SGD$), proposed in the seminal work of \citet{10.1214/aoms/1177729586},
is arguably one of the most classical methods. Due to its easy implementation
and computational efficiency, $\SGD$ is particularly popular when
$n$ is large, the standard case nowadays. More importantly, the convergence
guarantees of $\SGD$ have been extensively studied, yielding provable
rates in various settings \citep{Solr-KOHA-OAI-TEST:19722,doi:10.1137/16M1080173,lan2020first},
thereby providing a theoretical backbone for $\SGD$.

However, compared with the standard form of $\SGD$ analyzed in theory,
which uniformly samples a function to perform a gradient descent step
at each iteration (also known as with-replacement sampling), the practical
implementation differs subtly and is often referred to as Shuffling
Stochastic Gradient Descent ($\ShufflingSGD$). In $\ShufflingSGD$,
the optimization procedure is divided into $K$ epochs, and within
each epoch, the order in which functions are processed is determined
by a permutation $\pi$ of $\left\{ 1,\mydots,n\right\} $ (also known
as without-replacement sampling). A widely implemented strategy for
generating $\pi$ is Random Reshuffling ($\RR$), which, in each epoch,
independently and uniformly draws a new permutation from all possible
ones.

Although $\ShufflingSGD$ under $\RR$ has achieved great empirical
success across numerous experiments, it has long been considered a
heuristic due to a lack of theoretical support. Over the last decade,
beginning with the pioneering work of \citet{gurbuzbalaban2021random},
people have finally established provable convergence rates for $\RR$,
thereby justifying its observed superiority over standard $\SGD$.

In particular, for smooth convex optimization (i.e., each $f_{i}$
is convex and $L$-smooth\footnote{For simplicity, we adopt a uniform smoothness parameter in the discussion,
as in most of the existing literature.}), $\RR$ with a constant stepsize\footnote{We also use a constant stepsize in the discussion for convenience.}
$\eta$ is known to converge in expectation at the rate $\frac{D^{2}}{\eta nK}+\eta^{2}nL\sigma_{\star}^{2}$
(e.g., \citet{NEURIPS2020_c8cc6e90,JMLR:v22:20-1238}), provided that
the stepsize satisfies $\eta\lesssim\frac{1}{nL}$, where $D$ denotes
the distance between the initial point and the optimal solution, and
$\sigma_{\star}^{2}$ is the gradient variance at the minimizer. In
comparison, $\SGD$ under the same setting guarantees an in-expectation
convergence rate $\frac{D^{2}}{\eta nK}+\eta\sigma_{\star}^{2}$ \citep{garrigos2023handbook},
but only requires $\eta\leq\frac{1}{4L}$. Recently, two works \citep{cortild2025bias,NEURIPS2025_97335d45}
have further extended the convergence theory for $\SGD$ to the regime
$\eta<\frac{2}{L}$. To date, the best-known rate of $\SGD$ for any
$\eta<\frac{2}{L}$ is $\frac{D^{2}}{(2-\eta L)\eta nK}+\frac{\eta\sigma_{\star}^{2}}{(2-\eta L)^{2}}$
\citep{NEURIPS2025_97335d45}. Clearly, $\RR$ converges faster than
$\SGD$ in the regime $\eta\lesssim\frac{1}{nL}$, which has been
recognized as theoretical evidence demonstrating the strong performance
of $\RR$.

Despite the progress discussed above, some important issues remain
unaddressed. The most critical longstanding open question is that
people still do not understand what happens to $\RR$ when $\eta$
falls into the regime $\eta\gtrsim\frac{1}{nL}$. This point is critical
because, given that $n$ is typically large in modern tasks, the threshold
$\frac{1}{nL}$ can be extremely small or even vanish, whereas stepsizes
used in practice are usually at a constant level, thereby leaving
a significant gap between theory and practice. More crucially, even
if one temporarily assumes that the existing rate for $\RR$ mentioned
earlier could be extended to allow $\eta<\frac{2}{L}$ (though no
such theory has been established), it would still fail to explain
the advantage of $\RR$ over $\SGD$, as the term $\eta^{2}nL\sigma_{\star}^{2}$
for $\RR$ is worse than the term $\eta\sigma_{\star}^{2}$ for $\SGD$
when $\eta\gtrsim\frac{1}{nL}$. This hints that an analysis different
from existing ones may be needed.

Another issue implied by the above discussion is that the optimally
tuned rate for $\RR$ induced by the existing bound is only $\frac{LD^{2}}{K}+(\frac{L\sigma_{\star}^{2}D^{4}}{nK^{2}})^{\frac{1}{3}}$,
which is better than the best tuned rate $\frac{LD^{2}}{nK}+\frac{\sigma_{\star}D}{\sqrt{nK}}$
for $\SGD$ only when $K$ is larger than a threshold proportional
to $n$. Moreover, in the case of $\sigma_{\star}=0$ (i.e., all $f_{i}$'s
share a common optimal solution), the rate of $\RR$ reduces to only
$\frac{LD^{2}}{K}$, which is even worse than the $\frac{LD^{2}}{nK}$
rate of $\SGD$ by a factor of $n$.

The above restrictions on the stepsize $\eta$ or the number of epochs
$K$ heavily limit the applicability of existing theories and cannot
fully explain the favorable performance of $\ShufflingSGD$ under
$\RR$ compared with standard $\SGD$. Therefore, we are naturally
led to the following question:
\begin{center}
\textit{In smooth convex optimization, does $\RR$ dominate $\SGD$
without these two restrictions?}
\par\end{center}

\subsection{Our Contributions}

This work provides an affirmative answer to the above question.
\begin{itemize}
\item Concretely, we show that $\ShufflingSGD$ under $\RR$ provably converges
at a rate of $\frac{D^{2}}{(2-\eta L)\eta nK}+\frac{\min\left\{ 1,\eta nL\right\} \eta\sigma_{\star}^{2}}{(2-\eta L)^{2}}$
for any stepsize satisfying $\eta<\frac{2}{L}$ (see Theorem \ref{thm:RR-cvx-ub}
for the formal version with nonuniform smoothness parameters and dynamic
stepsizes that depend on the epoch number). We highlight that this
rate is not only the first provable result for $\RR$ that allows
$\eta\gtrsim\frac{1}{nL}$, but also provably dominates the bound
$\frac{D^{2}}{(2-\eta L)\eta nK}+\frac{\eta\sigma_{\star}^{2}}{(2-\eta L)^{2}}$
of $\SGD$ under any stepsize $\eta<\frac{2}{L}$. It is noteworthy
that this rate is not merely a simple extension of the previously
best bound for $\RR$, since the latter is slower than $\SGD$ when
$\eta\gtrsim\frac{1}{nL}$, as discussed before. Therefore, for the
first time, we prove that $\RR$ dominates $\SGD$ in smooth convex
optimization for any reasonable stepsize after any finite number of
epochs, resolving a longstanding open question. Consequently, the
optimally tuned rate for $\ShufflingSGD$ under $\RR$ is improved
to $\frac{LD^{2}}{nK}+\min\{\frac{\sigma_{\star}D}{\sqrt{nK}},(\frac{L\sigma_{\star}^{2}D^{4}}{nK^{2}})^{\frac{1}{3}}\}$
(see Corollary \ref{cor:RR-cvx-optimal} for the formal version with
nonuniform smoothness parameters), which dominates the best tuned
bound $\frac{LD^{2}}{nK}+\frac{\sigma_{\star}D}{\sqrt{nK}}$ of $\SGD$
for any finite $K$. Moreover, the rate reduces to $\frac{LD^{2}}{nK}$
when $\sigma_{\star}=0$ (i.e., the interpolation regime), improving
upon the best known result $\frac{LD^{2}}{K}$ by a factor of $\frac{1}{n}$.
\item In addition, we prove new algorithm-dependent lower bounds of $\frac{D^{2}}{\eta nK}+\min\left\{ 1,\eta nL\right\} \eta\sigma_{\star}^{2}$
for $\ShufflingSGD$ under $\RR$ in Theorem \ref{thm:RR-cvx-lb}
and $\frac{D^{2}}{\eta nK}+\eta\sigma_{\star}^{2}$ for $\SGD$ in
Theorem \ref{thm:SGD-cvx-lb}. Together with our new upper bound for
$\RR$ and the existing upper bound for $\SGD$, these results establish
(nearly) sharp worst-case characterizations of the convergence of
$\RR$ and $\SGD$ in the regime $\eta\leq\frac{1}{L}$ (and, more
generally, $\eta\leq\frac{C}{L}$ for any absolute constant $C\in\left(0,2\right)$).
\item Furthermore, building on the classical restart technique in the optimization
literature, we develop a restarted variant of $\ShufflingSGD$ termed
$\RShufflingSGD$. By applying the new upper bound for $\RR$ mentioned
above, we show that $\RShufflingSGD$ under $\RR$ can be used to
optimize functions satisfying the quadratic growth condition, which
is strictly weaker than strong convexity. In particular, we establish
a new complexity bound under this weaker assumption, improving upon
the best known results for $\RR$ that require either strong convexity
of $f$ or strong convexity of each $f_{i}$. Moreover, our new bound
is sharper than the corresponding bound for $\SGD$ across different
parameter regimes under the same setting.
\end{itemize}

\subsection{Related Work}

We provide a brief overview of $\ShufflingSGD$ under $\RR$ and defer
further details to Appendix \ref{sec:Additional-Related-Work}.

Over the past decades, the effectiveness of $\RR$ has been reported
in many works (e.g., \citet{bottou2009curiously,Bottou2012,Bengio2012}).
However, the theoretical understanding of it has long lagged behind.
The first breakthrough is by \citet{gurbuzbalaban2021random}, which
provides the first theoretical evidence that $\RR$ can beat $\SGD$
in smooth strongly convex optimization under certain additional assumptions.
Since then, $\RR$ has been extensively studied (e.g., \citet{pmlr-v97-nagaraj19a,pmlr-v97-haochen19a}).
To date, the best known rate in smooth convex optimization is $\frac{D^{2}}{\eta nK}+\eta^{2}nL\sigma_{\star}^{2}$
under the condition $\eta\lesssim\frac{1}{nL}$ \citep{NEURIPS2020_c8cc6e90,JMLR:v22:20-1238}.
The only existing lower bound in smooth convex optimization is $(\frac{L\sigma_{\star}^{2}D^{4}}{nK^{2}})^{\frac{1}{3}}$,
due to \citet{pmlr-v202-cha23a}, which holds for constant stepsizes
$\eta\lesssim\frac{1}{nL}$ and a large number of epochs at least
satisfying $K\gtrsim\frac{nL^{2}D^{2}}{\sigma_{\star}^{2}}$.

\section{Preliminaries}

\paragraph{Notation.}

$\N$ denotes the set of natural numbers (excluding $0$). Given $n\in\N$,
we write $\left[n\right]\defeq\left\{ 1,\mydots,n\right\} $. $\left\lfloor \cdot\right\rfloor $
and $\left\lceil \cdot\right\rceil $ denote the floor and ceiling
functions, respectively. $\left\langle \cdot,\cdot\right\rangle $
is the standard Euclidean inner product, and $\left\Vert \cdot\right\Vert \defeq\sqrt{\left\langle \cdot,\cdot\right\rangle }$
is the $\ell_{2}$ norm. Given a real-valued differentiable function
$h:\R^{d}\to\R$, $\nabla h(\bx)$ denotes the gradient at $\bx\in\R^{d}$.
The Bregman divergence induced by $h$ is defined as $\Breg_{h}(\bx,\by)\defeq h(\bx)-h(\by)-\left\langle \nabla h(\by),\bx-\by\right\rangle $,
which is nonnegative if $h$ is additionally convex.

\paragraph{Objective.}

We study the following finite-sum optimization problem in this work
\[
\inf_{\bx\in\R^{d}}f(\bx)\defeq\frac{1}{n}\sum_{i=1}^{n}f_{i}(\bx),
\]
where $n\in\N$ and each $f_{i}:\R^{d}\to\R$ is differentiable.
\begin{rem}
To ease notation, hereinafter we use $\Breg\defeq\Breg_{f}$ and $\Breg_{i}\defeq\Breg_{f_{i}}$
to denote the Bregman divergences induced by $f$ and $f_{i}$, respectively. 
\end{rem}

\paragraph{Assumptions.}

Our analysis relies on the following three assumptions.
\begin{assumption}[Minimizer]
\label{assu:minimizer}$\exists\bx_{\star}\in\R^{d}$ such that $f_{\star}\defeq f(\bx_{\star})=\inf_{\bx\in\R^{d}}f(\bx)\in\R$.
\end{assumption}
\begin{assumption}[Convexity]
\label{assu:convexity}Each $f_{i}$ is convex.
\end{assumption}
\begin{assumption}[Smoothness]
\label{assu:smoothness}Each $f_{i}$ is $L_{i}$-smooth, i.e., $\exists L_{i}>0$
such that $\left\Vert \nabla f_{i}(\bx)-\nabla f_{i}(\by)\right\Vert \leq L_{i}\left\Vert \bx-\by\right\Vert ,\forall\bx,\by\in\R^{d}$.
\end{assumption}
All three of the above assumptions are standard and commonly adopted
in the literature \citep{Solr-KOHA-OAI-TEST:19722,nesterov2018lectures,doi:10.1137/16M1080173,lan2020first}.
Notably, we do not impose any assumptions on the difference between
the individual gradient $\nabla f_{i}$ and the full gradient $\nabla f$,
such as the popular finite variance condition.

Next, we introduce three more notations to simplify the expressions
in the subsequent sections. $\sigma_{\star}^{2}$ denotes the variance
of the gradient at the minimizer $\bx_{\star}$, and $\Lavg$ (resp.
$\Lmax$) represents the average (resp. maximum) smoothness parameter,
i.e.,
\begin{eqnarray*}
\sigma_{\star}^{2}\defeq\frac{1}{n}\sum_{i=1}^{n}\left\Vert \nabla f_{i}(\bx_{\star})\right\Vert ^{2}, & \Lavg\defeq\frac{1}{n}\sum_{i=1}^{n}L_{i}, & \Lmax\defeq\max_{i\in\left[n\right]}L_{i}.
\end{eqnarray*}
We note that the quantity $\sigma_{\star}^{2}$ is widely used in
prior works on shuffling gradient methods (e.g., \citet{8514028,NEURIPS2020_c8cc6e90,JMLR:v22:20-1238})
and remains invariant even when $f$ has multiple minimizers (see
Lemma 4.17 of \citet{garrigos2023handbook}). In particular, $\sigma_{\star}^{2}=0$
corresponds to the case in which all $f_{i}$ share a common optimal
solution.

To finish this section, we state a classical result in convex optimization,
known as the co-coercivity property of smooth convex functions, which
serves as a key tool in our analysis. As for its proof, see, for example,
Theorem 2.15 of \citet{nesterov2018lectures}.
\begin{lem}[Co-coercivity]
\label{lem:co-coercivity-cvx}Let $h:\R^{d}\to\R$ be a differentiable
convex function that is also $L$-smooth, then we have, for any $\bx,\by\in\R^{d}$,
\begin{eqnarray*}
\left\Vert \nabla h(\bx)-\nabla h(\by)\right\Vert ^{2}\leq2L\Breg_{h}(\bx,\by) & \text{and} & \left\Vert \nabla h(\bx)-\nabla h(\by)\right\Vert ^{2}\leq L\left\langle \nabla h(\bx)-\nabla h(\by),\bx-\by\right\rangle .
\end{eqnarray*}
\end{lem}

\section{Shuffling Stochastic Gradient Descent}

\begin{algorithm}[H]
\caption{\label{alg:Shuffling-SGD}Shuffling Stochastic Gradient Descent ($\protect\ShufflingSGD$)}

\textbf{Input:} initial point $\bx_{1}^{1}\in\R^{d}$, stepsize $\eta_{k}>0$,
number of epochs $K\in\N$

\textbf{for} $k=1$ \textbf{to} $K$ \textbf{do}

$\quad$Generate a permutation $\pi_{k}$ of $\left[n\right]$

$\quad$\textbf{for} $i=1$ \textbf{to} $n$ \textbf{do}

$\qquad$$\bx_{k}^{i+1}=\bx_{k}^{i}-\eta_{k}\nabla f_{\pi_{k}^{i}}(\bx_{k}^{i})$

$\quad$\textbf{end for}

$\quad$$\bx_{k+1}^{1}=\bx_{k}^{n+1}$

\textbf{end for}
\end{algorithm}

The method studied in this work, Shuffling Stochastic Gradient Descent
($\ShufflingSGD$), is given in Algorithm \ref{alg:Shuffling-SGD}.
Compared with the standard $\SGD$ algorithm, which uniformly samples
a function to process at each step, $\ShufflingSGD$ determines the
order in which functions are passed in each epoch based on a permutation.
In particular, three strategies for generating permutations are popular
in practice, as illustrated in the following examples, where $S_{n}$
denotes the symmetric group of $\left[n\right]$.
\begin{example}[Random Reshuffling ($\RR$)]
\label{exa:RR}Each $\pi_{k}$ is drawn independently and uniformly
from $S_{n}$.
\end{example}
\begin{example}[Single Shuffling ($\SS$)]
\label{exa:SS}Each $\pi_{k}=\pi$, a permutation drawn uniformly
from $S_{n}$.
\end{example}
\begin{example}[Incremental Gradient ($\IG$)]
\label{exa:IG}Each $\pi_{k}=\pi$, a deterministic permutation from
$S_{n}$. 
\end{example}

\subsection{New Rate for $\protect\RR$}

We are now ready to provide the first main result, Theorem \ref{thm:RR-cvx-ub},
a new convergence rate for $\RR$.
\begin{thm}
\label{thm:RR-cvx-ub}Under Assumptions \ref{assu:minimizer}, \ref{assu:convexity},
and \ref{assu:smoothness}, suppose $\RR$ is employed with $\left\{ \eta_{k}<\frac{2}{\Lmax}\right\} _{k=1}^{K}$,
let
\begin{eqnarray}
\bar{\bx}_{K}\defeq\sum_{k=1}^{K}\sum_{i=1}^{n}\frac{(2-\eta_{k}\Lmax)\eta_{k}}{nH_{K}}\bx_{k}^{i} & \text{where} & H_{K}\defeq\sum_{k=1}^{K}(2-\eta_{k}\Lmax)\eta_{k},\label{eq:avg-x}
\end{eqnarray}
then $\ShufflingSGD$ (Algorithm \ref{alg:Shuffling-SGD}) guarantees
that
\[
\E\left[f(\bar{\bx}_{K})-f_{\star}\right]\lesssim\frac{\left\Vert \bx_{1}^{1}-\bx_{\star}\right\Vert ^{2}}{n\sum_{k=1}^{K}(2-\eta_{k}\Lmax)\eta_{k}}+\min\left\{ \frac{\sum_{k=1}^{K}\eta_{k}^{2}/(2-\eta_{k}\Lmax)}{\sum_{k=1}^{K}(2-\eta_{k}\Lmax)\eta_{k}},\frac{\sum_{k=1}^{K}\eta_{k}^{3}n\Lavg/(2-\eta_{k}\Lmax)}{\sum_{k=1}^{K}(2-\eta_{k}\Lmax)\eta_{k}}\right\} \sigma_{\star}^{2}.
\]
\end{thm}
\begin{rem}
As in the conference version \citep{pmlr-v336-liu26d}, we can also
consider the weighted average iterate $\bar{\bx}_{K}\defeq\sum_{k=1}^{K}\sum_{i=1}^{n}\frac{\eta_{k}}{nH_{K}}\bx_{k}^{i}$,
where $H_{K}\defeq\sum_{k=1}^{K}\eta_{k}$, and prove a similar rate
to that in Theorem \ref{thm:RR-cvx-ub}, with $2-\eta_{k}\Lmax$ replaced
by $2-\hat{\eta}\Lmax$, where $\hat{\eta}\defeq\max_{k\in\left[K\right]}\eta_{k}$.
\end{rem}
\begin{proof}
The proof is deferred to Subsection \ref{subsec:RR-cvx-ub-Proof}.
\end{proof}

To the best of our knowledge, in smooth convex optimization, Theorem
\ref{thm:RR-cvx-ub} shows the first convergence rate for $\RR$ under
any reasonable stepsize, i.e., $\eta_{k}<2/\Lmax$, thereby improving
existing results in different aspects. Previously, the best known
convergence rate for $\RR$ is $\frac{D^{2}}{n\sum_{k=1}^{K}\eta_{k}}+\frac{\sum_{k=1}^{K}\eta_{k}^{3}n\Lavg}{\sum_{k=1}^{K}\eta_{k}}\sigma_{\star}^{2}$
(where $D\defeq\left\Vert \bx_{1}^{1}-\bx_{\star}\right\Vert $) under
the condition $\eta_{k}\lesssim1/(n\sqrt{\Lavg\Lmax})$ \citep{pmlr-v235-liu24cg,NEURIPS2024_84d39572}.
However, this requirement on $\eta_{k}$ is highly unrealistic, as
$n$ is typically pretty large in modern machine learning tasks, meaning
that the stepsize has to be extremely small or even vanishing. This
contradicts the constant-level stepsizes commonly used in practice
and leads to a gap between theory and practice.

In comparison, Theorem \ref{thm:RR-cvx-ub} not only allows the stepsize
to lie in the constant regime $\eta_{k}<2/\Lmax$ but also shows a
fundamental improvement rather than a mere extension of the previously
best known rate, since the latter becomes slower than $\SGD$ once
$\eta_{k}\gtrsim1/(n\Lavg)$ (the rate of $\SGD$ is stated below),
while our Theorem \ref{thm:RR-cvx-ub} never does. 

In addition, compared with the rate $\frac{D^{2}}{n\sum_{k=1}^{K}(2-\eta_{k}\Lmax)\eta_{k}}+\frac{\sum_{k=1}^{K}\eta_{k}^{2}/(2-\eta_{k}\Lmax)}{\sum_{k=1}^{K}(2-\eta_{k}\Lmax)\eta_{k}}\sigma_{\star}^{2}$
of standard $\SGD$ \citep{NEURIPS2025_97335d45}\footnote{The work of \citet{NEURIPS2025_97335d45} does not provide a rate
for dynamic stepsizes. However, their proof can be directly extended
to this setting. For completeness, we include a separate theorem and
prove it in Appendix \ref{sec:SGD-cvx-ub}.}, our Theorem \ref{thm:RR-cvx-ub} is never worse for any reasonable
stepsize (i.e., $\eta_{k}<2/\Lmax$), and is strictly better when
$\eta_{k}\lesssim1/(n\Lavg)$. This feature has an important implication:
if both methods employ their own optimally tuned stepsizes, $\ShufflingSGD$
under $\RR$ provably achieves a better upper bound that dominates
the $\frac{\Lmax D^{2}}{nK}+\frac{\sigma_{\star}D}{\sqrt{nK}}$ rate
of $\SGD$ after any finite number of epochs, as evidenced by Corollary
\ref{cor:RR-cvx-optimal} below.
\begin{cor}
\label{cor:RR-cvx-optimal}Under the same setting as in Theorem \ref{thm:RR-cvx-ub},
with the optimally tuned constant stepsize $\left\{ \eta_{k}=\eta_{\star}\right\} _{k=1}^{K}$,
where $\eta_{\star}<\frac{2}{\Lmax}$, $\ShufflingSGD$ (Algorithm
\ref{alg:Shuffling-SGD}) guarantees that
\[
\E\left[f(\bar{\bx}_{K})-f_{\star}\right]\lesssim\frac{\Lmax\left\Vert \bx_{1}^{1}-\bx_{\star}\right\Vert ^{2}}{nK}+\min\left\{ \frac{\sigma_{\star}\left\Vert \bx_{1}^{1}-\bx_{\star}\right\Vert }{\sqrt{nK}},\left(\frac{\Lavg\sigma_{\star}^{2}\left\Vert \bx_{1}^{1}-\bx_{\star}\right\Vert ^{4}}{nK^{2}}\right)^{\frac{1}{3}}\right\} .
\]
\end{cor}
\begin{proof}
The proof is deferred to Subsection \ref{subsec:RR-cvx-ub-Proof}.
\end{proof}

In contrast, the known optimally tuned rate for $\RR$ is only $\frac{\sqrt{\Lavg\Lmax}D^{2}}{K}+(\frac{\Lavg\sigma_{\star}^{2}D^{4}}{nK^{2}})^{\frac{1}{3}}$.
To better understand the differences among these optimally tuned rates,
let $\Order_{\SGD}(K)$ and $\Order_{\RR}^{\Old}(K)$ denote the dominant
terms in the existing optimally tuned rates for $\SGD$ and $\ShufflingSGD$
under $\RR$, respectively, i.e.,
\begin{eqnarray*}
\Order_{\SGD}(K)\defeq\max\left\{ \frac{\Lmax D^{2}}{nK},\frac{\sigma_{\star}D}{\sqrt{nK}}\right\}  & \text{and} & \Order_{\RR}^{\Old}(K)\defeq\max\left\{ \frac{\sqrt{\Lavg\Lmax}D^{2}}{K},\left(\frac{\Lavg\sigma_{\star}^{2}D^{4}}{nK^{2}}\right)^{\frac{1}{3}}\right\} .
\end{eqnarray*}
Similarly, $\Order_{\RR}^{\New}(K)$ denotes the dominant term in
the rate obtained in Corollary \ref{cor:RR-cvx-optimal}, i.e.,
\[
\Order_{\RR}^{\New}(K)\defeq\max\left\{ \frac{\Lmax D^{2}}{nK},\min\left\{ \frac{\sigma_{\star}D}{\sqrt{nK}},\left(\frac{\Lavg\sigma_{\star}^{2}D^{4}}{nK^{2}}\right)^{\frac{1}{3}}\right\} \right\} .
\]

In the nondegenerate case $\sigma_{\star}\neq0$ (which implies that
$n\geq2$), the following comparisons hold:
\[
\begin{cases}
\Order_{\RR}^{\New}(K)=\Order_{\SGD}(K)\overset{(a)}{\leq}\Order_{\RR}^{\Old}(K), & K\leq\frac{n\Lavg^{2}D^{2}}{\sigma_{\star}^{2}},\\
\Order_{\RR}^{\New}(K)<\Order_{\SGD}(K)\overset{(b)}{\leq}\Order_{\RR}^{\Old}(K), & \frac{n\Lavg^{2}D^{2}}{\sigma_{\star}^{2}}<K\leq\frac{n\Lavg\Lmax D^{2}}{\sigma_{\star}^{2}},\\
\Order_{\RR}^{\New}(K)\overset{(c)}{\leq}\Order_{\RR}^{\Old}(K)<\Order_{\SGD}(K), & \frac{n\Lavg\Lmax D^{2}}{\sigma_{\star}^{2}}<K\leq\frac{n\Lavg^{1/2}\Lmax^{3/2}D^{2}}{\sigma_{\star}^{2}},\\
\Order_{\RR}^{\New}(K)=\Order_{\RR}^{\Old}(K)<\Order_{\SGD}(K), & K>\frac{n\Lavg^{1/2}\Lmax^{3/2}D^{2}}{\sigma_{\star}^{2}},
\end{cases}
\]
where $(a)$ becomes an equality if and only if $K=\frac{n\Lavg^{2}D^{2}}{\sigma_{\star}^{2}}$
and $\Lavg=\Lmax$, $(b)$ becomes an equality if and only if $K=\frac{n\Lavg\Lmax D^{2}}{\sigma_{\star}^{2}}$,
and $(c)$ becomes an equality if and only if $K=\frac{n\Lavg^{1/2}\Lmax^{3/2}D^{2}}{\sigma_{\star}^{2}}$.
As one can see, unlike the existing optimally tuned bound for $\RR$
in the literature, which can be slower than standard $\SGD$ when
$K\leq\frac{n\Lavg\Lmax D^{2}}{\sigma_{\star}^{2}}$, the new result
in Corollary \ref{cor:RR-cvx-optimal} is never worse than the optimally
tuned rate of standard $\SGD$ and strictly beats it once $K>\frac{n\Lavg^{2}D^{2}}{\sigma_{\star}^{2}}$,
thereby improving the threshold from $\frac{n\Lavg\Lmax D^{2}}{\sigma_{\star}^{2}}$
to $\frac{n\Lavg^{2}D^{2}}{\sigma_{\star}^{2}}$.

Moreover, in the special case where all $f_{i}$'s share a common
minimizer, or equivalently when $\sigma_{\star}=0$, Corollary \ref{cor:RR-cvx-optimal}
achieves the same $\frac{\Lmax D^{2}}{nK}$ rate as $\SGD$. However,
the prior bound for $\RR$ can only be reduced to a slower rate in
the order of $\frac{\sqrt{\Lavg\Lmax}D^{2}}{K}$ due to the stepsize
restriction $\eta_{k}\lesssim1/(n\sqrt{\Lavg\Lmax})$, as discussed
before. Formally, once $n\geq2$, we have
\[
\Order_{\RR}^{\New}(K)=\Order_{\SGD}(K)<\Order_{\RR}^{\Old}(K),\forall K\in\N.
\]

In summary, Theorem \ref{thm:RR-cvx-ub} and Corollary \ref{cor:RR-cvx-optimal}
together imply that $\RR$ dominates $\SGD$ in smooth convex optimization
under any reasonable stepsize after any finite number of epochs.

\subsection{New Lower Bounds for $\protect\RR$ and $\protect\SGD$}

The careful reader may notice that our previous comparison focuses
only on upper bounds, which alone may not fully substantiate the claim
that $\RR$ dominates $\SGD$. Motivated by this, we establish new
algorithm-dependent lower bounds for both $\RR$ and standard $\SGD$
under a constant stepsize $\left\{ \eta_{k}=\eta\right\} _{k=1}^{K}$.

Before stating these lower bounds, we introduce some notation for
clarity. Given $d\in\N$, $n\geq2$, $D>0$, $L_{1:n}\defeq\left\{ L_{i}>0\right\} _{i=1}^{n}$,
and $\sigma_{\star}\geq0$, we let
\[
\F(d,n,D,L_{1:n},\sigma_{\star})\defeq\left\{ f=\frac{1}{n}\sum_{i=1}^{n}f_{i}:\begin{array}{c}
\exists\bx_{\star}\in\R^{d}\text{ such that }f_{\star}\defeq f(\bx_{\star})=\inf_{\bx\in\R^{d}}f(\bx)\text{ and }\left\Vert \bx_{\star}\right\Vert \leq D,\\
\text{each }f_{i}:\R^{d}\to\R\text{ is differentiable, convex, and }L_{i}\text{-smooth},\\
\frac{1}{n}\sum_{i=1}^{n}\left\Vert \nabla f_{i}(\bx_{\star})\right\Vert ^{2}\leq\sigma_{\star}^{2}.
\end{array}\right\} .
\]
In words, the class $\F(d,n,D,L_{1:n},\sigma_{\star})$ contains all
finite-sum objectives satisfying the required assumptions.

We are now ready to present the new lower bounds. The first is Theorem
\ref{thm:RR-cvx-lb}, a new lower bound for $\ShufflingSGD$ with
$\RR$.
\begin{thm}
\label{thm:RR-cvx-lb}Given $d\in\N$, $n\geq2$, $D>0$, $\left\{ L_{i}>0\right\} _{i=1}^{n}$,
and $\sigma_{\star}\geq0$, suppose $\RR$ is employed with $\left\{ \eta_{k}=\eta\right\} _{k=1}^{K}$
for $\eta>0$, let $\bar{\bx}_{K}\defeq\frac{1}{nK}\sum_{k=1}^{K}\sum_{i=1}^{n}\bx_{k}^{i}$
denote the average iterate generated by $\ShufflingSGD$ (Algorithm
\ref{alg:Shuffling-SGD}), initialized at $\bx_{1}^{1}=\bzero$, over
$K\in\N$ epochs, then we have
\begin{eqnarray*}
\sup_{f\in\F(d,n,D,L_{1:n},\sigma_{\star})}\E\left[f(\bar{\bx}_{K})-f_{\star}\right]\gtrsim\frac{D^{2}}{\eta nK}+\min\left\{ 1,\eta n\Lavg\right\} \eta\sigma_{\star}^{2} & \text{when} & K\eta n\Lavg\geq1.
\end{eqnarray*}
\end{thm}
\begin{rem}
\label{rem:RR-cvx-lb-condition}Although the theorem is stated under
the condition $K\eta n\Lavg\geq1$, our proof applies to any $\eta>0$
and $K\in\N$. We restrict our attention to the regime $K\eta n\Lavg\geq1$,
since our analysis indicates that, in the regime $K\eta n\Lavg<1$,
$\sup_{f\in\F(d,n,D,L_{1:n},\sigma_{\star})}\E\left[f(\bar{\bx}_{K})-f_{\star}\right]\gtrsim\Lavg D^{2}$,
suggesting that no meaningful convergence occurs for $\RR$ in this
situation.
\end{rem}
\begin{proof}
The proof is deferred to Subsection \ref{subsec:RR-cvx-lb-Proof}.
\end{proof}

Theorem \ref{thm:RR-cvx-lb} provides a new algorithm-dependent lower
bound for $\RR$, extending the existing lower bound $(\frac{L\sigma_{\star}^{2}D^{4}}{nK^{2}})^{\frac{1}{3}}$
for $\left\{ L_{i}=L\right\} _{i=1}^{n}$ under the conditions $\eta\lesssim1/(nL)$
and $K\gtrsim nL^{2}D^{2}/\sigma_{\star}^{2}$ \citep{pmlr-v202-cha23a}
to every prespecified $\eta>0$ and $K\in\N$ for nonuniform $\left\{ L_{i}\right\} _{i=1}^{n}$
satisfying $K\eta n\Lavg\geq1$, the meaningful regime discussed in
Remark \ref{rem:RR-cvx-lb-condition}. Importantly, this allows a
direct comparison with stepsize-dependent upper bounds and provides
a finer characterization of the limitations of $\RR$. As one can
observe, there remains a gap between the upper and lower bounds by
factors polynomial in $1/(2-\eta\Lmax)$ when $\eta<2/\Lmax$, and
resolving this discrepancy is left for future work.

We highlight that Theorem \ref{thm:RR-cvx-lb} accommodates nonuniform
smoothness parameters, therefore requiring a more careful analysis.
At a high level, the analysis is divided into two parts. In the first
part, we establish the term $\frac{D^{2}}{\eta nK}$ to capture the
error induced by the initial distance between $\bx_{1}^{1}$ and $\bx_{\star}$.
In the second part, inspired by \citet{pmlr-v119-rajput20a,yun2022minibatch,pmlr-v202-cha23a},
we also construct asymmetric piecewise quadratic functions that allow
us to explicitly track the $\RR$ iterates and establish the second
term $\min\left\{ 1,\eta n\Lavg\right\} \eta\sigma_{\star}^{2}$,
which captures the lower bound arising from $\sigma_{\star}^{2}$.
For more details, please refer to Section \ref{sec:RR-cvx-lb}.
\begin{thm}
\label{thm:SGD-cvx-lb}Given $d\in\N$, $n\geq2$, $D>0$, $\left\{ L_{i}>0\right\} _{i=1}^{n}$,
and $\sigma_{\star}\geq0$, suppose a constant stepsize $\left\{ \eta_{k}=\eta\right\} _{k=1}^{K}$
for $\eta>0$ is employed, let $\bar{\bx}_{K}\defeq\frac{1}{nK}\sum_{k=1}^{K}\sum_{i=1}^{n}\bx_{k}^{i}$
denote the average iterate generated by $\SGD$ (Algorithm \ref{alg:Shuffling-SGD},
with $\pi_{k}^{i}$ sampled independently and uniformly from $\left[n\right]$),
initialized at $\bx_{1}^{1}=\bzero$, over $K\in\N$ epochs, then
we have
\begin{eqnarray*}
\sup_{f\in\F(d,n,D,L_{1:n},\sigma_{\star})}\E\left[f(\bar{\bx}_{K})-f_{\star}\right]\gtrsim\frac{D^{2}}{\eta nK}+\eta\sigma_{\star}^{2} & \text{when} & K\eta n\Lavg\geq1.
\end{eqnarray*}
\end{thm}
\begin{rem}
Similar to Theorem \ref{thm:RR-cvx-lb}, we focus on the regime $K\eta n\Lavg\geq1$,
since our proof suggests that, for $\SGD$ as well, $\sup_{f\in\F(d,n,D,L_{1:n},\sigma_{\star})}\E\left[f(\bar{\bx}_{K})-f_{\star}\right]\gtrsim\Lavg D^{2}$
when $K\eta n\Lavg<1$.
\end{rem}
\begin{proof}
The proof is deferred to Section \ref{sec:SGD-cvx-lb}.
\end{proof}

Theorem \ref{thm:SGD-cvx-lb} establishes a new algorithm-dependent
lower bound for $\SGD$, which, again, holds for every prespecified
$\eta>0$ and $K\in\N$ satisfying $K\eta n\Lavg\geq1$. As in the
case of $\RR$, the lower bound in Theorem \ref{thm:SGD-cvx-lb} also
differs from the upper bound in Theorem \ref{thm:SGD-cvx-ub} by factors
polynomial in $1/(2-\eta\Lmax)$ when $\eta<2/\Lmax$, and closing
this gap remains open.

The proof of Theorem \ref{thm:SGD-cvx-lb} is relatively simpler than
that of Theorem \ref{thm:RR-cvx-lb} due to the additional independence
in the sampling process. In particular, leveraging this independence,
we can analyze the iterates of $\SGD$ via a direct recursion and
evaluate it in expectation, without having to track the dependence
induced by sampling without replacement. For more details, we refer
the interested reader to Section \ref{sec:SGD-cvx-lb}.

\subsection{Worst-Case Characterization}

Combining the new upper and lower bounds, we obtain (nearly) sharp
worst-case characterizations for both $\RR$ and $\SGD$ with a constant
stepsize $\eta$ in smooth convex optimization. Concretely, given
$\Alg\in\left\{ \RR,\SGD\right\} $, $K\in\N$, $\eta\in(0,1/\Lmax]$,
and $f\in\F(d,n,D,L_{1:n},\sigma_{\star})$, where $d\in\N$, $n\geq2$,
$D>0$, $\left\{ L_{i}>0\right\} _{i=1}^{n}$, $\sigma_{\star}\geq0$,
we define
\[
\Rate_{\Alg}^{\Avg}(K,\eta,f)\defeq f(\bar{\bx}_{K})-f_{\star},
\]
where, as before, $\bar{\bx}_{K}=\frac{1}{nK}\sum_{k=1}^{K}\sum_{i=1}^{n}\bx_{k}^{i}$
denotes the average iterate generated by $\Alg$, initialized at $\bx_{1}^{1}=\bzero$,
over $K$ epochs.
\begin{rem}
For simplicity, we restrict our attention to the regime $\eta\in(0,1/\Lmax]$
to drop the term $1/(2-\eta\Lmax)$ in the upper bounds for $\RR$
(see Theorem \ref{thm:RR-cvx-ub}) and $\SGD$ (see Theorem \ref{thm:SGD-cvx-ub}).
\end{rem}
With this notation, Theorems \ref{thm:RR-cvx-ub} and \ref{thm:RR-cvx-lb}
for $\RR$, together with Theorems \ref{thm:SGD-cvx-lb} and \ref{thm:SGD-cvx-ub}
for $\SGD$, give
\begin{equation}
\sup_{f\in\F(d,n,D,L_{1:n},\sigma_{\star})}\E\left[\Rate_{\Alg}^{\Avg}(K,\eta,f)\right]\asymp\Rate_{\Alg}^{\Avg,\E}(K,\eta),\forall K\in\N,\eta\in(0,1/\Lmax],K\eta n\Lavg\geq1,\label{eq:worst-case-rate}
\end{equation}
where
\[
\Rate_{\Alg}^{\Avg,\E}(K,\eta)\defeq\begin{cases}
\frac{D^{2}}{\eta nK}+\min\left\{ 1,\eta n\Lavg\right\} \eta\sigma_{\star}^{2}, & \Alg=\RR,\\
\frac{D^{2}}{\eta nK}+\eta\sigma_{\star}^{2}, & \Alg=\SGD.
\end{cases}
\]
(\ref{eq:worst-case-rate}) immediately yields that, under the stated
conditions,
\[
\sup_{f\in\F(d,n,D,L_{1:n},\sigma_{\star})}\E\left[\Rate_{\RR}^{\Avg}(K,\eta,f)\right]\lesssim\sup_{f\in\F(d,n,D,L_{1:n},\sigma_{\star})}\E\left[\Rate_{\SGD}^{\Avg}(K,\eta,f)\right].
\]

To the best of our knowledge, Theorem 7 of \citet{pmlr-v97-haochen19a}\footnote{The theorem numbering follows the arXiv version of \citet{pmlr-v97-haochen19a}.}
is the only analogous worst-case dominance result in the literature.
Their result considers a different setting in which each $f_{i}$
is $\mu_{i}$-strongly convex and $L_{i}$-smooth, all $f_{i}$'s
share a common minimizer (equivalently, $\sigma_{\star}=0$), and
the stepsize satisfies $\eta\in(0,\min_{i\in[n]}2/(L_{i}+\mu_{i})]$.
Under these conditions, their theorem states that the worst-case squared-distance
error of $\RR$ is no larger than that of $\SGD$ after any number
of iterations.

\paragraph{Discussion of \citet{NEURIPS2020_42299f06}.}

Readers familiar with the literature of $\RR$ may notice that, at
first glance, our result appears to contradict the work of \citet{NEURIPS2020_42299f06},
which shows that $\RR$ is not uniformly better than $\SGD$. We note
that this apparent discrepancy arises from the use of different convergence
measurements. Similarly, for $\Alg\in\left\{ \RR,\SGD\right\} $,
$K\in\N$, $\eta\in(0,1/\Lmax]$, and $f\in\F(d,n,D,L_{1:n},\sigma_{\star})$,
we introduce
\[
\Rate_{\Alg}^{\Last}(K,\eta,f)\defeq f(\bx_{K}^{n+1})-f_{\star}.
\]

Under this notation, the finding of \citet{NEURIPS2020_42299f06}
should be understood as an instance-dependent result. Concretely,
Section 5 of \citet{NEURIPS2020_42299f06} constructs a function $\tilde{f}\in\F(d=3,n=4,D,\{L_{1}=(2+\sqrt{3})/2,L_{2}=1+\cos(\pi/7),L_{3}=L_{4}=2\},\sigma_{\star}=0)$,
where $D>0$ can be arbitrary, such that
\begin{eqnarray*}
\Pr\left[\Rate_{\RR}^{\Last}(K,1/\Lmax,\tilde{f})>0,\forall K\in\N\right]=1 & \text{and} & \Pr\left[\Rate_{\SGD}^{\Last}(K,1/\Lmax,\tilde{f})=0\text{ eventually}\right]=1.
\end{eqnarray*}
We emphasize two subtle points here. First, the above result is established
for the last iterate rather than the average iterate considered above.
Second, and more importantly, this almost-sure convergence dominance
does not necessarily imply, even for the constructed $\tilde{f}$,
either
\begin{eqnarray*}
\E\left[\Rate_{\SGD}^{\Last}(K,1/\Lmax,\tilde{f})\right]<\E\left[\Rate_{\RR}^{\Last}(K,1/\Lmax,\tilde{f})\right] & \text{or} & \E\left[\Rate_{\SGD}^{\Avg}(K,1/\Lmax,\tilde{f})\right]<\E\left[\Rate_{\RR}^{\Avg}(K,1/\Lmax,\tilde{f})\right].
\end{eqnarray*}
Therefore, the conclusion of \citet{NEURIPS2020_42299f06} does not
contradict ours, as the two results concern fundamentally different
settings: worst-case comparison versus instance-dependent comparison.

\subsection{Application to Functions with Quadratic Growth}

In this part, we provide an application of our new Theorem \ref{thm:RR-cvx-ub}
to functions satisfying the quadratic growth condition in Assumption
\ref{assu:quadratic-growth} below.
\begin{assumption}[Quadratic growth]
\label{assu:quadratic-growth}$f$ satisfies the $\mu$-quadratic
growth condition, i.e., $f(\bx)-f_{\star}\geq\frac{\mu}{2}\dist^{2}(\bx,\X_{\star}),\forall\bx\in\R^{d}$
for some $\mu>0$, where $\X_{\star}\defeq\arg\min_{\bx\in\R^{d}}f(\bx)$
and $\dist(\bx,\X_{\star})\defeq\inf_{\bx_{\star}\in\X_{\star}}\left\Vert \bx-\bx_{\star}\right\Vert $.
\end{assumption}
The quadratic growth condition has a long history in the optimization
literature \citep{doi:10.1137/0804001,doi:10.1287/moor.20.4.801,doi:10.1137/S1052623499359178,necoara2019linear}.
It is a strict relaxation of strong convexity, i.e., every strongly
convex function satisfies the quadratic growth condition, whereas
the converse does not hold. It is also closely related to several
other well-known conditions, including the Polyak-\L{}ojasiewicz
inequality \citep{POLYAK1963864}, error bounds \citep{luo1993error},
and the restricted secant inequality \citep{zhang2013gradient}. We
refer the interested reader to, e.g., \citet{10.1007/978-3-319-46128-1_50}
and \citet{pmlr-v130-guille-escuret21a} for systematic studies of
these relationships.

\begin{algorithm}[H]
\caption{\label{alg:Restarted-Shuffling-SGD}Restarted Shuffling Stochastic
Gradient Descent ($\protect\RShufflingSGD$)}

\textbf{Input:} initial point $\bx_{1}\in\R^{d}$, stepsize $\eta>0$,
number of epochs $K\in\N$, number of restarts $T\in\N$

\textbf{for} $t=1$ \textbf{to} $T$ \textbf{do}

$\quad$$\bx_{t+1}=\ShufflingSGD(\bx_{t},\eta,K)$, where $\bx_{t+1}$
is constructed based on (\ref{eq:avg-x})

\textbf{end for}
\end{algorithm}

To handle Assumption \ref{assu:quadratic-growth}, inspired by the
classical restart technique in the literature, e.g., \citet{nesterov1983method},
we consider a restarted variant of $\ShufflingSGD$, termed $\RShufflingSGD$,
as presented in Algorithm \ref{alg:Restarted-Shuffling-SGD}. Specifically,
we use $\ShufflingSGD$ as a subroutine and restart it after several
epochs. At each restart, the output of the previous stage is used
as the initialization for the next run. Theorem \ref{thm:RR-QG-ub}
below provides a general convergence guarantee for this restarted
scheme.
\begin{thm}
\label{thm:RR-QG-ub}Under Assumptions \ref{assu:minimizer}, \ref{assu:convexity},
\ref{assu:smoothness}, and \ref{assu:quadratic-growth}, suppose
$\RR$ is employed by $\ShufflingSGD$ with $\eta<\frac{2}{\Lmax}$
and $K\asymp K_{\star}\defeq\frac{1}{(2-\eta\Lmax)\eta n\mu}$, then
$\RShufflingSGD$ (Algorithm \ref{alg:Restarted-Shuffling-SGD}) guarantees
that
\[
\E\left[f(\bx_{T+1})-f_{\star}\right]\lesssim\frac{\mu\dist^{2}(\bx_{1},\X_{\star})}{2^{T}}+\frac{\min\left\{ 1,\eta n\Lavg\right\} \eta\sigma_{\star}^{2}}{(2-\eta\Lmax)^{2}}.
\]
\end{thm}
\begin{proof}
The proof is deferred to Subsection \ref{subsec:RR-QG-ub-Proof}.
\end{proof}

Theorem \ref{thm:RR-QG-ub} shows that, with a constant stepsize $\eta$
and a suitable choice of the number of epochs $K$ between restarts,
$\RShufflingSGD$ with $\RR$ reduces the function value gap at a
linear rate up to an additive stochastic error. In particular, when
$\sigma_{\star}=0$, i.e., the interpolation regime, the function
value gap $\E\left[f(\bx_{T+1})-f_{\star}\right]$ converges linearly
to zero, and consequently so does $\E\left[\dist^{2}(\bx_{T+1},\X_{\star})\right]$.

At a high level, the intuition behind Theorem \ref{thm:RR-QG-ub}
is as follows. At each restart stage, Theorem \ref{thm:RR-cvx-ub}
allows us to upper bound $\E\left[f(\bx_{t+1})-f_{\star}\right]$
in terms of $\E\left[\dist^{2}(\bx_{t},\X_{\star})\right]$ up to
an additive error. The quadratic growth condition, in turn, lower
bounds $f(\bx_{t+1})-f_{\star}$ by $\dist^{2}(\bx_{t+1},\X_{\star})$.
With an appropriate choice of the restart period, these two inequalities
yield a constant-factor contraction of $\E\left[\dist^{2}(\bx_{t},\X_{\star})\right]$
at each stage. Repeating this argument leads to the geometric decay
stated in Theorem \ref{thm:RR-QG-ub}, while keeping the stochastic
error at the same order as in Theorem \ref{thm:RR-cvx-ub}.
\begin{cor}
\label{cor:RR-QG-complexity}Under the same setting as in Theorem
\ref{thm:RR-QG-ub}, given $0<\varepsilon<\mu\dist^{2}(\bx_{1},\X_{\star})$
sufficiently small, with $\eta=\eta_{\star}\defeq\min\left\{ \frac{1}{\Lmax},\max\left\{ \frac{\varepsilon}{\sigma_{\star}^{2}},\frac{1}{\sigma_{\star}}\sqrt{\frac{\varepsilon}{n\Lavg}}\right\} \right\} $
and $T=T_{\star}\defeq\left\lceil \log_{2}\frac{\mu\dist^{2}(\bx_{1},\X_{\star})}{\varepsilon}\right\rceil $,
$\RShufflingSGD$ (Algorithm \ref{alg:Restarted-Shuffling-SGD}) guarantees
that
\[
\E\left[f(\bx_{T_{\star}+1})-f_{\star}\right]\lesssim\varepsilon,
\]
using the following number of individual gradient evaluations:
\[
\left(\frac{\Lmax}{\mu}+\min\left\{ \frac{\sigma_{\star}^{2}}{\mu\varepsilon},\frac{\sigma_{\star}}{\mu}\sqrt{\frac{n\Lavg}{\varepsilon}}\right\} \right)\log\frac{\mu\dist^{2}(\bx_{1},\X_{\star})}{\varepsilon}.
\]
\end{cor}
\begin{proof}
The proof is deferred to Subsection \ref{subsec:RR-QG-complexity-Proof}.
\end{proof}

With Theorem \ref{thm:RR-QG-ub} in hand, we further establish Corollary
\ref{cor:RR-QG-complexity}, which provides the number of individual
gradient evaluations required to reach an $\varepsilon$-accurate
solution. To the best of our knowledge, when each $f_{i}$ is smooth
and convex, our sample complexity bound is the sharpest currently
available among shuffling-based methods, improving upon the best known
results obtained under assumptions strictly stronger than the quadratic
growth condition: $\left(\frac{n\sqrt{\Lavg\Lmax}}{\mu}A+\frac{\sigma_{\star}}{\mu}\sqrt{\frac{n\Lavg}{\varepsilon}}B\right)\log\left(\frac{\mu\dist^{2}(\bx_{1},\X_{\star})}{\varepsilon}\right)$
for $\mu$-strongly convex $f$ \citep{pmlr-v235-liu24cg}, where
$A\asymp\sqrt{\log\left(\frac{\sqrt{\Lavg\Lmax}}{\mu}\log\left(\frac{\mu\dist^{2}(\bx_{1},\X_{\star})}{\varepsilon}\right)\right)}$
and $B\asymp\sqrt{\log\left(\frac{\sigma_{\star}}{\mu}\sqrt{\frac{\Lavg}{n\varepsilon}}\log\left(\frac{\mu\dist^{2}(\bx_{1},\X_{\star})}{\varepsilon}\right)\right)}$,
and $\left(\frac{\Lmax}{\mu}+\frac{\sigma_{\star}}{\mu}\sqrt{\frac{n\Lmax}{\varepsilon}}\right)\log\frac{\mu\dist^{2}(\bx_{1},\X_{\star})}{\varepsilon}$
for $\mu$-strongly convex $f_{i}$ \citep{NEURIPS2020_c8cc6e90},
both of which are achieved by $\ShufflingSGD$ with $\RR$.

We would like to emphasize that, although our bound has a better order,
it is achieved through a restarted scheme, namely $\RShufflingSGD$.
In contrast, the results of \citet{pmlr-v235-liu24cg} and \citet{NEURIPS2020_c8cc6e90}
mentioned above do not require restarts. Therefore, an interesting
direction for future work is to determine whether $\ShufflingSGD$
without restarts can achieve a similar guarantee for functions satisfying
the quadratic growth condition.

For $\SGD$, as far as we know, the best sample complexity bound in
the literature that applies to our setting is $\frac{\Lmax}{\mu}\log\frac{\Lavg(f(\bx_{1})-f_{\star})}{\mu\varepsilon}+\frac{\Lmax\sigma_{\star}^{2}}{\mu^{2}\varepsilon}$,
as implied by Corollary 3 of \citet{khaled2022better}. Since $\mu\dist^{2}(\bx_{1},\X_{\star})\lesssim f(\bx_{1})-f_{\star}$
and $\mu\leq\Lavg$, their first term $\frac{\Lmax}{\mu}\log\frac{\Lavg(f(\bx_{1})-f_{\star})}{\mu\varepsilon}$
is hence worse than the corresponding term $\frac{\Lmax}{\mu}\log\frac{\mu\dist^{2}(\bx_{1},\X_{\star})}{\varepsilon}$
in our bound. Although their second term $\frac{\Lmax\sigma_{\star}^{2}}{\mu^{2}\varepsilon}$
is not directly comparable in all parameter regimes with our $\min\left\{ \frac{\sigma_{\star}^{2}}{\mu\varepsilon},\frac{\sigma_{\star}}{\mu}\sqrt{\frac{n\Lavg}{\varepsilon}}\right\} \log\frac{\mu\dist^{2}(\bx_{1},\X_{\star})}{\varepsilon}$,
our bound is better both in the low-accuracy regime when $\log\frac{\mu\dist^{2}(\bx_{1},\X_{\star})}{\varepsilon}\leq\frac{\Lmax}{\mu}\Leftrightarrow\frac{\mu\dist^{2}(\bx_{1},\X_{\star})}{\exp(\Lmax/\mu)}\leq\varepsilon$
and asymptotically in the high-accuracy regime as $\varepsilon\to0$.
Taken together, these comparisons indicate that our Corollary \ref{cor:RR-QG-complexity}
improves on the best known bound for $\SGD$ under the quadratic growth
condition in different regimes.

\section{Theoretical Analysis: Upper Bound for $\protect\RR$}

In this section, we lay the groundwork for proving Theorem \ref{thm:RR-cvx-ub}
and complete its proof at the end. The section is organized into four
parts. First, we provide the three most important lemmas in Subsection
\ref{subsec:Core-Lemmas}. Next, in Subsection \ref{subsec:RR-cvx-ub-I},
we establish an upper bound in Theorem \ref{thm:RR-cvx-ub-I}, which
indicates that $\ShufflingSGD$ under $\RR$ never converges more
slowly than $\SGD$ under any reasonable stepsize. Then, Theorem \ref{thm:RR-cvx-ub-II}
in Subsection \ref{subsec:RR-cvx-ub-II} presents an alternative convergence
rate for $\ShufflingSGD$ under $\RR$, which demonstrates that Algorithm
\ref{alg:Shuffling-SGD} under $\RR$ provably converges faster than
$\SGD$ when the stepsize is sufficiently small. Finally, in Subsection
\ref{subsec:RR-cvx-ub-Proof}, we conclude Theorem \ref{thm:RR-cvx-ub}
and then use it to prove Corollary \ref{cor:RR-cvx-optimal}.

\subsection{Three Core Lemmas\label{subsec:Core-Lemmas}}

This subsection contains three core lemmas, all of which are critical
to our analysis. 

Before presenting the lemmas, we introduce two notions. Given a permutation
$\pi$ of $\left[n\right]$ and two indices $i,j\in\left[n\right]$
satisfying $j\leq i$, we define the following new permutation
\begin{equation}
\pi(i,j)\defeq\left(\pi^{1},\mydots,\pi^{j-1},\pi^{i},\pi^{j+1},\mydots,\pi^{i-1},\pi^{j},\pi^{i+1},\mydots,\pi^{n}\right).\label{eq:main-virtual-pi}
\end{equation}
In words, $\pi(i,j)$ is the permutation generated by exchanging the
elements $\pi^{i}$ and $\pi^{j}$ in $\pi$. Equipped with the notion
of $\pi(i,j)$, we introduce the following virtual sequence, for any
given $k\in\left[K\right]$,
\begin{eqnarray}
 & \bx_{k}^{l+1}(i,j)\defeq\bx_{k}^{l}(i,j)-\eta_{k}\nabla f_{\pi_{k}^{l}(i,j)}(\bx_{k}^{l}(i,j)),\forall l\in\left[n\right], & \text{where}\enspace\bx_{k}^{1}(i,j)\defeq\bx_{k}^{1}.\label{eq:main-virtual-x}
\end{eqnarray}
This means that the sequence $\bx_{k}^{l}(i,j),\forall l\in\left[n+1\right]$
denotes the trajectory of the $k$-th epoch starting from $\bx_{k}^{1}$,
produced by Algorithm \ref{alg:Shuffling-SGD}, but under the permutation
$\pi_{k}(i,j)$. This virtual iterate can be viewed as a coupled sequence
of the real output and, to the best of our knowledge, was first introduced
by \citet{NEURIPS2021_107030ca}. It plays a fundamental role in our
proof, as will become clear.

With these two concepts in hand, we proceed to state the three core
lemmas. The first is Lemma \ref{lem:core-function}, which is based
on Lemma 2 of \citet{NEURIPS2021_107030ca}. For completeness, we
reproduce the proof of Lemma \ref{lem:core-function} in Appendix
\ref{sec:Core-Lemmas}.
\begin{lem}
\label{lem:core-function}Given an arbitrary finite-sum function $\ell(\bx)=\frac{1}{n}\sum_{i=1}^{n}\ell_{i}(\bx)$,
suppose $\RR$ is employed, then for any $k\in\left[K\right]$ and
$i\in\left[n\right]$, $\ShufflingSGD$ (Algorithm \ref{alg:Shuffling-SGD})
guarantees that
\[
\E\left[\ell(\bx_{k}^{i})-\ell_{\pi_{k}^{i}}(\bx_{k}^{i})\right]=\frac{1}{n}\sum_{j=1}^{i-1}\E\left[\ell_{\pi_{k}^{i}}(\bx_{k}^{i}(i,j))-\ell_{\pi_{k}^{i}}(\bx_{k}^{i})\right],
\]
where $\bx_{k}(i,j)$ is defined in (\ref{eq:main-virtual-x}).
\end{lem}
In the analysis of shuffling gradient methods, a well-known major
challenge is to properly bound the term $\E\left[f(\bx_{k}^{i})-f_{\pi_{k}^{i}}(\bx_{k}^{i})\right]$,
unlike in $\SGD$, which no longer equals $0$ due to the nature of
shuffling-based algorithms. Lemma \ref{lem:core-function} provides
a possible approach by relating the term we want to control (in a
slightly more general form, applicable to any finite-sum function
$\ell$) to another quantity involving the virtual sequence introduced
earlier in (\ref{eq:main-virtual-x}).

For the convenience of the discussion, temporarily assume $\ell_{i}=f_{i}$
in Lemma \ref{lem:core-function}. Then, under the smoothness assumption,
one would expect the difference between $\ell_{\pi_{k}^{i}}(\bx_{k}^{i}(i,j))$
and $\ell_{\pi_{k}^{i}}(\bx_{k}^{i})$ to be small whenever $\bx_{k}^{i}(i,j)$
and $\bx_{k}^{i}$ are close. This observation naturally leads us
to another core Lemma \ref{lem:core-stability} stated below.
\begin{lem}
\label{lem:core-stability}Under Assumptions \ref{assu:convexity}
and \ref{assu:smoothness}, suppose $\eta_{k}\leq\frac{2}{\Lmax},\forall k\in\left[K\right]$,
then for any $k\in\left[K\right]$, $i\in\left[n\right]$, and $j\in\left[i-1\right]$,
$\ShufflingSGD$ (Algorithm \ref{alg:Shuffling-SGD}) guarantees that
\[
\left\Vert \bx_{k}^{i}(i,j)-\bx_{k}^{i}\right\Vert \leq\eta_{k}\left\Vert \nabla f_{\pi_{k}^{i}}(\bx_{k}^{j})-\nabla f_{\pi_{k}^{j}}(\bx_{k}^{j})\right\Vert ,
\]
where $\bx_{k}(i,j)$ is defined in (\ref{eq:main-virtual-x}).
\end{lem}
Lemma \ref{lem:core-stability} quantifies how close the virtual iterate
and the true trajectory can be, under the widely required condition
of $\eta_{k}\leq2/\Lmax$ in smooth optimization. Although the inequality
does not directly offer a bound on the distance between $\bx_{k}^{i}(i,j)$
and $\bx_{k}^{i}$ that depends only on deterministic terms (e.g.,
the stepsize $\eta_{k}$), it is sufficient for our proof when combined
with a careful analysis.

The proof of Lemma \ref{lem:core-stability} is given in Appendix
\ref{sec:Core-Lemmas} and relies on the co-coercivity property (i.e.,
Lemma \ref{lem:co-coercivity-cvx}), which is closely related to the
nonexpansiveness of the update rule in gradient methods \citep{Solr-KOHA-OAI-TEST:19722,nesterov2018lectures}.

Equipped with Lemmas \ref{lem:core-function} and \ref{lem:core-stability},
we are ready to present the final core result, Lemma \ref{lem:core-Bregman}.
\begin{lem}
\label{lem:core-Bregman}Under Assumptions \ref{assu:minimizer},
\ref{assu:convexity}, and \ref{assu:smoothness}, suppose $\RR$
is employed with $\left\{ \eta_{k}\leq\frac{2}{\Lmax}\right\} _{k=1}^{K}$,
then for any $k\in\left[K\right]$ and $i\in\left[n\right]$, $\ShufflingSGD$
(Algorithm \ref{alg:Shuffling-SGD}) guarantees that
\[
\E\left[\Breg(\bx_{k}^{i},\bx_{\star})\right]\leq2\E\left[\Breg_{\pi_{k}^{i}}(\bx_{k}^{i},\bx_{\star})\right]+\frac{64}{n}\sum_{j=1}^{i-1}\E\left[\Breg_{\pi_{k}^{j}}(\bx_{k}^{j},\bx_{\star})\right]+\frac{12(i-1)}{n}\eta_{k}^{2}\Lmax\sigma_{\star}^{2}.
\]
\end{lem}
\begin{rem}
We make no effort to optimize the constants in the bounds obtained
in this work.
\end{rem}
To the best of our knowledge, we are the first to obtain an inequality
in such a form. Intuitively, Lemma \ref{lem:core-Bregman} says that
each $\E\left[\Breg(\bx_{k}^{i},\bx_{\star})\right]$ (which now equals
$\E\left[f(\bx_{k}^{i})-f_{\star}\right]$ because $\nabla f(\bx_{\star})=\bzero$)
differs from its stochastic counterpart $\E\left[\Breg_{\pi_{k}^{i}}(\bx_{k}^{i},\bx_{\star})\right]$
by at most a multiplicative constant, together with an average of
the preceding terms $\E\left[\Breg_{\pi_{k}^{j}}(\bx_{k}^{j},\bx_{\star})\right]$
satisfying $j<i$ (not an exact average due to the coefficient $64/n$),
plus an additional term involving $\sigma_{\star}^{2}$. In particular,
summing the inequality in Lemma \ref{lem:core-Bregman} from $i=1$
to $n$ yields a meaningful lower bound on $\E\left[\sum_{i=1}^{n}\Breg_{\pi_{k}^{i}}(\bx_{k}^{i},\bx_{\star})\right]$
based on $\E\left[\sum_{i=1}^{n}\Breg(\bx_{k}^{i},\bx_{\star})\right]$
and $\sigma_{\star}^{2}$.

The proof of Lemma \ref{lem:core-Bregman} is rather technical, so
we skip the discussion here. For details, we kindly refer the interested
reader to Appendix \ref{sec:Core-Lemmas}.

\subsection{Bound I: Never Worse than $\protect\SGD$ under Reasonable Stepsize\label{subsec:RR-cvx-ub-I}}

In this subsection, we give the first bound for $\ShufflingSGD$ under
$\RR$, stated in Theorem \ref{thm:RR-cvx-ub-I} below.
\begin{thm}
\label{thm:RR-cvx-ub-I}Under Assumptions \ref{assu:minimizer}, \ref{assu:convexity},
and \ref{assu:smoothness}, suppose $\RR$ is employed with $\left\{ \eta_{k}<\frac{2}{\Lmax}\right\} _{k=1}^{K}$,
then $\ShufflingSGD$ (Algorithm \ref{alg:Shuffling-SGD}) guarantees
that
\[
\E\left[f(\bar{\bx}_{K})-f_{\star}\right]\lesssim\frac{\left\Vert \bx_{1}^{1}-\bx_{\star}\right\Vert ^{2}}{n\sum_{k=1}^{K}(2-\eta_{k}\Lmax)\eta_{k}}+\frac{\sum_{k=1}^{K}\eta_{k}^{2}/(2-\eta_{k}\Lmax)}{\sum_{k=1}^{K}(2-\eta_{k}\Lmax)\eta_{k}}\sigma_{\star}^{2},
\]
where $\bar{\bx}_{K}$ is defined in (\ref{eq:avg-x}).
\end{thm}

\paragraph{Discussion on Theorem \ref{thm:RR-cvx-ub-I}.}

To the best of our knowledge, Theorem \ref{thm:RR-cvx-ub-I} offers
the first theoretical evidence that $\ShufflingSGD$ under $\RR$
shares surprising similarities with $\SGD$, as reflected in the two
aspects elaborated below.

First, Theorem \ref{thm:RR-cvx-ub-I} states that, similar to $\SGD$,
$\ShufflingSGD$ under $\RR$ does converge under any reasonable stepsize
(i.e., $\eta_{k}<2/\Lmax$). In contrast, as far as we know, all prior
works that provide provable rates of multi-epoch $\RR$ for smooth
convex optimization require the stepsize $\eta_{k}$ to be smaller
than a threshold proportional to $1/n$, with only one exception \citep{pmlr-v97-nagaraj19a},
which, however, assumes that each $f_{i}$ is additionally Lipschitz,
thereby limiting the applicability of their theory and even excluding
common quadratic optimization problems over $\R^{d}$.

Second, we highlight that Theorem \ref{thm:RR-cvx-ub-I} gives the
same convergence upper bound (up to constant factors) as $\SGD$ \citep{garrigos2023handbook}
in smooth convex optimization, while allowing a stepsize that depends
on the epoch number. This result thus fills a gap in the literature.

Putting these together, Theorem \ref{thm:RR-cvx-ub-I} indicates that,
for smooth convex optimization, $\RR$ under any reasonable stepsize
never converges more slowly than $\SGD$.

\paragraph{Analysis.}

In the following, we present the analysis for Theorem \ref{thm:RR-cvx-ub-I}
and finally prove it. The core idea underlying the proof is, as one
might expect, to analyze $\ShufflingSGD$ in a manner analogous to
$\SGD$. In other words, we aim to quantify the progress made by Algorithm
\ref{alg:Shuffling-SGD} at each iteration. Although this perspective
is natural, it has been less explored in prior studies. The main reason
is that, as discussed earlier, the permutation in $\ShufflingSGD$
causes the most important property of $\SGD$, unbiasedness, to no
longer hold. To overcome this barrier, we develop a novel analysis
that avoids any additional assumptions, such as Lipschitz continuity
considered in \citet{pmlr-v97-nagaraj19a}.

We start with the following Lemma \ref{lem:RR-cvx-ub-I-descent},
which characterizes the per-iterate progress of $\ShufflingSGD$ (or
$\SGD$). Lemma \ref{lem:RR-cvx-ub-I-descent} is inspired by Theorem
8 of \citet{NEURIPS2025_97335d45}, which was established only for
$\SGD$. In contrast, our result slightly generalizes and tightens
their bound. Indeed, if $\pi_{k}^{i}$ is sampled uniformly from $\left[n\right]$,
after taking expectations on both sides and dropping the nonnegative
term $\Breg_{\pi_{k}^{i}}(\bx_{\star},\bx_{k}^{i})$, Lemma \ref{lem:RR-cvx-ub-I-descent}
recovers their result. We include the proof of Lemma \ref{lem:RR-cvx-ub-I-descent}
in Appendix \ref{sec:RR-cvx-ub-I}.
\begin{lem}
\label{lem:RR-cvx-ub-I-descent}Under Assumptions \ref{assu:minimizer},
\ref{assu:convexity}, and \ref{assu:smoothness}, suppose $\eta_{k}<\frac{2}{\Lmax},\forall k\in\left[K\right]$,
then for any $k\in\left[K\right]$ and $i\in\left[n\right]$, $\ShufflingSGD$
(Algorithm \ref{alg:Shuffling-SGD}) guarantees that
\begin{align*}
\frac{2-\eta_{k}\Lmax}{2}\eta_{k}\Breg_{\pi_{k}^{i}}(\bx_{k}^{i},\bx_{\star})\leq & \left\Vert \bx_{k}^{i}-\bx_{\star}\right\Vert ^{2}-\left\Vert \bx_{k}^{i+1}-\bx_{\star}\right\Vert ^{2}+\frac{2+\eta_{k}\Lmax}{2-\eta_{k}\Lmax}\eta_{k}^{2}\left\Vert \nabla f_{\pi_{k}^{i}}(\bx_{\star})\right\Vert ^{2}\\
 & -2\eta_{k}\left\langle \nabla f_{\pi_{k}^{i}}(\bx_{\star}),\bx_{k}^{i}-\bx_{\star}\right\rangle -\frac{2-\eta_{k}\Lmax}{2}\eta_{k}\Breg_{\pi_{k}^{i}}(\bx_{\star},\bx_{k}^{i}).
\end{align*}
\end{lem}
In $\SGD$, the uniform sampling scheme ensures that $\E\left[\Breg_{\pi_{k}^{i}}(\bx_{k}^{i},\bx_{\star})\right]=\E\left[\Breg(\bx_{k}^{i},\bx_{\star})\right]=\E\left[f(\bx_{k}^{i})-f_{\star}\right]$
and $\E\left[\left\langle \nabla f_{\pi_{k}^{i}}(\bx_{\star}),\bx_{k}^{i}-\bx_{\star}\right\rangle \right]=\E\left[\left\langle \nabla f(\bx_{\star}),\bx_{k}^{i}-\bx_{\star}\right\rangle \right]=0$.
In contrast, the analysis of $\ShufflingSGD$ presents two challenges.
The first is to relate $\E\left[\Breg_{\pi_{k}^{i}}(\bx_{k}^{i},\bx_{\star})\right]$
to $\E\left[\Breg(\bx_{k}^{i},\bx_{\star})\right]$, which can be
handled by the previous Lemma \ref{lem:core-Bregman}. The second,
and the main remaining difficulty, is to derive a suitable upper bound
on $\E\left[-\left\langle \nabla f_{\pi_{k}^{i}}(\bx_{\star}),\bx_{k}^{i}-\bx_{\star}\right\rangle \right]$.

We address this latter challenge by establishing the following novel
inequality in Lemma \ref{lem:RR-cvx-ub-I-core}.
\begin{lem}
\label{lem:RR-cvx-ub-I-core}Under Assumptions \ref{assu:minimizer},
\ref{assu:convexity}, and \ref{assu:smoothness}, suppose $\RR$
is employed with $\left\{ \eta_{k}<\frac{2}{\Lmax}\right\} _{k=1}^{K}$,
then for any $k\in\left[K\right]$ and $i\in\left[n\right]$, $\ShufflingSGD$
(Algorithm \ref{alg:Shuffling-SGD}) guarantees that
\[
\E\left[-\left\langle \nabla f_{\pi_{k}^{i}}(\bx_{\star}),\bx_{k}^{i}-\bx_{\star}\right\rangle \right]\leq\frac{2-\eta_{k}\Lmax}{4n}\sum_{j=1}^{i-1}\E\left[\Breg_{\pi_{k}^{j}}(\bx_{\star},\bx_{k}^{j})\right]+\frac{2(i-1)}{n}\cdot\frac{2+3\eta_{k}\Lmax}{2-\eta_{k}\Lmax}\eta_{k}\sigma_{\star}^{2}.
\]
\end{lem}
Lemma \ref{lem:RR-cvx-ub-I-core} provides a novel inequality to
upper bound $\E\left[-\left\langle \nabla f_{\pi_{k}^{i}}(\bx_{\star}),\bx_{k}^{i}-\bx_{\star}\right\rangle \right]$
by the second moment of the stochastic gradients up to time $i$.
Note that the second term on the R.H.S. can be absorbed by the R.H.S.
of the inequality in Lemma \ref{lem:RR-cvx-ub-I-descent}. For the
remaining term, the coefficient $\eta_{k}/n$ is key to the final
proof, which ensures that the accumulated error in one epoch is controlled
by $\E[\eta_{k}\sum_{i=1}^{n}\Vert\nabla f_{\pi_{k}^{i}}(\bx_{k}^{i})\Vert^{2}]$.

The proof of Lemma \ref{lem:RR-cvx-ub-I-core} builds on the two core
results, Lemmas \ref{lem:core-function} and \ref{lem:core-stability},
presented before. To save space, we defer it to Appendix \ref{sec:RR-cvx-ub-I}.

\paragraph{Final proof.}

With Lemmas \ref{lem:RR-cvx-ub-I-descent} and \ref{lem:RR-cvx-ub-I-core}
stated above, we are ready to prove Theorem \ref{thm:RR-cvx-ub-I}.

\begin{proof}[Proof of Theorem \ref{thm:RR-cvx-ub-I}]
We sum the inequality in Lemma \ref{lem:RR-cvx-ub-I-descent} from
$i=1$ to $n$ and use $\bx_{k+1}^{1}=\bx_{k}^{n+1}$ to obtain that
\begin{align*}
\frac{2-\eta_{k}\Lmax}{2}\eta_{k}\sum_{i=1}^{n}\Breg_{\pi_{k}^{i}}(\bx_{k}^{i},\bx_{\star})\leq & \left\Vert \bx_{k}^{1}-\bx_{\star}\right\Vert ^{2}-\left\Vert \bx_{k+1}^{1}-\bx_{\star}\right\Vert ^{2}+\frac{2+\eta_{k}\Lmax}{2-\eta_{k}\Lmax}\eta_{k}^{2}\sum_{i=1}^{n}\left\Vert \nabla f_{\pi_{k}^{i}}(\bx_{\star})\right\Vert ^{2}\\
 & -2\eta_{k}\sum_{i=1}^{n}\left\langle \nabla f_{\pi_{k}^{i}}(\bx_{\star}),\bx_{k}^{i}-\bx_{\star}\right\rangle -\frac{2-\eta_{k}\Lmax}{2}\eta_{k}\sum_{i=1}^{n}\Breg_{\pi_{k}^{i}}(\bx_{\star},\bx_{k}^{i}).
\end{align*}
Take expectations on both sides and note that $\E\left[\left\Vert \nabla f_{\pi_{k}^{i}}(\bx_{\star})\right\Vert ^{2}\right]=\sigma_{\star}^{2}$
to yield that
\begin{align}
\frac{2-\eta_{k}\Lmax}{2}\eta_{k}\E\left[\sum_{i=1}^{n}\Breg_{\pi_{k}^{i}}(\bx_{k}^{i},\bx_{\star})\right]\leq & \E\left[\left\Vert \bx_{k}^{1}-\bx_{\star}\right\Vert ^{2}\right]-\E\left[\left\Vert \bx_{k+1}^{1}-\bx_{\star}\right\Vert ^{2}\right]+\frac{2+\eta_{k}\Lmax}{2-\eta_{k}\Lmax}\eta_{k}^{2}n\sigma_{\star}^{2}\nonumber \\
 & +2\eta_{k}\sum_{i=1}^{n}\E\left[-\left\langle \nabla f_{\pi_{k}^{i}}(\bx_{\star}),\bx_{k}^{i}-\bx_{\star}\right\rangle \right]-\frac{2-\eta_{k}\Lmax}{2}\eta_{k}\sum_{i=1}^{n}\E\left[\Breg_{\pi_{k}^{i}}(\bx_{\star},\bx_{k}^{i})\right].\label{eq:I-rate-1}
\end{align}

Next, we invoke Lemma \ref{lem:RR-cvx-ub-I-core} and sum it up from
$i=1$ to $n$ to have
\begin{equation}
\sum_{i=1}^{n}\E\left[-\left\langle \nabla f_{\pi_{k}^{i}}(\bx_{\star}),\bx_{k}^{i}-\bx_{\star}\right\rangle \right]\leq\frac{2-\eta_{k}\Lmax}{4}\sum_{i=1}^{n}\E\left[\Breg_{\pi_{k}^{i}}(\bx_{\star},\bx_{k}^{i})\right]+\frac{2+3\eta_{k}\Lmax}{2-\eta_{k}\Lmax}\eta_{k}n\sigma_{\star}^{2}.\label{eq:I-rate-2}
\end{equation}
Combine (\ref{eq:I-rate-1}) and (\ref{eq:I-rate-2}) to obtain that
\begin{equation}
\frac{2-\eta_{k}\Lmax}{2}\eta_{k}\E\left[\sum_{i=1}^{n}\Breg_{\pi_{k}^{i}}(\bx_{k}^{i},\bx_{\star})\right]\leq\E\left[\left\Vert \bx_{k}^{1}-\bx_{\star}\right\Vert ^{2}\right]-\E\left[\left\Vert \bx_{k+1}^{1}-\bx_{\star}\right\Vert ^{2}\right]+\frac{6+7\eta_{k}\Lmax}{2-\eta_{k}\Lmax}\eta_{k}^{2}n\sigma_{\star}^{2}.\label{eq:I-rate-3}
\end{equation}

One more step, we observe that, by Lemma \ref{lem:core-Bregman},
\begin{equation}
\E\left[\sum_{i=1}^{n}\Breg(\bx_{k}^{i},\bx_{\star})\right]\leq66\E\left[\sum_{i=1}^{n}\Breg_{\pi_{k}^{i}}(\bx_{k}^{i},\bx_{\star})\right]+6\eta_{k}^{2}n\Lmax\sigma_{\star}^{2}\leq66\E\left[\sum_{i=1}^{n}\Breg_{\pi_{k}^{i}}(\bx_{k}^{i},\bx_{\star})\right]+12\eta_{k}n\sigma_{\star}^{2},\label{eq:I-rate-4}
\end{equation}
where the last step is due to $\eta_{k}\leq\frac{2}{\Lmax}$.

Finally, we plug (\ref{eq:I-rate-4}) back into (\ref{eq:I-rate-3}),
use $\eta_{k}\Lmax\leq2$, sum over $k=1$ to $K$, divide both sides
by $n\sum_{k=1}^{K}(2-\eta_{k}\Lmax)\eta_{k}$, apply the convexity
of $f$, and use the definition of $\bar{\bx}_{K}$ in (\ref{eq:avg-x})
to conclude.
\end{proof}

\subsection{Bound II: Always Better than $\protect\SGD$ under Small Stepsize\label{subsec:RR-cvx-ub-II}}

This subsection presents the other rate of $\ShufflingSGD$ under
$\RR$, as shown in Theorem \ref{thm:RR-cvx-ub-II} below.
\begin{thm}
\label{thm:RR-cvx-ub-II}Under Assumptions \ref{assu:minimizer},
\ref{assu:convexity}, and \ref{assu:smoothness}, suppose $\RR$
is employed with $\left\{ \eta_{k}<\frac{2}{\Lmax}\right\} _{k=1}^{K}$,
then $\ShufflingSGD$ (Algorithm \ref{alg:Shuffling-SGD}) guarantees
that
\[
\E\left[f(\bar{\bx}_{K})-f_{\star}\right]\lesssim\frac{\left\Vert \bx_{1}^{1}-\bx_{\star}\right\Vert ^{2}}{n\sum_{k=1}^{K}(2-\eta_{k}\Lmax)\eta_{k}}+\frac{\sum_{k=1}^{K}\eta_{k}^{3}n\Lavg/(2-\eta_{k}\Lmax)}{\sum_{k=1}^{K}(2-\eta_{k}\Lmax)\eta_{k}}\sigma_{\star}^{2},
\]
where $\bar{\bx}_{K}$ is defined in (\ref{eq:avg-x}).
\end{thm}

\paragraph{Discussion on Theorem \ref{thm:RR-cvx-ub-II}.}

Readers familiar with the literature on shuffling gradient methods
may readily figure out that the rate given in Theorem \ref{thm:RR-cvx-ub-II}
perfectly matches the known bound for $\ShufflingSGD$ under $\RR$
in smooth convex optimization (e.g., \citet{NEURIPS2020_c8cc6e90,JMLR:v22:20-1238}).
However, we emphasize a key difference here, that is, the stepsize
in our Theorem \ref{thm:RR-cvx-ub-II} is allowed to satisfy $\eta_{k}<2/\Lmax$,
in contrast to all existing results that require $\eta_{k}$ to be
at most inversely proportional to $n$.

More importantly, in the setting of nonuniform $L_{i}$ considered
in this work, the largest threshold on the stepsize in the literature
that guarantees a rate similar to Theorem \ref{thm:RR-cvx-ub-II}
is in the order of $1/(n\sqrt{\Lavg\Lmax})$ \citep{pmlr-v235-liu24cg,NEURIPS2024_84d39572}.
But as indicated by our Theorem \ref{thm:RR-cvx-ub-II}, the superiority
of $\RR$ over $\SGD$ already exists once $\eta_{k}\lesssim1/(n\Lavg)$.
Especially, this improvement can be significant when a dominant smoothness
parameter exists, leading to $\Lmax\approx n\Lavg$.

Therefore, Theorem \ref{thm:RR-cvx-ub-II} is the first result to
extend the known bound in smooth convex optimization to any reasonable
stepsize while preserving the favorable property of Algorithm \ref{alg:Shuffling-SGD},
i.e., $\ShufflingSGD$ under $\RR$ provably converges faster than
$\SGD$ when the stepsize is sufficiently small.

\paragraph{Analysis.}

The roadmap for establishing Theorem \ref{thm:RR-cvx-ub-II} differs
substantially from that used before to prove Theorem \ref{thm:RR-cvx-ub-I}.
This time, our proof strategy is to check how close $\ShufflingSGD$
can be to Gradient Descent. More concretely, we will view each epoch
of Algorithm \ref{alg:Shuffling-SGD} (containing $n$ iterations)
as a single step and analyze the progress made by it at once. This
kind of approach has appeared in different previous works (e.g., \citet{NEURIPS2020_c8cc6e90,JMLR:v22:20-1238})
and always yields a convergence rate in the order of $1/(n^{1/3}K^{2/3})$,
faster than $\SGD$ when $K$ is large.

However, all works that follow the idea described above share the
same issue, that is, they require the stepsize to be in the order
of $1/n$, which is, however, not enough for our purpose. In the following,
we develop a new analysis to bypass this critical obstacle.

We now formally begin the analysis by introducing another virtual
sequence, defined as follows, for any $k\in\left[K\right]$,
\begin{eqnarray}
\by_{k}^{i+1}\defeq\by_{k}^{i}-\eta_{k}\nabla f_{\pi_{k}^{i}}(\bx_{\star}),\forall i\in\left[n\right], & \by_{k+1}^{1}\defeq\by_{k}^{n+1}, & \text{where}\enspace\by_{1}^{1}\defeq\bx_{\star}.\label{eq:main-virtual-y-def}
\end{eqnarray}
Under the above definition, and noting that $\nabla f(\bx_{\star})=\bzero$,
one can find
\[
\by_{k}^{n+1}=\by_{k}^{1}-\eta_{k}\sum_{i=1}^{n}\nabla f_{\pi_{k}^{i}}(\bx_{\star})=\by_{k}^{1}-\eta_{k}n\nabla f(\bx_{\star})=\by_{k}^{1},\forall k\in\left[K\right].
\]
Combine the above line and $\by_{1}^{1}=\bx_{\star}$ as defined in
(\ref{eq:main-virtual-y-def}) to have
\begin{equation}
\by_{k}^{n+1}=\by_{k}^{1}=\bx_{\star},\forall k\in\left[K\right].\label{eq:main-virtual-y-prop}
\end{equation}
The above virtual sequence is inspired by the work of \citet{NEURIPS2020_c8cc6e90},
which, as far as we know, was the first to propose a similar term
under the constant stepsize. Here, we slightly extend their idea to
accommodate the case where the stepsize can depend on the current
epoch number.
\begin{rem}
We note that \citet{NEURIPS2020_c8cc6e90} introduced the virtual
sequence to handle the case of individual strong convexity, i.e.,
each $f_{i}$ is required to be strongly convex. However, in our
setting, only individual convexity is assumed. This means that their
proof cannot be applied. As such, our analysis substantially departs
from the existing approach.
\end{rem}
Equipped with the new virtual sequence introduced above, we first
present the following Lemma \ref{lem:RR-cvx-ub-II-descent}.
\begin{lem}
\label{lem:RR-cvx-ub-II-descent}Under Assumptions \ref{assu:minimizer},
\ref{assu:convexity}, and \ref{assu:smoothness}, suppose $\eta_{k}<\frac{2}{\Lmax},\forall k\in\left[K\right]$,
then for any $k\in\left[K\right]$, $\ShufflingSGD$ (Algorithm \ref{alg:Shuffling-SGD})
guarantees that
\[
\frac{4-2\eta_{k}\Lmax}{2+\eta_{k}\Lmax}\eta_{k}\sum_{i=1}^{n}\Breg_{\pi_{k}^{i}}(\bx_{k}^{i},\bx_{\star})\leq\left\Vert \bx_{k}^{1}-\bx_{\star}\right\Vert ^{2}-\left\Vert \bx_{k+1}^{1}-\bx_{\star}\right\Vert ^{2}+\frac{4+\eta_{k}^{2}\Lmax^{2}}{4-\eta_{k}^{2}\Lmax^{2}}\eta_{k}\sum_{i=1}^{n}L_{\pi_{k}^{i}}\left\Vert \by_{k}^{i}-\bx_{\star}\right\Vert ^{2},
\]
where $\by_{k}^{i}$ is defined in (\ref{eq:main-virtual-y-def}).
\end{lem}
As discussed earlier, we intentionally treat each epoch of $\ShufflingSGD$
as a single step. Hence, compared with Lemma \ref{lem:RR-cvx-ub-I-descent}
used to prove the first bound, Lemma \ref{lem:RR-cvx-ub-II-descent}
has a different flavor, which shows the progress made by Algorithm
\ref{alg:Shuffling-SGD} over an entire epoch. Based on the form
of Lemma \ref{lem:RR-cvx-ub-II-descent}, two tasks naturally arise.
The first is to lower bound the term $\sum_{i=1}^{n}\Breg_{\pi_{k}^{i}}(\bx_{k}^{i},\bx_{\star})$
on the L.H.S. by $\sum_{i=1}^{n}\Breg(\bx_{k}^{i},\bx_{\star})=\sum_{i=1}^{n}f(\bx_{k}^{i})-f_{\star}$.
The second is to upper bound the residual term $\sum_{i=1}^{n}L_{\pi_{k}^{i}}\left\Vert \by_{k}^{i}-\bx_{\star}\right\Vert ^{2}$
on the R.H.S.

Fortunately, the first task can be directly addressed by summing the
inequality in Lemma \ref{lem:core-Bregman} over $i\in\left[n\right]$.
Therefore, we only need to focus on the residual term $\sum_{i=1}^{n}L_{\pi_{k}^{i}}\left\Vert \by_{k}^{i}-\bx_{\star}\right\Vert ^{2}$
in Lemma \ref{lem:RR-cvx-ub-II-descent}, which is bounded in Lemma
\ref{lem:RR-cvx-ub-II-residual} below.
\begin{lem}
\label{lem:RR-cvx-ub-II-residual}Under Assumptions \ref{assu:minimizer}
and \ref{assu:smoothness}, suppose $\RR$ is employed, then for any
$k\in\left[K\right]$, $\ShufflingSGD$ (Algorithm \ref{alg:Shuffling-SGD})
guarantees that
\[
\E\left[\sum_{i=1}^{n}L_{\pi_{k}^{i}}\left\Vert \by_{k}^{i}-\bx_{\star}\right\Vert ^{2}\right]\leq\frac{\eta_{k}^{2}n^{2}\Lavg\sigma_{\star}^{2}}{3},
\]
where $\by_{k}^{i}$ is defined in (\ref{eq:main-virtual-y-def}).
\end{lem}
Lemma \ref{lem:RR-cvx-ub-II-residual} can be derived in a relatively
easy way, as it can be deduced from existing works. In particular,
thanks to \citet{pmlr-v235-liu24cg}, we obtain a bound that depends
only on the average smoothness parameter $\Lavg$ rather than the
maximum smoothness parameter $\Lmax$.

\paragraph{Final proof.}

Armed with Lemmas \ref{lem:RR-cvx-ub-II-descent} and \ref{lem:RR-cvx-ub-II-residual}
above, we are finally able to prove Theorem \ref{thm:RR-cvx-ub-II}.

\begin{proof}[Proof of Theorem \ref{thm:RR-cvx-ub-II}]
First, we invoke Lemma \ref{lem:RR-cvx-ub-II-descent} to have
\[
\frac{4-2\eta_{k}\Lmax}{2+\eta_{k}\Lmax}\eta_{k}\sum_{i=1}^{n}\Breg_{\pi_{k}^{i}}(\bx_{k}^{i},\bx_{\star})\leq\left\Vert \bx_{k}^{1}-\bx_{\star}\right\Vert ^{2}-\left\Vert \bx_{k+1}^{1}-\bx_{\star}\right\Vert ^{2}+\frac{4+\eta_{k}^{2}\Lmax^{2}}{4-\eta_{k}^{2}\Lmax^{2}}\eta_{k}\sum_{i=1}^{n}L_{\pi_{k}^{i}}\left\Vert \by_{k}^{i}-\bx_{\star}\right\Vert ^{2},
\]
where $\by_{k}^{i}$ is defined in (\ref{eq:main-virtual-y-def}).
Take expectations on both sides and apply Lemma \ref{lem:RR-cvx-ub-II-residual}
to yield
\begin{equation}
\frac{4-2\eta_{k}\Lmax}{2+\eta_{k}\Lmax}\eta_{k}\E\left[\sum_{i=1}^{n}\Breg_{\pi_{k}^{i}}(\bx_{k}^{i},\bx_{\star})\right]\leq\E\left[\left\Vert \bx_{k}^{1}-\bx_{\star}\right\Vert ^{2}\right]-\E\left[\left\Vert \bx_{k+1}^{1}-\bx_{\star}\right\Vert ^{2}\right]+\frac{4+\eta_{k}^{2}\Lmax^{2}}{3(4-\eta_{k}^{2}\Lmax^{2})}\eta_{k}^{3}n^{2}\Lavg\sigma_{\star}^{2}.\label{eq:II-rate-1}
\end{equation}

Next, we multiply both sides of the inequality in Lemma \ref{lem:core-Bregman}
by $\frac{4-2\eta_{k}\Lmax}{2+\eta_{k}\Lmax}\eta_{k}$ and sum it
up from $i=1$ to $n$ to obtain
\begin{align}
\frac{4-2\eta_{k}\Lmax}{2+\eta_{k}\Lmax}\eta_{k}\E\left[\sum_{i=1}^{n}\Breg(\bx_{k}^{i},\bx_{\star})\right] & \leq66\frac{4-2\eta_{k}\Lmax}{2+\eta_{k}\Lmax}\eta_{k}\E\left[\sum_{i=1}^{n}\Breg_{\pi_{k}^{i}}(\bx_{k}^{i},\bx_{\star})\right]+\frac{24-12\eta_{k}\Lmax}{2+\eta_{k}\Lmax}\eta_{k}^{3}n\Lmax\sigma_{\star}^{2}\nonumber \\
 & \leq66\frac{4-2\eta_{k}\Lmax}{2+\eta_{k}\Lmax}\eta_{k}\E\left[\sum_{i=1}^{n}\Breg_{\pi_{k}^{i}}(\bx_{k}^{i},\bx_{\star})\right]+\frac{24-12\eta_{k}\Lmax}{2+\eta_{k}\Lmax}\eta_{k}^{3}n^{2}\Lavg\sigma_{\star}^{2},\label{eq:II-rate-2}
\end{align}
where the last step is due to $\Lmax\leq n\Lavg$. Combine (\ref{eq:II-rate-1})
and (\ref{eq:II-rate-2}) to have
\[
\frac{4-2\eta_{k}\Lmax}{2+\eta_{k}\Lmax}\eta_{k}\E\left[\sum_{i=1}^{n}\Breg(\bx_{k}^{i},\bx_{\star})\right]\leq66\E\left[\left\Vert \bx_{k}^{1}-\bx_{\star}\right\Vert ^{2}\right]-66\E\left[\left\Vert \bx_{k+1}^{1}-\bx_{\star}\right\Vert ^{2}\right]+\frac{12(2-\eta_{k}\Lmax)^{2}+22(4+\eta_{k}^{2}\Lmax^{2})}{4-\eta_{k}^{2}\Lmax^{2}}\eta_{k}^{3}n^{2}\Lavg\sigma_{\star}^{2},
\]
which implies that, by $\eta_{k}\Lmax\leq2$,
\begin{align*}
\E\left[\sum_{i=1}^{n}(2-\eta_{k}\Lmax)\eta_{k}\Breg(\bx_{k}^{i},\bx_{\star})\right] & \leq132\E\left[\left\Vert \bx_{k}^{1}-\bx_{\star}\right\Vert ^{2}\right]-132\E\left[\left\Vert \bx_{k+1}^{1}-\bx_{\star}\right\Vert ^{2}\right]+\frac{112\eta_{k}^{3}n^{2}\Lavg}{2-\eta_{k}\Lmax}\sigma_{\star}^{2}\\
\Rightarrow\E\left[\sum_{k=1}^{K}\sum_{i=1}^{n}(2-\eta_{k}\Lmax)\eta_{k}\Breg(\bx_{k}^{i},\bx_{\star})\right] & \leq132\left\Vert \bx_{1}^{1}-\bx_{\star}\right\Vert ^{2}+112\sum_{k=1}^{K}\frac{\eta_{k}^{3}n^{2}\Lavg}{2-\eta_{k}\Lmax}\sigma_{\star}^{2}.
\end{align*}

Finally, we observe that $\Breg(\bx_{k}^{i},\bx_{\star})=f(\bx_{k}^{i})-f_{\star}$,
divide both sides by $n\sum_{k=1}^{K}(2-\eta_{k}\Lmax)\eta_{k}$,
apply the convexity of $f$, and use the definition of $\bar{\bx}_{K}$
in (\ref{eq:avg-x}) to conclude.
\end{proof}

\subsection{Proofs of Theorem \ref{thm:RR-cvx-ub} and Corollary \ref{cor:RR-cvx-optimal}\label{subsec:RR-cvx-ub-Proof}}

With the previous preparation, Theorem \ref{thm:RR-cvx-ub} follows
immediately.

\begin{proof}[Proof of Theorem \ref{thm:RR-cvx-ub}]
Combine Theorems \ref{thm:RR-cvx-ub-I} and \ref{thm:RR-cvx-ub-II}
to conclude.
\end{proof}

We next derive Corollary \ref{cor:RR-cvx-optimal} directly from Theorem
\ref{thm:RR-cvx-ub}.

\begin{proof}[Proof of Corollary \ref{cor:RR-cvx-optimal}]
With a constant stepsize $\eta_{k}=\eta<\frac{2}{\Lmax},\forall k\in\left[K\right]$,
Theorem \ref{thm:RR-cvx-ub} reduces to
\[
\E\left[f(\bar{\bx}_{K})-f_{\star}\right]\lesssim\frac{\left\Vert \bx_{1}^{1}-\bx_{\star}\right\Vert ^{2}}{(2-\eta\Lmax)\eta nK}+\frac{\min\left\{ 1,\eta n\Lavg\right\} \eta\sigma_{\star}^{2}}{(2-\eta\Lmax)^{2}}.
\]
Optimizing the R.H.S. of the above inequality over $0<\eta<\frac{2}{\Lmax}$
yields the desired result.
\end{proof}

\section{Theoretical Analysis: Lower Bound for $\protect\RR$\label{sec:RR-cvx-lb}}

In this section, our goal is to prove Theorem \ref{thm:RR-cvx-lb}.
Our strategy is again to separate the analysis into two complementary
parts. In Subsection \ref{subsec:RR-cvx-lb-I}, we study the error
arising from $D$, the initial distance between $\bx_{1}^{1}$ and
$\bx_{\star}$, and derive the corresponding lower bound. In Subsection
\ref{subsec:RR-cvx-lb-II}, we turn to the stochastic error induced
by $\sigma_{\star}^{2}$, the variance of the gradient at $\bx_{\star}$,
and establish the second lower bound. Finally, in Subsection \ref{subsec:RR-cvx-lb-Proof},
we combine these two lower bounds to complete the proof of Theorem
\ref{thm:RR-cvx-lb}.

\subsection{Bound I: Induced by the Initial Distance\label{subsec:RR-cvx-lb-I}}

In this subsection, we prove the first lower bound in Theorem \ref{thm:RR-cvx-lb-I},
showing that the error induced by the initial distance $D$ cannot
decay faster than $\frac{D^{2}}{\eta nK}$ under the condition $K\eta n\Lavg\geq1$.
Importantly, this lower bound applies to any shuffling strategy, rather
than being specific to $\RR$.

To prove the result, we construct a simple quadratic hard instance
whose smoothness parameters for the individual components are chosen
according to the prescribed stepsize $\eta$ and parameters $L_{1},\mydots,L_{n}$.
The resulting dynamics can be tracked explicitly, allowing us to analyze
the iterates under an arbitrary ordering of the component functions
and derive the claimed lower bound.
\begin{thm}
\label{thm:RR-cvx-lb-I}Given $d\in\N$, $n\geq2$, $D>0$, $\left\{ L_{i}>0\right\} _{i=1}^{n}$,
and $\sigma_{\star}\geq0$, suppose any shuffling strategy is employed
with $\left\{ \eta_{k}=\eta\right\} _{k=1}^{K}$ for $\eta>0$, let
$\bar{\bx}_{K}\defeq\frac{1}{nK}\sum_{k=1}^{K}\sum_{i=1}^{n}\bx_{k}^{i}$
denote the average iterate generated by $\ShufflingSGD$ (Algorithm
\ref{alg:Shuffling-SGD}), initialized at $\bx_{1}^{1}=\bzero$, over
$K\in\N$ epochs, then we have
\begin{eqnarray*}
\sup_{f\in\F(d,n,D,L_{1:n},\sigma_{\star})}\E\left[f(\bar{\bx}_{K})-f_{\star}\right]\gtrsim\frac{D^{2}}{\eta nK} & \text{when} & K\eta n\Lavg\geq1.
\end{eqnarray*}
\end{thm}
\begin{rem}
In fact, the lower bound is stronger than stated: for the constructed
hard instance, the lower bound holds for every sample path, rather
than merely in expectation. See the proof for details.
\end{rem}
\begin{proof}
In the following proof, $\bv$ denotes an arbitrary unit vector in
$\R^{d}$. We first let
\[
\beta_{i}\defeq\min\left\{ \frac{1}{2\eta K},n\Lavg\right\} \frac{L_{i}}{n\Lavg},\forall i\in\left[n\right].
\]
Clearly, each $\beta_{i}$ satisfies $0<\beta_{i}\leq\min\left\{ \frac{1}{\eta},L_{i}\right\} $.
Moreover, we write
\begin{equation}
\bar{\beta}\defeq\frac{1}{n}\sum_{i=1}^{n}\beta_{i}=\min\left\{ \frac{1}{2\eta nK},\Lavg\right\} .\label{eq:RR-cvx-lb-I-beta}
\end{equation}
Next, we define
\[
f_{i}(\bx)\defeq\frac{\beta_{i}\left\Vert \bx-D\bv\right\Vert ^{2}}{2},\forall i\in\left[n\right].
\]
The corresponding finite-sum objective is
\begin{equation}
f(\bx)\defeq\frac{1}{n}\sum_{i=1}^{n}f_{i}(\bx)=\frac{\bar{\beta}\left\Vert \bx-D\bv\right\Vert ^{2}}{2}.\label{eq:RR-cvx-lb-I-f}
\end{equation}
Note that $\bx_{\star}=D\bv\Rightarrow\left\Vert \bx_{\star}\right\Vert =D$,
each $f_{i}$ is differentiable, convex, and $L_{i}$-smooth (since
$\beta_{i}\leq L_{i}$), and $\frac{1}{n}\sum_{i=1}^{n}\left\Vert \nabla f_{i}(\bx_{\star})\right\Vert ^{2}=0\leq\sigma_{\star}^{2}$.
Therefore, $f$ belongs to the class $\F(d,n,D,L_{1:n},\sigma_{\star})$.

Now, we can find that
\[
\nabla f_{\pi_{k}^{i}}(\bx_{k}^{i})=\beta_{\pi_{k}^{i}}(\bx_{k}^{i}-D\bv),\forall k\in\left[K\right],i\in\left[n\right],
\]
which implies that
\[
\bx_{k}^{i+1}=(1-\eta\beta_{\pi_{k}^{i}})\bx_{k}^{i}+\eta\beta_{\pi_{k}^{i}}D\bv,\forall k\in\left[K\right],i\in\left[n\right].
\]
Since $\bx_{1}^{1}=\bzero$, unrolling the recursion yields
\[
\bx_{k}^{i}=\left(1-\rho^{k-1}\prod_{j=1}^{i-1}(1-\eta\beta_{\pi_{k}^{j}})\right)D\bv,\forall k\in\left[K\right],i\in\left[n\right],
\]
where
\begin{equation}
\rho\defeq\prod_{i=1}^{n}(1-\eta\beta_{i})\overset{\eta\beta_{i}\in\left[0,1\right],\forall i\in\left[n\right]}{\geq}1-\eta\sum_{i=1}^{n}\beta_{i}=1-\eta n\bar{\beta}\overset{(\ref{eq:RR-cvx-lb-I-beta})}{\geq}1-\frac{1}{2K}\geq\frac{1}{2}.\label{eq:RR-cvx-lb-I-rho}
\end{equation}
As such, we have
\begin{equation}
\bar{\bx}_{K}=\frac{1}{nK}\sum_{k=1}^{K}\sum_{i=1}^{n}\bx_{k}^{i}=\left(1-\frac{1}{nK}\sum_{k=1}^{K}\rho^{k-1}\sum_{i=1}^{n}\prod_{j=1}^{i-1}(1-\eta\beta_{\pi_{k}^{j}})\right)D\bv.\label{eq:RR-cvx-lb-I-x}
\end{equation}

Finally, we know
\begin{align*}
f(\bar{\bx}_{K})-f_{\star} & \overset{(\ref{eq:RR-cvx-lb-I-f}),(\ref{eq:RR-cvx-lb-I-x}),\left\Vert \bv\right\Vert =1}{=}\frac{\bar{\beta}\left(\sum_{k=1}^{K}\rho^{k-1}\sum_{i=1}^{n}\prod_{j=1}^{i-1}(1-\eta\beta_{\pi_{k}^{j}})\right)^{2}D^{2}}{2n^{2}K^{2}}\overset{(a)}{\geq}\frac{(1-\exp(-K\eta n\bar{\beta}))^{2}D^{2}}{8\eta^{2}n^{2}K^{2}\bar{\beta}}\\
 & \overset{(b)}{\geq}\frac{\bar{\beta}D^{2}}{32}\overset{(\ref{eq:RR-cvx-lb-I-beta})}{=}\frac{D^{2}}{32\eta nK}\min\left\{ \frac{1}{2},K\eta n\Lavg\right\} \geq\frac{D^{2}}{64\eta nK}\min\left\{ 1,K\eta n\Lavg\right\} ,
\end{align*}
where $(a)$ is due to
\[
\prod_{j=1}^{i-1}(1-\eta\beta_{\pi_{k}^{j}})\geq\prod_{j=1}^{n}(1-\eta\beta_{j})\overset{(\ref{eq:RR-cvx-lb-I-rho})}{=}\rho,\forall i\in\left[n\right],
\]
which implies that
\[
\sum_{k=1}^{K}\rho^{k-1}\sum_{i=1}^{n}\prod_{j=1}^{i-1}(1-\eta\beta_{\pi_{k}^{j}})\geq n\rho\sum_{k=1}^{K}\rho^{k-1}\overset{(\ref{eq:RR-cvx-lb-I-rho})}{\geq}\frac{n}{2}\sum_{k=1}^{K}(1-\eta n\bar{\beta})^{k-1}=\frac{1-(1-\eta n\bar{\beta})^{K}}{2\eta\bar{\beta}}\geq\frac{1-\exp(-K\eta n\bar{\beta})}{2\eta\bar{\beta}},
\]
and $(b)$ follows from $1-\exp(-x)\geq\frac{x}{2},\forall x\in\left[0,\frac{1}{2}\right]$
and $0\leq K\eta n\bar{\beta}\overset{(\ref{eq:RR-cvx-lb-I-beta})}{\leq}\frac{1}{2}$.
Therefore, for the constructed $f\in\F(d,n,D,L_{1:n},\sigma_{\star})$,
any shuffling strategy satisfies
\[
f(\bar{\bx}_{K})-f_{\star}\gtrsim\frac{D^{2}}{\eta nK}\min\left\{ 1,K\eta n\Lavg\right\} .
\]
Consequently,
\begin{eqnarray*}
\sup_{f\in\F(d,n,D,L_{1:n},\sigma_{\star})}\E\left[f(\bar{\bx}_{K})-f_{\star}\right]\gtrsim\frac{D^{2}}{\eta nK} & \text{when} & K\eta n\Lavg\geq1.
\end{eqnarray*}
\end{proof}

\subsection{Bound II: Induced by the Gradient Variance at the Optimum\label{subsec:RR-cvx-lb-II}}

This subsection establishes the second lower bound in Theorem \ref{thm:RR-cvx-lb-II},
which captures the stochastic error $\sigma_{\star}^{2}$ induced
by the variance of the gradient at $\bx_{\star}$. In contrast to
the first lower bound in Theorem \ref{thm:RR-cvx-lb-I}, which holds
for any shuffling strategy, Theorem \ref{thm:RR-cvx-lb-II} is specific
to $\RR$ and exploits the randomness of the permutation within each
epoch. 
\begin{thm}
\label{thm:RR-cvx-lb-II}Given $d\in\N$, $n\geq2$, $D>0$, $\left\{ L_{i}>0\right\} _{i=1}^{n}$,
and $\sigma_{\star}\geq0$, suppose $\RR$ is employed with $\left\{ \eta_{k}=\eta\right\} _{k=1}^{K}$
for $\eta>0$, let $\bar{\bx}_{K}\defeq\frac{1}{nK}\sum_{k=1}^{K}\sum_{i=1}^{n}\bx_{k}^{i}$
denote the average iterate generated by $\ShufflingSGD$ (Algorithm
\ref{alg:Shuffling-SGD}), initialized at $\bx_{1}^{1}=\bzero$, over
$K\in\N$ epochs, then we have
\begin{eqnarray*}
\sup_{f\in\F(d,n,D,L_{1:n},\sigma_{\star})}\E\left[f(\bar{\bx}_{K})-f_{\star}\right]\gtrsim\min\left\{ 1,\eta n\Lavg\right\} \eta\sigma_{\star}^{2} & \text{when} & K\eta n\Lavg\geq1.
\end{eqnarray*}
\end{thm}
The analysis required for Theorem \ref{thm:RR-cvx-lb-II} is considerably
more involved. We therefore treat the cases $n\geq3$ and $n=2$ separately,
establishing the corresponding results in Theorems \ref{thm:RR-cvx-lb-II-n>2}
and \ref{thm:RR-cvx-lb-II-n=00003D2} below, respectively. 

At a high level, the proof is based on carefully designed hard instances
whose dynamics under $\RR$ can be analyzed explicitly. For $n\geq3$,
a key ingredient is the asymmetric piecewise quadratic structure,
allowing the randomness of the permutation to generate a persistent
drift across epochs. Together with the nonuniform smoothness parameters
of each component, this structure makes the effect of the random ordering
analytically tractable. In contrast, the case $n=2$ admits a simpler
direct analysis. In both cases, the key is to quantify how the stochastic
effect introduced within each epoch accumulates under $\RR$ and remains
nonnegligible in the average iterate, ultimately yielding the claimed
lower bound.

\paragraph{Analysis when $n\protect\geq3$.}

The key part of the proof is to design $n$ asymmetric piecewise quadratic
functions inspired by \citet{pmlr-v119-rajput20a,yun2022minibatch,pmlr-v202-cha23a},
in which sampling without replacement generates a persistent bias
even when the algorithm starts at $\bx_{\star}$. Controlling this
bias across epochs then yields the desired term $\min\left\{ 1,\eta n\Lavg\right\} \eta\sigma_{\star}^{2}$.

The most important technical component of our proof builds on Lemma
10 of \citet{pmlr-v119-rajput20a}, Lemma 13 of \citet{yun2022minibatch},
and Lemma B.4 of \citet{pmlr-v202-cha23a}. Concretely, given $k\in\left[K\right]$,
we couple the true iterates $\left\{ \bx_{k}^{i}\right\} _{i=1}^{n+1}$
generated by the asymmetric piecewise quadratic hard instance $\left\{ f_{i}\right\} _{i=1}^{n}$
(see (\ref{eq:RR-cvx-lb-II-n>2-f_i})) with $\left\{ \bz_{k}^{i}\right\} _{i=1}^{n+1}$
generated by $\left\{ h_{i}\right\} _{i=1}^{n}$ (see (\ref{eq:RR-cvx-lb-II-n>2-h_i})),
the symmetric counterpart of $\left\{ f_{i}\right\} _{i=1}^{n}$,
and then apply a novel comparison lemma (see Lemma \ref{lem:RR-cvx-lb-II-auxiliary})
to derive a recursion in expectation for the true iterates. However,
the existing analyses consider only uniform smoothness parameters,
i.e., $L_{i}=L,\forall i\in\left[n\right]$. Therefore, our proof
requires several modifications to accommodate the nonuniform case.
Moreover, for our purpose, we must explicitly track the dependence
on the stepsize $\eta$.
\begin{thm}
\label{thm:RR-cvx-lb-II-n>2}Given $d\in\N$, $n\geq3$, $D>0$, $\left\{ L_{i}>0\right\} _{i=1}^{n}$,
and $\sigma_{\star}\geq0$, suppose $\RR$ is employed with $\left\{ \eta_{k}=\eta\right\} _{k=1}^{K}$
for $\eta>0$, let $\bar{\bx}_{K}\defeq\frac{1}{nK}\sum_{k=1}^{K}\sum_{i=1}^{n}\bx_{k}^{i}$
denote the average iterate generated by $\ShufflingSGD$ (Algorithm
\ref{alg:Shuffling-SGD}), initialized at $\bx_{1}^{1}=\bzero$, over
$K\in\N$ epochs, then we have
\begin{eqnarray*}
\sup_{f\in\F(d,n,D,L_{1:n},\sigma_{\star})}\E\left[f(\bar{\bx}_{K})-f_{\star}\right]\gtrsim\min\left\{ 1,\eta n\Lavg\right\} \eta\sigma_{\star}^{2} & \text{when} & K\eta n\Lavg\geq1.
\end{eqnarray*}
\end{thm}
\begin{proof}
W.l.o.g., we assume $L_{1}\geq L_{2}\geq\mydots\geq L_{n-1}\geq L_{n}$
in the following proof. We first define $S\defeq\sum_{i=1}^{n-2}L_{i}$
and note that
\begin{equation}
3S\geq S+L_{n-1}+L_{n}=n\Lavg\Rightarrow S\geq\frac{n\Lavg}{3}.\label{eq:RR-cvx-lb-II-n>2-S}
\end{equation}
We now introduce
\begin{equation}
\beta_{i}\defeq\begin{cases}
\min\left\{ \frac{1}{2\eta},S\right\} \frac{L_{i}}{S} & i\in\left[n-2\right]\\
0 & i\in\left\{ n-1,n\right\} 
\end{cases}.\label{eq:RR-cvx-lb-II-n>2-beta_i}
\end{equation}
By construction, each $\beta_{i}$ satisfies $0\leq\beta_{i}\leq\min\left\{ \frac{1}{\eta},L_{i}\right\} $
and at least one $\beta_{i}$ is positive since $n\geq3$. Moreover,
we write
\begin{equation}
\bar{\beta}\defeq\frac{1}{n}\sum_{i=1}^{n}\beta_{i}=\min\left\{ \frac{1}{2\eta n},\frac{S}{n}\right\} \overset{(\ref{eq:RR-cvx-lb-II-n>2-S})}{\geq}\min\left\{ \frac{1}{2\eta n},\frac{\Lavg}{3}\right\} \geq\frac{1}{3}\min\left\{ \frac{1}{\eta n},\Lavg\right\} .\label{eq:RR-cvx-lb-II-n>2-beta}
\end{equation}
We then define
\begin{equation}
f_{i}(\bx)\defeq\left(\beta_{i}\1\left[\bx[1]<0\right]+\frac{\beta_{i}}{2}\1\left[\bx[1]\geq0\right]\right)\frac{\left(\bx[1]\right)^{2}}{2}+\chi_{i}\sqrt{\frac{n}{2}}\sigma_{\star}\bx[1],\forall i\in\left[n\right],\label{eq:RR-cvx-lb-II-n>2-f_i}
\end{equation}
where $\chi_{i}\defeq\begin{cases}
0 & i\in\left[n-2\right]\\
-1 & i=n-1\\
1 & i=n
\end{cases}$. The corresponding finite-sum objective is
\begin{equation}
f(\bx)\defeq\frac{1}{n}\sum_{i=1}^{n}f_{i}(\bx)=\left(\bar{\beta}\1\left[\bx[1]<0\right]+\frac{\bar{\beta}}{2}\1\left[\bx[1]\geq0\right]\right)\frac{\left(\bx[1]\right)^{2}}{2}\overset{(\ref{eq:RR-cvx-lb-II-n>2-beta})}{\geq}\frac{\min\left\{ \frac{1}{\eta n},\Lavg\right\} \left(\bx[1]\right)^{2}}{12}.\label{eq:RR-cvx-lb-II-n>2-f}
\end{equation}
Note that $\bx_{\star}\defeq\bzero\in\arg\min_{\bx\in\R^{d}}f(\bx)\Rightarrow\left\Vert \bx_{\star}\right\Vert =0<D$,
each $f_{i}$ is differentiable, convex, and $L_{i}$-smooth since
$\frac{\beta_{1}}{2}\leq\beta_{i}\leq L_{i}$, and $\frac{1}{n}\sum_{i=1}^{n}\left\Vert \nabla f_{i}(\bx_{\star})\right\Vert ^{2}=\frac{1}{2}\sum_{i=1}^{n}\chi_{i}^{2}\sigma_{\star}^{2}=\sigma_{\star}^{2}$.
Therefore, $f\in\F(d,n,D,L_{1:n},\sigma_{\star})$.

Next, we introduce
\begin{equation}
h_{i}(\bx)\defeq\beta_{i}\frac{\left(\bx[1]\right)^{2}}{2}+\chi_{i}\sqrt{\frac{n}{2}}\sigma_{\star}\bx[1],\forall i\in\left[n\right].\label{eq:RR-cvx-lb-II-n>2-h_i}
\end{equation}
Moreover, given $k\in\left[K\right]$, we recursively define
\begin{eqnarray}
\bz_{k}^{i+1}\defeq\bz_{k}^{i}-\eta\nabla h_{\pi_{k}^{i}}(\bz_{k}^{i}),\forall i\in\left[n\right] & \text{where} & \bz_{k}^{1}\defeq\bx_{k}^{1}.\label{eq:RR-cvx-lb-II-n>2-z-def}
\end{eqnarray}
Furthermore, we let
\[
\Delta_{k}^{i}\defeq\bx_{k}^{i}[1]-\bz_{k}^{i}[1],\forall i\in\left[n\right].
\]
We use induction to prove that, for all $i\in\left[n+1\right]$,
\begin{equation}
\Delta_{k}^{i}\geq\frac{\eta}{2}\sum_{j=1}^{i-1}\beta_{\pi_{k}^{j}}(\bz_{k}^{j}[1])_{+}\prod_{l=j+1}^{i-1}(1-\eta\beta_{\pi_{k}^{l}}),\label{eq:RR-cvx-lb-II-n>2-hypothesis}
\end{equation}
where $(x)_{+}\defeq x\1\left[x\geq0\right]$. When $i=1$, (\ref{eq:RR-cvx-lb-II-n>2-hypothesis})
is true by the definition of $\bz_{k}^{1}=\bx_{k}^{1}$. Suppose (\ref{eq:RR-cvx-lb-II-n>2-hypothesis})
holds for some $i\in\left[n\right]$, then for the index $i+1$, we
have
\begin{align*}
\bx_{k}^{i+1}[1] & \overset{(\ref{eq:RR-cvx-lb-II-n>2-f_i})}{=}\left(1-\eta\beta_{\pi_{k}^{i}}\1\left[\bx_{k}^{i}[1]<0\right]-\eta\frac{\beta_{\pi_{k}^{i}}}{2}\1\left[\bx_{k}^{i}[1]\geq0\right]\right)\bx_{k}^{i}[1]-\eta\chi_{\pi_{k}^{i}}\sqrt{\frac{n}{2}}\sigma_{\star},\\
\bz_{k}^{i+1}[1] & \overset{(\ref{eq:RR-cvx-lb-II-n>2-h_i}),(\ref{eq:RR-cvx-lb-II-n>2-z-def})}{=}(1-\eta\beta_{\pi_{k}^{i}})\bz_{k}^{i}[1]-\eta\chi_{\pi_{k}^{i}}\sqrt{\frac{n}{2}}\sigma_{\star},
\end{align*}
which together imply that
\[
\Delta_{k}^{i+1}=(1-\eta\beta_{\pi_{k}^{i}})\Delta_{k}^{i}+\frac{\eta}{2}\beta_{\pi_{k}^{i}}(\bx_{k}^{i}[1])_{+}\overset{(a)}{\geq}(1-\eta\beta_{\pi_{k}^{i}})\Delta_{k}^{i}+\frac{\eta}{2}\beta_{\pi_{k}^{i}}(\bz_{k}^{i}[1])_{+}\overset{(b)}{\geq}\frac{\eta}{2}\sum_{j=1}^{i}\beta_{\pi_{k}^{j}}(\bz_{k}^{j}[1])_{+}\prod_{l=j+1}^{i}(1-\eta\beta_{\pi_{k}^{l}}),
\]
where $(a)$ is due to $\Delta_{k}^{i}\geq0$ (by the inductive hypothesis)
and the monotonicity of $(\cdot)_{+}$, and $(b)$ is by $1-\eta\beta_{\pi_{k}^{i}}\overset{(\ref{eq:RR-cvx-lb-II-n>2-beta_i})}{\geq}1-\frac{1}{2}>0$
and the inductive hypothesis. Therefore, we complete the induction.

By (\ref{eq:RR-cvx-lb-II-n>2-hypothesis}), we have
\begin{equation}
\Delta_{k}^{i}\geq\frac{\eta}{2}\sum_{j=1}^{i-1}\beta_{\pi_{k}^{j}}(\bz_{k}^{j}[1])_{+}\prod_{l=j+1}^{i-1}(1-\eta\beta_{\pi_{k}^{l}})\geq\frac{\rho\eta}{2}\sum_{j=1}^{i-1}\beta_{\pi_{k}^{j}}(\bz_{k}^{j}[1])_{+},\forall i\in\left[n+1\right],\label{eq:RR-cvx-lb-II-n>2-1}
\end{equation}
where 
\begin{equation}
\rho\defeq\prod_{i=1}^{n}(1-\eta\beta_{i})\overset{\eta\beta_{i}\in\left[0,1\right],\forall i\in\left[n\right]}{\geq}1-\eta\sum_{i=1}^{n}\beta_{i}=1-\eta n\bar{\beta}\overset{(\ref{eq:RR-cvx-lb-II-n>2-beta})}{=}1-\min\left\{ \frac{1}{2},\eta S\right\} \geq\frac{1}{2}.\label{eq:RR-cvx-lb-II-n>2-rho}
\end{equation}
Taking expectations on both sides of (\ref{eq:RR-cvx-lb-II-n>2-1})
yields that, for all $k\in\left[K\right]$ and $i\in\left[n+1\right]$,
\begin{equation}
\E\left[\bx_{k}^{i}[1]\right]\geq\E\left[\bz_{k}^{i}[1]\right]+\frac{\rho\eta}{2}\E\left[\sum_{j=1}^{i-1}\beta_{\pi_{k}^{j}}(\bz_{k}^{j}[1])_{+}\right]=\E\left[\prod_{j=1}^{i-1}(1-\eta\beta_{\pi_{k}^{j}})\right]\E\left[\bx_{k}^{1}[1]\right]+\frac{\rho\eta}{2}\E\left[\sum_{j=1}^{i-1}\beta_{\pi_{k}^{j}}(\bz_{k}^{j}[1])_{+}\right],\label{eq:RR-cvx-lb-II-n>2-2}
\end{equation}
where the last step is due to the first equation in Lemma \ref{lem:RR-cvx-lb-II-auxiliary}.

Finally, we consider two cases:
\begin{casenv}
\item $K\geq2$. In this case, we take $i=n+1$ in (\ref{eq:RR-cvx-lb-II-n>2-2})
to obtain, for all $k\in\left[K\right]$,
\begin{align*}
\E\left[\bx_{k+1}^{1}[1]\right] & =\E\left[\bx_{k}^{n+1}[1]\right]\geq\E\left[\prod_{j=1}^{n}(1-\eta\beta_{\pi_{k}^{j}})\right]\E\left[\bx_{k}^{1}[1]\right]+\frac{\rho\eta}{2}\E\left[\sum_{i=1}^{n}\beta_{\pi_{k}^{i}}(\bz_{k}^{i}[1])_{+}\right]\\
 & \overset{(\ref{eq:RR-cvx-lb-II-n>2-rho})}{=}\rho\E\left[\bx_{k}^{1}[1]\right]+\frac{\rho\eta}{2}\E\left[\sum_{i=1}^{n}\beta_{\pi_{k}^{i}}(\bz_{k}^{i}[1])_{+}\right]\\
 & \geq\rho\E\left[\bx_{k}^{1}[1]\right]+\frac{\rho^{2}\eta n\bar{\beta}}{12}\E\left[\left(\bx_{k}^{1}[1]+\eta\sqrt{\frac{n}{2}}\sigma_{\star}\right)_{+}\right],
\end{align*}
where the last step is by the second inequality in Lemma \ref{lem:RR-cvx-lb-II-auxiliary}.
Since $\bx_{1}^{1}=\bzero$, a simple induction gives that $\E\left[\bx_{k}^{1}[1]\right]\geq0,\forall k\in\left[K\right]$,
which further implies that, for any $k\in\left[K\right]$,
\[
\E\left[\bx_{k+1}^{1}[1]\right]\geq\rho\E\left[\bx_{k}^{1}[1]\right]+\frac{\rho^{2}\eta^{2}n^{\frac{3}{2}}\bar{\beta}\sigma_{\star}}{12\sqrt{2}}\overset{(\ref{eq:RR-cvx-lb-II-n>2-rho})}{\geq}(1-\eta n\bar{\beta})\E\left[\bx_{k}^{1}[1]\right]+\frac{\rho^{2}\eta^{2}\bar{\beta}n^{\frac{3}{2}}\sigma_{\star}}{12\sqrt{2}}.
\]
Unrolling this recursion yields
\begin{equation}
\E\left[\bx_{k}^{1}[1]\right]\geq\frac{\rho^{2}\eta\sqrt{n}\sigma_{\star}}{12\sqrt{2}}(1-(1-\eta n\bar{\beta})^{k-1})\overset{(\ref{eq:RR-cvx-lb-II-n>2-rho})}{\geq}\frac{\eta\sqrt{n}\sigma_{\star}}{48\sqrt{2}}(1-(1-\eta n\bar{\beta})^{k-1}),\forall k\in\left[K\right].\label{eq:RR-cvx-lb-II-n>2-3}
\end{equation}
Next, we invoke (\ref{eq:RR-cvx-lb-II-n>2-2}) and drop the second
nonnegative term to have
\begin{equation}
\E\left[\bx_{k}^{i}[1]\right]\geq\E\left[\prod_{j=1}^{i-1}(1-\eta\beta_{\pi_{k}^{j}})\right]\E\left[\bx_{k}^{1}[1]\right]\geq\frac{\E\left[\bx_{k}^{1}[1]\right]}{2},\forall k\in\left[K\right],i\in\left[n+1\right],\label{eq:RR-cvx-lb-II-n>2-4}
\end{equation}
where the last step is due to $\E\left[\bx_{k}^{1}[1]\right]\geq0,\forall k\in\left[K\right]$
and $\prod_{j=1}^{i-1}(1-\eta\beta_{\pi_{k}^{j}})\overset{0\leq\eta\beta_{l}\leq1,\forall l\in\left[n\right]}{\geq}\prod_{j=1}^{n}(1-\eta\beta_{\pi_{k}^{j}})\overset{(\ref{eq:RR-cvx-lb-II-n>2-rho})}{\geq}\frac{1}{2}$.
Therefore, we can find that
\begin{align}
\E\left[\bar{\bx}_{K}[1]\right] & =\frac{1}{nK}\sum_{k=1}^{K}\sum_{i=1}^{n}\E\left[\bx_{k}^{i}[1]\right]\overset{(\ref{eq:RR-cvx-lb-II-n>2-4})}{\geq}\frac{1}{2K}\sum_{k=1}^{K}\E\left[\bx_{k}^{1}[1]\right]\overset{(\ref{eq:RR-cvx-lb-II-n>2-3})}{\geq}\frac{\eta\sqrt{n}\sigma_{\star}}{96\sqrt{2}K}\sum_{k=1}^{K}1-(1-\eta n\bar{\beta})^{k-1}\nonumber \\
 & \geq\frac{\eta\sqrt{n}\sigma_{\star}}{96\sqrt{2}K}\sum_{k=1}^{K}1-\exp(-\eta n\bar{\beta}(k-1))\geq\frac{\eta\sqrt{n}\sigma_{\star}}{96\sqrt{2}K}\sum_{k=\left\lceil K/2\right\rceil +1}^{K}1-\exp(-\eta n\bar{\beta}(k-1))\nonumber \\
 & \geq\frac{\eta\sqrt{n}\sigma_{\star}}{96\sqrt{2}}\cdot\frac{\left\lfloor K/2\right\rfloor }{K}\cdot(1-\exp(-K\eta n\bar{\beta}/2))\overset{K\geq2}{\geq}\frac{\eta\sqrt{n}\sigma_{\star}}{288\sqrt{2}}(1-\exp(-K\eta n\bar{\beta}/2))\nonumber \\
 & \overset{(a)}{\geq}\frac{\min\left\{ 1-\exp(-1/2),\min\left\{ \frac{K}{2},\frac{K\eta n\Lavg}{3}\right\} /4\right\} }{288\sqrt{2}}\eta\sqrt{n}\sigma_{\star}\geq\frac{\min\left\{ 1,K\eta n\Lavg\right\} }{3456\sqrt{2}}\eta\sqrt{n}\sigma_{\star},\label{eq:RR-cvx-lb-II-n>2-5}
\end{align}
where $(a)$ follows from $1-\exp(-x)\geq\min\left\{ 1-\exp(-1/2),x/2\right\} ,\forall x\geq0$
and $K\eta n\bar{\beta}\overset{(\ref{eq:RR-cvx-lb-II-n>2-beta})}{=}\min\left\{ \frac{K}{2},\frac{K\eta n\Lavg}{3}\right\} $.
Hence, we have
\begin{align*}
\E\left[f(\bar{\bx}_{K})-f_{\star}\right] & \geq f(\E\left[\bar{\bx}_{K}\right])-f_{\star}\overset{(\ref{eq:RR-cvx-lb-II-n>2-f})}{\geq}\frac{\min\left\{ \frac{1}{\eta n},\Lavg\right\} \left(\E\left[\bar{\bx}_{K}[1]\right]\right)^{2}}{12}\\
 & \overset{(\ref{eq:RR-cvx-lb-II-n>2-5})}{\gtrsim}\min\left\{ 1,\eta n\Lavg\right\} \left(\min\left\{ 1,K\eta n\Lavg\right\} \right)^{2}\eta\sigma_{\star}^{2},
\end{align*}
where the first step is by Jensen's inequality.
\item $K=1$. In this case, we take $k=1$ in (\ref{eq:RR-cvx-lb-II-n>2-2})
and use $\bx_{1}^{1}=\bzero$ to obtain
\begin{equation}
\E\left[\bx_{1}^{i}[1]\right]\geq\frac{\rho\eta}{2}\E\left[\sum_{j=1}^{i-1}\beta_{\pi_{1}^{j}}(\bz_{1}^{j}[1])_{+}\right],\forall i\in\left[n+1\right].\label{eq:RR-cvx-lb-II-n>2-6}
\end{equation}
Therefore, we can find that
\begin{align}
\E\left[\bar{\bx}_{K}[1]\right] & =\E\left[\bar{\bx}_{1}[1]\right]=\frac{1}{n}\sum_{i=1}^{n}\E\left[\bx_{1}^{i}[1]\right]\overset{(\ref{eq:RR-cvx-lb-II-n>2-6})}{\geq}\frac{\rho\eta}{2n}\sum_{i=1}^{n}\E\left[\sum_{j=1}^{i-1}\beta_{\pi_{1}^{j}}(\bz_{1}^{j}[1])_{+}\right]\nonumber \\
 & =\frac{\rho\eta}{2n}\E\left[\sum_{i=1}^{n}(n-i)\beta_{\pi_{1}^{i}}(\bz_{1}^{i}[1])_{+}\right]\overset{(a)}{\geq}\frac{\rho^{2}\eta(n-1)\bar{\beta}}{24}\E\left[\left(\bx_{1}^{1}[1]+\eta\sqrt{\frac{n}{2}}\sigma_{\star}\right)_{+}\right]\nonumber \\
 & \overset{\bx_{1}^{1}[1]=0}{=}\frac{\rho^{2}\eta^{2}(n-1)\sqrt{n}\bar{\beta}\sigma_{\star}}{24\sqrt{2}}\overset{(\ref{eq:RR-cvx-lb-II-n>2-rho})}{\geq}\frac{\eta^{2}(n-1)\sqrt{n}\bar{\beta}\sigma_{\star}}{96\sqrt{2}}\overset{n-1\geq\frac{n}{2}}{\geq}\frac{\eta^{2}n^{\frac{3}{2}}\bar{\beta}\sigma_{\star}}{192\sqrt{2}}\nonumber \\
 & \overset{(\ref{eq:RR-cvx-lb-II-n>2-rho})}{\geq}\frac{\eta\sqrt{n}\min\left\{ 1,\eta n\Lavg\right\} \sigma_{\star}}{576\sqrt{2}}.\label{eq:RR-cvx-lb-II-n>2-7}
\end{align}
Hence, we have
\begin{align*}
\E\left[f(\bar{\bx}_{K})-f_{\star}\right] & \geq f(\E\left[\bar{\bx}_{K}\right])-f_{\star}\overset{(\ref{eq:RR-cvx-lb-II-n>2-f})}{\geq}\frac{\min\left\{ \frac{1}{\eta n},\Lavg\right\} \left(\E\left[\bar{\bx}_{K}[1]\right]\right)^{2}}{12}\\
 & \overset{(\ref{eq:RR-cvx-lb-II-n>2-7})}{\gtrsim}\left(\min\left\{ 1,\eta n\Lavg\right\} \right)^{3}\eta\sigma_{\star}^{2}\overset{K=1}{=}\min\left\{ 1,\eta n\Lavg\right\} \left(\min\left\{ 1,K\eta n\Lavg\right\} \right)^{2}\eta\sigma_{\star}^{2},
\end{align*}
where the first step is by Jensen's inequality.
\end{casenv}
Combining the two cases above, we conclude that, for $n\geq3$,
\[
\sup_{f\in\F(d,n,D,L_{1:n},\sigma_{\star})}\E\left[f(\bar{\bx}_{K})-f_{\star}\right]\gtrsim\min\left\{ 1,\eta n\Lavg\right\} \left(\min\left\{ 1,K\eta n\Lavg\right\} \right)^{2}\eta\sigma_{\star}^{2}.
\]
Consequently, for $n\geq3$,
\begin{eqnarray*}
\sup_{f\in\F(d,n,D,L_{1:n},\sigma_{\star})}\E\left[f(\bar{\bx}_{K})-f_{\star}\right]\gtrsim\min\left\{ 1,\eta n\Lavg\right\} \eta\sigma_{\star}^{2} & \text{when} & K\eta n\Lavg\geq1.
\end{eqnarray*}
\end{proof}

Lemma \ref{lem:RR-cvx-lb-II-auxiliary} collects three technical results
that characterize key properties of the auxiliary sequence $\left\{ \bz_{k}^{i}\right\} _{i=1}^{n+1}$
and serve as important ingredients for proving Theorem \ref{thm:RR-cvx-lb-II-n>2}.
In particular, these results provide the identities and lower bounds
needed to control the drift induced by $\RR$ within each epoch.
\begin{lem}
\label{lem:RR-cvx-lb-II-auxiliary}Let $\left\{ h_{i}\right\} _{i=1}^{n}$
be defined as in (\ref{eq:RR-cvx-lb-II-n>2-h_i}) and, for a fixed
$k\in\left[K\right]$, let $\left\{ \bz_{k}^{i}\right\} _{i=1}^{n+1}$
be defined as in (\ref{eq:RR-cvx-lb-II-n>2-z-def}), then we have
\begin{itemize}
\item $\E\left[\bz_{k}^{i}[1]\right]=\E\left[\prod_{j=1}^{i-1}(1-\eta\beta_{\pi_{k}^{j}})\right]\E\left[\bx_{k}^{1}[1]\right],\forall i\in\left[n+1\right]$.
\item $\E\left[\sum_{i=1}^{n}\beta_{\pi_{k}^{i}}(\bz_{k}^{i}[1])_{+}\right]\geq\frac{\rho n\bar{\beta}}{6}\E\left[\left(\bx_{k}^{1}[1]+\eta\sqrt{\frac{n}{2}}\sigma_{\star}\right)_{+}\right]$.
\item $\E\left[\sum_{i=1}^{n}(n-i)\beta_{\pi_{k}^{i}}(\bz_{k}^{i}[1])_{+}\right]\geq\frac{\rho n(n-1)\bar{\beta}}{12}\E\left[\left(\bx_{k}^{1}[1]+\eta\sqrt{\frac{n}{2}}\sigma_{\star}\right)_{+}\right]$.
\end{itemize}
\end{lem}
\begin{proof}
By definition,
\[
\bz_{k}^{i+1}[1]\overset{(\ref{eq:RR-cvx-lb-II-n>2-h_i}),(\ref{eq:RR-cvx-lb-II-n>2-z-def})}{=}(1-\eta\beta_{\pi_{k}^{i}})\bz_{k}^{i}[1]-\eta\chi_{\pi_{k}^{i}}\sqrt{\frac{n}{2}}\sigma_{\star},\forall i\in\left[n\right].
\]
Unrolling this recursion yields
\begin{equation}
\bz_{k}^{i+1}[1]=\left(\prod_{j=1}^{i}(1-\eta\beta_{\pi_{k}^{j}})\right)\bx_{k}^{1}[1]-\eta\sqrt{\frac{n}{2}}\sigma_{\star}\sum_{j=1}^{i}\chi_{\pi_{k}^{j}}\prod_{l=j+1}^{i}(1-\eta\beta_{\pi_{k}^{l}}),\forall i\in\left[n\right].\label{eq:RR-cvx-lb-II-auxiliary-z}
\end{equation}

\subparagraph{The first equation.}

For $\left\{ a,b\right\} \subset\left[n\right]$, we define the operator
$\T_{a\leftrightarrow b}:S_{n}\to S_{n},o\mapsto\T_{a\leftrightarrow b}(o)$,
where $\T_{a\leftrightarrow b}(o)$ is obtained by swapping the positions
of $a$ and $b$ in $o\in S_{n}$. Recall that $\beta_{n-1}=\beta_{n}$
(see (\ref{eq:RR-cvx-lb-II-n>2-beta_i})) and $\chi_{i}=\begin{cases}
0 & i\in\left[n-2\right]\\
-1 & i=n-1\\
1 & i=n
\end{cases}$, we therefore have, for any $1\leq j\leq i\leq n$,
\[
\chi_{\T_{n-1\leftrightarrow n}^{j}(\pi_{k})}\prod_{l=j+1}^{i}(1-\eta\beta_{\T_{n-1\leftrightarrow n}^{l}(\pi_{k})})=\chi_{\T_{n-1\leftrightarrow n}^{j}(\pi_{k})}\prod_{l=j+1}^{i}(1-\eta\beta_{\pi_{k}^{l}})=-\chi_{\pi_{k}^{j}}\prod_{l=j+1}^{i}(1-\eta\beta_{\pi_{k}^{l}}).
\]
Consequently,
\begin{equation}
\E\left[\chi_{\pi_{k}^{j}}\prod_{l=j+1}^{i}(1-\eta\beta_{\pi_{k}^{l}})\right]=-\E\left[\chi_{\T_{n-1\leftrightarrow n}^{j}(\pi_{k})}\prod_{l=j+1}^{i}(1-\eta\beta_{\T_{n-1\leftrightarrow n}^{l}(\pi_{k})})\right]=-\E\left[\chi_{\pi_{k}^{j}}\prod_{l=j+1}^{i}(1-\eta\beta_{\pi_{k}^{l}})\right],\label{eq:RR-cvx-lb-II-auxiliary-1}
\end{equation}
where the last step holds because, when $\pi_{k}$ is uniformly sampled
from $S_{n}$, $\pi_{k}$ equals $\T_{n-1\leftrightarrow n}(\pi_{k})$
in distribution. Hence, (\ref{eq:RR-cvx-lb-II-auxiliary-1}) implies
that
\[
\E\left[\chi_{\pi_{k}^{j}}\prod_{l=j+1}^{i}(1-\eta\beta_{\pi_{k}^{l}})\right]=0,\forall1\leq j\leq i\leq n.
\]
Finally, we plug the above equation back into (\ref{eq:RR-cvx-lb-II-auxiliary-z})
to have
\[
\E\left[\bz_{k}^{i+1}[1]\right]=\E\left[\left(\prod_{j=1}^{i}(1-\eta\beta_{\pi_{k}^{j}})\right)\bx_{k}^{1}[1]\right]=\E\left[\prod_{j=1}^{i}(1-\eta\beta_{\pi_{k}^{j}})\right]\E\left[\bx_{k}^{1}[1]\right],\forall i\in\left[n\right],
\]
where the last step follows from the independence of $\pi_{k}$ and
$\bx_{k}^{1}$. Relabeling the indices completes the proof.

\subparagraph{The second inequality.}

In what follows, for $i\in\left[n\right]$, let $\mathring{\pi}_{k}^{i}\in\left[n\right]$
denote the unique index such that $\pi_{k}^{\mathring{\pi}_{k}^{i}}=i$.
In other words, $\mathring{\pi}_{k}$ is the inverse permutation of
$\pi_{k}$, and hence is also uniformly distributed over $S_{n}$.
Moreover, we introduce the event $E_{k}^{i}\defeq\left\{ \mathring{\pi}_{k}^{n-1}<\mathring{\pi}_{k}^{i}<\mathring{\pi}_{k}^{n}\right\} ,\forall i\in\left[n-2\right]$.
Since $\mathring{\pi}_{k}$ is a random permutation, we have
\begin{equation}
\Pr\left[E_{k}^{i}\right]=\frac{1}{6},\forall i\in\left[n-2\right].\label{eq:RR-cvx-lb-II-auxiliary-2}
\end{equation}

Now, we observe that
\begin{equation}
\sum_{i=1}^{n}\beta_{\pi_{k}^{i}}(\bz_{k}^{i}[1])_{+}=\sum_{i=1}^{n}\beta_{i}(\bz_{k}^{\mathring{\pi}_{k}^{i}}[1])_{+}\overset{(a)}{=}\sum_{i=1}^{n-2}\beta_{i}(\bz_{k}^{\mathring{\pi}_{k}^{i}}[1])_{+}\geq\sum_{i=1}^{n-2}\beta_{i}(\bz_{k}^{\mathring{\pi}_{k}^{i}}[1])_{+}\1\left[E_{k}^{i}\right],\label{eq:RR-cvx-lb-II-auxiliary-3}
\end{equation}
where $(a)$ follows from $\beta_{n-1}=\beta_{n}=0$ (see (\ref{eq:RR-cvx-lb-II-n>2-beta_i})).
Next, given $i\in\left[n-2\right]$, we can find that
\begin{align}
(\bz_{k}^{\mathring{\pi}_{k}^{i}}[1])_{+}\1\left[E_{k}^{i}\right] & \overset{(\ref{eq:RR-cvx-lb-II-auxiliary-z})}{=}\left(\left[\prod_{j=1}^{\mathring{\pi}_{k}^{i}-1}(1-\eta\beta_{\pi_{k}^{j}})\right]\bx_{k}^{1}[1]-\eta\sqrt{\frac{n}{2}}\sigma_{\star}\sum_{j=1}^{\mathring{\pi}_{k}^{i}-1}\chi_{\pi_{k}^{j}}\prod_{l=j+1}^{\mathring{\pi}_{k}^{i}-1}(1-\eta\beta_{\pi_{k}^{l}})\right)_{+}\1\left[E_{k}^{i}\right]\nonumber \\
 & \overset{(b)}{=}\left(\left[\prod_{j=1}^{\mathring{\pi}_{k}^{i}-1}(1-\eta\beta_{\pi_{k}^{j}})\right]\bx_{k}^{1}[1]+\eta\sqrt{\frac{n}{2}}\sigma_{\star}\prod_{l=\mathring{\pi}_{k}^{n-1}+1}^{\mathring{\pi}_{k}^{i}-1}(1-\eta\beta_{\pi_{k}^{l}})\right)_{+}\1\left[E_{k}^{i}\right]\nonumber \\
 & \overset{(c)}{\geq}\left[\prod_{j=1}^{n}(1-\eta\beta_{\pi_{k}^{j}})\right]\left(\bx_{k}^{1}[1]+\eta\sqrt{\frac{n}{2}}\sigma_{\star}\right)_{+}\1\left[E_{k}^{i}\right]\overset{(\ref{eq:RR-cvx-lb-II-n>2-rho})}{=}\rho\left(\bx_{k}^{1}[1]+\eta\sqrt{\frac{n}{2}}\sigma_{\star}\right)_{+}\1\left[E_{k}^{i}\right],\label{eq:RR-cvx-lb-II-auxiliary-4}
\end{align}
where $(b)$ holds because, under the event $E_{k}^{i}$, $\chi_{\pi_{k}^{j}}=0$
for all $j\in\left\{ 1,\mydots,\mathring{\pi}_{k}^{n-1}-1\right\} \cup\left\{ \mathring{\pi}_{k}^{n-1}+1,\mydots,\mathring{\pi}_{k}^{i}-1\right\} $
and $\chi_{\pi_{k}^{j}}=-1$ when $j=\mathring{\pi}_{k}^{n-1}$, and
$(c)$ follows from combining the following two cases:
\begin{casenv}
\item $\bx_{k}^{1}[1]\geq0$. In this case, we have
\begin{align*}
 & \left(\left[\prod_{j=1}^{\mathring{\pi}_{k}^{i}-1}(1-\eta\beta_{\pi_{k}^{j}})\right]\bx_{k}^{1}[1]+\eta\sqrt{\frac{n}{2}}\sigma_{\star}\prod_{l=\mathring{\pi}_{k}^{n-1}+1}^{\mathring{\pi}_{k}^{i}-1}(1-\eta\beta_{\pi_{k}^{l}})\right)_{+}\\
= & \left[\prod_{j=1}^{\mathring{\pi}_{k}^{i}-1}(1-\eta\beta_{\pi_{k}^{j}})\right]\bx_{k}^{1}[1]+\eta\sqrt{\frac{n}{2}}\sigma_{\star}\prod_{l=\mathring{\pi}_{k}^{n-1}+1}^{\mathring{\pi}_{k}^{i}-1}(1-\eta\beta_{\pi_{k}^{l}})\\
\geq & \left[\prod_{j=1}^{n}(1-\eta\beta_{\pi_{k}^{j}})\right]\left(\bx_{k}^{1}[1]+\eta\sqrt{\frac{n}{2}}\sigma_{\star}\right)=\left[\prod_{j=1}^{n}(1-\eta\beta_{\pi_{k}^{j}})\right]\left(\bx_{k}^{1}[1]+\eta\sqrt{\frac{n}{2}}\sigma_{\star}\right)_{+}.
\end{align*}
\item $\bx_{k}^{1}[1]<0$. In this case, we have
\begin{align*}
 & \left(\left[\prod_{j=1}^{\mathring{\pi}_{k}^{i}-1}(1-\eta\beta_{\pi_{k}^{j}})\right]\bx_{k}^{1}[1]+\eta\sqrt{\frac{n}{2}}\sigma_{\star}\prod_{l=\mathring{\pi}_{k}^{n-1}+1}^{\mathring{\pi}_{k}^{i}-1}(1-\eta\beta_{\pi_{k}^{l}})\right)_{+}\\
= & \left[\prod_{j=\mathring{\pi}_{k}^{n-1}+1}^{\mathring{\pi}_{k}^{i}-1}(1-\eta\beta_{\pi_{k}^{j}})\right]\left(\left[\prod_{j=1}^{\mathring{\pi}_{k}^{n-1}}(1-\eta\beta_{\pi_{k}^{j}})\right]\bx_{k}^{1}[1]+\eta\sqrt{\frac{n}{2}}\sigma_{\star}\right)_{+}\\
\geq & \left[\prod_{j=\mathring{\pi}_{k}^{n-1}+1}^{\mathring{\pi}_{k}^{i}-1}(1-\eta\beta_{\pi_{k}^{j}})\right]\left(\bx_{k}^{1}[1]+\eta\sqrt{\frac{n}{2}}\sigma_{\star}\right)_{+}\geq\left[\prod_{j=1}^{n}(1-\eta\beta_{\pi_{k}^{j}})\right]\left(\bx_{k}^{1}[1]+\eta\sqrt{\frac{n}{2}}\sigma_{\star}\right)_{+}.
\end{align*}
\end{casenv}
Finally, we plug (\ref{eq:RR-cvx-lb-II-auxiliary-4}) back into (\ref{eq:RR-cvx-lb-II-auxiliary-3})
to have
\[
\sum_{i=1}^{n}\beta_{\pi_{k}^{i}}(\bz_{k}^{i}[1])_{+}\geq\rho\sum_{i=1}^{n-2}\beta_{i}\left(\bx_{k}^{1}[1]+\eta\sqrt{\frac{n}{2}}\sigma_{\star}\right)_{+}\1\left[E_{k}^{i}\right],
\]
which implies that
\begin{align*}
\E\left[\sum_{i=1}^{n}\beta_{\pi_{k}^{i}}(\bz_{k}^{i}[1])_{+}\right] & \geq\rho\sum_{i=1}^{n-2}\beta_{i}\E\left[\left(\bx_{k}^{1}[1]+\eta\sqrt{\frac{n}{2}}\sigma_{\star}\right)_{+}\1\left[E_{k}^{i}\right]\right]\\
 & \overset{(d)}{=}\rho\sum_{i=1}^{n-2}\beta_{i}\Pr\left[E_{k}^{i}\right]\E\left[\left(\bx_{k}^{1}[1]+\eta\sqrt{\frac{n}{2}}\sigma_{\star}\right)_{+}\right]\\
 & \overset{(\ref{eq:RR-cvx-lb-II-n>2-beta}),(\ref{eq:RR-cvx-lb-II-auxiliary-2})}{=}\frac{\rho n\bar{\beta}}{6}\E\left[\left(\bx_{k}^{1}[1]+\eta\sqrt{\frac{n}{2}}\sigma_{\star}\right)_{+}\right],
\end{align*}
where $(d)$ follows from the independence of $\pi_{k}$ and $\bx_{k}^{1}$.

\subparagraph{The third inequality.}

Using an argument similar to that for the second inequality, we obtain
\begin{align*}
\E\left[\sum_{i=1}^{n}(n-i)\beta_{\pi_{k}^{i}}(\bz_{k}^{i}[1])_{+}\right] & \geq\rho\sum_{i=1}^{n-2}\beta_{i}\E\left[(n-\mathring{\pi}_{k}^{i})\1\left[E_{k}^{i}\right]\right]\E\left[\left(\bx_{k}^{1}[1]+\eta\sqrt{\frac{n}{2}}\sigma_{\star}\right)_{+}\right]\\
 & =\frac{\rho n(n-1)\bar{\beta}}{12}\E\left[\left(\bx_{k}^{1}[1]+\eta\sqrt{\frac{n}{2}}\sigma_{\star}\right)_{+}\right],
\end{align*}
where the last step is due to (\ref{eq:RR-cvx-lb-II-n>2-beta}) and
the fact that, for any $i\in\left[n-2\right]$,
\begin{align*}
\E\left[(n-\mathring{\pi}_{k}^{i})\1\left[E_{k}^{i}\right]\right] & =n\Pr\left[E_{k}^{i}\right]-\E\left[\mathring{\pi}_{k}^{i}\1\left[E_{k}^{i}\right]\right]\overset{(\ref{eq:RR-cvx-lb-II-auxiliary-2})}{=}\frac{n}{6}-\E\left[\mathring{\pi}_{k}^{i}\1\left[E_{k}^{i}\right]\right]\\
 & =\frac{n}{6}-\sum_{l=2}^{n-1}\frac{l(l-1)(n-l)}{n(n-1)(n-2)}=\frac{n-1}{12}.
\end{align*}
\end{proof}

\paragraph{Analysis when $n=2$.}

We now turn to the case $n=2$. Compared with the case $n\geq3$,
the construction and analysis are considerably simpler, as the problem
involves only two component functions and the effect of the random
permutation can be characterized explicitly. This allows us to establish
the desired lower bound through a direct calculation.
\begin{thm}
\label{thm:RR-cvx-lb-II-n=00003D2}Given $d\in\N$, $n=2$, $D>0$,
$\left\{ L_{i}>0\right\} _{i=1}^{n}$, and $\sigma_{\star}\geq0$,
suppose $\RR$ is employed with $\left\{ \eta_{k}=\eta\right\} _{k=1}^{K}$
for $\eta>0$, let $\bar{\bx}_{K}\defeq\frac{1}{nK}\sum_{k=1}^{K}\sum_{i=1}^{n}\bx_{k}^{i}$
denote the average iterate generated by $\ShufflingSGD$ (Algorithm
\ref{alg:Shuffling-SGD}), initialized at $\bx_{1}^{1}=\bzero$, over
$K\in\N$ epochs, then we have
\begin{eqnarray*}
\sup_{f\in\F(d,n,D,L_{1:n},\sigma_{\star})}\E\left[f(\bar{\bx}_{K})-f_{\star}\right]\gtrsim\min\left\{ 1,\eta n\Lavg\right\} \eta\sigma_{\star}^{2} & \text{when} & K\eta n\Lavg\geq1.
\end{eqnarray*}
\end{thm}
\begin{proof}
W.l.o.g., we assume $L_{1}\geq L_{2}$ in the following proof. We
first define 
\begin{equation}
\beta\defeq\min\left\{ \frac{1}{\eta n},\Lavg\right\} =\min\left\{ \frac{1}{2\eta},\Lavg\right\} .\label{eq:RR-cvx-lb-II-n=00003D2-beta}
\end{equation}
Next, we introduce
\begin{eqnarray*}
f_{1}(\bx)\defeq\beta\frac{\left(\bx[1]\right)^{2}}{2}-\sigma_{\star}\bx[1] & \text{and} & f_{2}(\bx)\defeq\sigma_{\star}\bx[1].
\end{eqnarray*}
The corresponding finite-sum objective is
\begin{equation}
f(\bx)\defeq\frac{f_{1}(\bx)+f_{2}(\bx)}{2}=\frac{\beta}{4}\left(\bx[1]\right)^{2}\overset{(\ref{eq:RR-cvx-lb-II-n=00003D2-beta})}{=}\frac{\min\left\{ \frac{1}{n\eta},\Lavg\right\} \left(\bx[1]\right)^{2}}{4}.\label{eq:RR-cvx-lb-II-n=00003D2-f}
\end{equation}
As before, one can verify that $f\in\F(d,n,D,L_{1:n},\sigma_{\star})$.

Again, we consider two cases:
\begin{casenv}
\item $K\geq2$. In this case, a direct calculation gives us
\[
\E\left[\bx_{k+1}^{1}[1]\right]=(1-\eta\beta)\E\left[\bx_{k}^{1}[1]\right]+\frac{\eta^{2}\beta\sigma_{\star}}{2},\forall k\in\left[K\right].
\]
Since $\bx_{1}^{1}=\bzero$, unrolling this recursion yields
\begin{equation}
\E\left[\bx_{k}^{1}[1]\right]=\frac{\eta\sigma_{\star}}{2}(1-(1-\eta\beta)^{k-1}),\forall k\in\left[K+1\right].\label{eq:RR-cvx-lb-II-n=00003D2-1}
\end{equation}
Next, we observe
\begin{equation}
\E\left[\bx_{k}^{2}[1]\right]=\left(1-\frac{\eta\beta}{2}\right)\E\left[\bx_{k}^{1}[1]\right]\overset{\eta\beta\leq\frac{1}{2}}{\geq}\frac{3}{4}\E\left[\bx_{k}^{1}[1]\right].\label{eq:RR-cvx-lb-II-n=00003D2-2}
\end{equation}
Therefore, we can obtain that
\begin{align}
\E\left[\bar{\bx}_{K}[1]\right] & =\frac{1}{2K}\sum_{k=1}^{K}\sum_{i=1}^{2}\E\left[\bx_{k}^{i}[1]\right]\overset{(\ref{eq:RR-cvx-lb-II-n=00003D2-2})}{\geq}\frac{7}{8K}\sum_{k=1}^{K}\E\left[\bx_{k}^{1}[1]\right]\overset{(\ref{eq:RR-cvx-lb-II-n=00003D2-1})}{\geq}\frac{7\eta\sigma_{\star}}{16K}\sum_{k=1}^{K}1-(1-\eta\beta)^{k-1}\nonumber \\
 & \geq\frac{7\eta\sigma_{\star}}{16K}\sum_{k=1}^{K}1-\exp(-\eta\beta(k-1))\geq\frac{7\eta\sigma_{\star}}{16K}\sum_{k=\left\lceil K/2\right\rceil +1}^{K}1-\exp(-\eta\beta(k-1))\nonumber \\
 & \geq\frac{7\eta\sigma_{\star}}{16}\cdot\frac{\left\lfloor K/2\right\rfloor }{K}\cdot(1-\exp(-K\eta\beta/2))\overset{K\geq2}{\geq}\frac{7\eta\sigma_{\star}}{48}(1-\exp(-K\eta\beta/2))\nonumber \\
 & \overset{(a)}{\geq}\frac{7\min\left\{ 1-\exp(-1/2),\min\left\{ \frac{K}{2},\frac{K\eta n\Lavg}{2}\right\} /4\right\} }{48\sqrt{2}}\eta\sqrt{n}\sigma_{\star}\geq\frac{7\min\left\{ 1,K\eta n\Lavg\right\} }{384\sqrt{2}}\eta\sqrt{n}\sigma_{\star},\label{eq:RR-cvx-lb-II-n=00003D2-3}
\end{align}
where $(a)$ follows from $n=2$, $1-\exp(-x)\ge\min\left\{ 1-\exp(-1/2),x/2\right\} ,\forall x\geq0$,
and $K\eta\beta\overset{(\ref{eq:RR-cvx-lb-II-n=00003D2-beta})}{=}\min\left\{ \frac{K}{2},K\eta\Lavg\right\} \overset{n=2}{=}\min\left\{ \frac{K}{2},\frac{K\eta n\Lavg}{2}\right\} $.
Hence, we have
\begin{align*}
\E\left[f(\bar{\bx}_{K})-f_{\star}\right] & \geq f(\E\left[\bar{\bx}_{K}\right])-f_{\star}\overset{(\ref{eq:RR-cvx-lb-II-n=00003D2-f})}{\geq}\frac{\min\left\{ \frac{1}{n\eta},\Lavg\right\} \left(\E\left[\bar{\bx}_{K}[1]\right]\right)^{2}}{4}\\
 & \overset{(\ref{eq:RR-cvx-lb-II-n=00003D2-3})}{\gtrsim}\min\left\{ 1,\eta n\Lavg\right\} \left(\min\left\{ 1,K\eta n\Lavg\right\} \right)^{2}\eta\sigma_{\star}^{2},
\end{align*}
where the first step is by Jensen's inequality.
\item $K=1$. In this case, we know $\bar{\bx}_{K}[1]=\bar{\bx}_{1}[1]=\frac{\bx_{1}^{1}[1]+\bx_{1}^{2}[1]}{2}\overset{\bx_{1}^{1}[1]=0}{=}\frac{\bx_{1}^{2}[1]}{2}$
and
\[
\Pr\left[\bx_{1}^{2}[1]=\eta\sigma_{\star}\right]=\Pr\left[\bx_{1}^{2}[1]=-\eta\sigma_{\star}\right]=\frac{1}{2}.
\]
Hence, we have
\begin{align*}
\E\left[f(\bar{\bx}_{K})-f_{\star}\right] & \overset{(\ref{eq:RR-cvx-lb-II-n=00003D2-f})}{\geq}\frac{\min\left\{ \frac{1}{\eta n},\Lavg\right\} \E\left[\left(\bar{\bx}_{K}[1]\right)^{2}\right]}{4}=\frac{\min\left\{ \frac{1}{\eta n},\Lavg\right\} \eta^{2}\sigma_{\star}^{2}}{16}\\
 & \overset{n=2}{=}\frac{\min\left\{ 1,\eta n\Lavg\right\} \eta\sigma_{\star}^{2}}{32}\gtrsim\min\left\{ 1,\eta n\Lavg\right\} \eta\sigma_{\star}^{2}\\
 & \geq\min\left\{ 1,\eta n\Lavg\right\} \left(\min\left\{ 1,K\eta n\Lavg\right\} \right)^{2}\eta\sigma_{\star}^{2}.
\end{align*}
\end{casenv}
Combining the two cases above, we conclude that, for $n=2$,
\[
\sup_{f\in\F(d,n,D,L_{1:n},\sigma_{\star})}\E\left[f(\bar{\bx}_{K})-f_{\star}\right]\gtrsim\min\left\{ 1,\eta n\Lavg\right\} \left(\min\left\{ 1,K\eta n\Lavg\right\} \right)^{2}\eta\sigma_{\star}^{2}.
\]
Consequently, for $n=2$,
\begin{eqnarray*}
\sup_{f\in\F(d,n,D,L_{1:n},\sigma_{\star})}\E\left[f(\bar{\bx}_{K})-f_{\star}\right]\gtrsim\min\left\{ 1,\eta n\Lavg\right\} \eta\sigma_{\star}^{2} & \text{when} & K\eta n\Lavg\geq1.
\end{eqnarray*}
\end{proof}

\paragraph{Final proof.}

With Theorems \ref{thm:RR-cvx-lb-II-n>2} and \ref{thm:RR-cvx-lb-II-n=00003D2}
in place, we can immediately prove Theorem \ref{thm:RR-cvx-lb-II}.

\begin{proof}[Proof of Theorem \ref{thm:RR-cvx-lb-II}]
The result follows directly from Theorems \ref{thm:RR-cvx-lb-II-n>2}
and \ref{thm:RR-cvx-lb-II-n=00003D2}.
\end{proof}

\subsection{Proof of Theorem \ref{thm:RR-cvx-lb}\label{subsec:RR-cvx-lb-Proof}}

The preceding results allow us to establish Theorem \ref{thm:RR-cvx-lb}
now.

\begin{proof}[Proof of Theorem \ref{thm:RR-cvx-lb}]
Combining Theorems \ref{thm:RR-cvx-lb-I} and \ref{thm:RR-cvx-lb-II}
completes the proof.
\end{proof}

\section{Theoretical Analysis: Lower Bound for $\protect\SGD$\label{sec:SGD-cvx-lb}}

In this section, we prove the algorithm-dependent lower bound for
$\SGD$ with a fixed stepsize $\eta$, i.e., Theorem \ref{thm:SGD-cvx-lb}.
Our proof proceeds by constructing a family of quadratic finite-sum
problems tailored to $\eta$ and showing that, for any such stepsize,
$\SGD$ necessarily incurs an optimization error of the desired order.
Similar to the previous proof for $\RR$, the construction separates
the effects of $D$ (the initial distance to the minimizer) and $\sigma_{\star}^{2}$
(the variance of the gradient at the optimum), allowing us to derive
the corresponding lower-bound terms individually and then combine
them.

Compared with the analysis of the lower bound for $\RR$, the case
of $\SGD$ is relatively straightforward because the sampled indices
are independent across iterations, making the calculations much easier.
Nevertheless, the result is of independent interest: to the best of
our knowledge, it provides the first algorithm-dependent lower bound
for $\SGD$ that explicitly characterizes its performance for each
prescribed stepsize $\eta$. In particular, rather than optimizing
over the stepsize in advance, our lower bound holds for every fixed
stepsize and thus directly captures the stepsize-dependent limitations
of $\SGD$.

\begin{proof}[Proof of Theorem \ref{thm:SGD-cvx-lb}]
For simplicity, we write the update rule of $\SGD$ as
\begin{eqnarray*}
\bx_{t+1}=\bx_{t}-\eta\nabla f_{\pi_{t}}(\bx_{t}),\forall t\in\left[T\right] & \text{where} & T\defeq nK,
\end{eqnarray*}
where $\bx_{1}=\bzero$, $\eta>0$, and $\pi_{t}$ is sampled independently
and uniformly from $\left[n\right]$. In the following proof, $\bv$
denotes an arbitrary unit vector in $\R^{d}$.

We first define
\[
\beta_{i}\defeq\min\left\{ \frac{1}{2\eta T},\Lavg\right\} \frac{L_{i}}{\Lavg},\forall i\in\left[n\right].
\]
Clearly, each $\beta_{i}$ satisfies $0<\beta_{i}\leq L_{i}$. Moreover,
we write
\begin{equation}
\bar{\beta}\defeq\frac{1}{n}\sum_{i=1}^{n}\beta_{i}=\min\left\{ \frac{1}{2\eta T},\Lavg\right\} .\label{eq:SGD-cvx-lb-beta}
\end{equation}
We then fix $R\in\left[0,D\right]$ and $\nu\in\left[0,\sigma_{\star}\right]$,
both of which will be decided later, and define
\[
f_{i}(\bx)\defeq\frac{\beta_{i}\left\Vert \bx-R\bv\right\Vert ^{2}}{2}+\chi_{i}\nu\left\langle \bv,\bx\right\rangle ,\forall i\in\left[n\right],
\]
where $\chi_{i}\defeq\begin{cases}
1 & i\leq n\text{ is odd}\\
-1 & i\leq n\text{ is even}
\end{cases}$ for even $n$, and $\chi_{i}\defeq\begin{cases}
1 & i\leq n-1\text{ is odd}\\
-1 & i\leq n-1\text{ is even}\\
0 & i=n
\end{cases}$ for odd $n$. The corresponding finite-sum objective is
\begin{equation}
f(\bx)\defeq\frac{1}{n}\sum_{i=1}^{n}f_{i}(\bx)=\frac{\bar{\beta}\left\Vert \bx-R\bv\right\Vert ^{2}}{2}.\label{eq:SGD-cvx-lb-f}
\end{equation}
Note that $\bx_{\star}=R\bv\Rightarrow\left\Vert \bx_{\star}\right\Vert =R\leq D$,
each $f_{i}$ is differentiable, convex, and $L_{i}$-smooth (since
$\beta_{i}\leq L_{i}$), and $\frac{1}{n}\sum_{i=1}^{n}\left\Vert \nabla f_{i}(\bx_{\star})\right\Vert ^{2}=\frac{1}{n}\sum_{i=1}^{n}\chi_{i}^{2}\nu^{2}\leq\sigma_{\star}^{2}$.
Therefore, $f$ belongs to the class $\F(d,n,D,L_{1:n},\sigma_{\star})$.

We then have
\[
\nabla f_{\pi_{t}}(\bx_{t})=\beta_{\pi_{t}}(\bx_{t}-R\bv)+\chi_{\pi_{t}}\nu\bv,\forall t\in\left[T\right],
\]
which implies that
\[
\bx_{t+1}=(1-\eta\beta_{\pi_{t}})\bx_{t}+\eta(\beta_{\pi_{t}}R-\chi_{\pi_{t}}\nu)\bv,\forall t\in\left[T\right].
\]
Since $\bx_{1}=\bzero$, unrolling the recursion yields that
\[
\bx_{t}=\left[\left(1-\prod_{s=1}^{t-1}(1-\eta\beta_{\pi_{s}})\right)R-\sum_{s=1}^{t-1}\chi_{\pi_{s}}\prod_{r=s+1}^{t-1}(1-\eta\beta_{\pi_{r}})\eta\nu\right]\bv,\forall t\in\left[T\right].
\]
Therefore, we have
\begin{equation}
\bar{\bx}_{K}=\frac{1}{T}\sum_{t=1}^{T}\bx_{t}=\left[\left(1-\frac{1}{T}\sum_{t=1}^{T}\prod_{s=1}^{t-1}(1-\eta\beta_{\pi_{s}})\right)R-\frac{1}{T}\sum_{t=1}^{T-1}\chi_{\pi_{t}}\sum_{s=t+1}^{T}\prod_{r=t+1}^{s-1}(1-\eta\beta_{\pi_{r}})\eta\nu\right]\bv.\label{eq:SGD-cvx-lb-x}
\end{equation}

Next, we consider two cases:
\begin{casenv}
\item $R=D$ and $\nu=0$. In this case, we have
\begin{align}
\bar{\bx}_{K} & \overset{(\ref{eq:SGD-cvx-lb-x})}{=}\left(1-\frac{1}{T}\sum_{t=1}^{T}\prod_{s=1}^{t-1}(1-\eta\beta_{\pi_{s}})\right)D\bv\nonumber \\
\Rightarrow\E\left[\bar{\bx}_{K}\right] & =\left(1-\frac{1}{T}\sum_{t=1}^{T}\prod_{s=1}^{t-1}(1-\eta\bar{\beta})\right)D\bv=\left(1-\frac{1-(1-\eta\bar{\beta})^{T}}{\eta\bar{\beta}T}\right)D\bv.\label{eq:SGD-cvx-lb-det-1}
\end{align}
Therefore, we know
\begin{align*}
\E\left[f(\bar{\bx}_{K})-f_{\star}\right] & \overset{(a)}{\geq}f(\E\left[\bar{\bx}_{K}\right])-f_{\star}\overset{(\ref{eq:SGD-cvx-lb-f}),(\ref{eq:SGD-cvx-lb-det-1}),\left\Vert \bv\right\Vert =1}{=}\frac{(1-(1-\eta\bar{\beta})^{T})^{2}D^{2}}{2\eta^{2}T^{2}\bar{\beta}}\overset{\eta\bar{\beta}\in\left[0,1\right]}{\geq}\frac{(1-\exp(-T\eta\bar{\beta}))^{2}D^{2}}{2\eta^{2}T^{2}\bar{\beta}}\\
 & \overset{(b)}{\geq}\frac{\bar{\beta}D^{2}}{8}\overset{(\ref{eq:SGD-cvx-lb-beta})}{=}\frac{D^{2}}{8}\min\left\{ \frac{1}{2\eta T},\Lavg\right\} \geq\frac{D^{2}}{16\eta T}\min\left\{ 1,T\eta\Lavg\right\} =\frac{D^{2}}{16\eta nK}\min\left\{ 1,K\eta n\Lavg\right\} ,
\end{align*}
where $(a)$ is due to Jensen's inequality and $(b)$ follows from
$1-\exp(-x)\geq\frac{x}{2},\forall x\in\left[0,\frac{1}{2}\right]$
and $0\leq T\eta\bar{\beta}\overset{(\ref{eq:SGD-cvx-lb-beta})}{\leq}\frac{1}{2}$.
\item $R=0$ and $\nu=\sigma_{\star}$. In this case, we have
\begin{equation}
\bar{\bx}_{K}\overset{(\ref{eq:SGD-cvx-lb-x})}{=}-\frac{1}{T}\sum_{t=1}^{T-1}\chi_{\pi_{t}}\sum_{s=t+1}^{T}\prod_{r=t+1}^{s-1}(1-\eta\beta_{\pi_{r}})\eta\sigma_{\star}\bv=-\frac{1}{T}\sum_{t=1}^{T-1}\chi_{\pi_{t}}\alpha_{t}\eta\sigma_{\star}\bv,\label{eq:SGD-cvx-lb-stoc-1}
\end{equation}
where $\alpha_{t}\defeq\sum_{s=t+1}^{T}\prod_{r=t+1}^{s-1}(1-\eta\beta_{\pi_{r}}),\forall t\in\left[T\right]$.
Therefore, we know
\begin{align*}
\E\left[f(\bar{\bx}_{K})-f_{\star}\right] & \overset{(\ref{eq:SGD-cvx-lb-f}),(\ref{eq:SGD-cvx-lb-stoc-1}),\left\Vert \bv\right\Vert =1}{=}\frac{\bar{\beta}\E\left[(\sum_{t=1}^{T-1}\chi_{\pi_{t}}\alpha_{t})^{2}\right]\eta^{2}\sigma_{\star}^{2}}{2T^{2}}\overset{(c)}{=}\frac{\bar{\beta}\sum_{t=1}^{T-1}\E\left[\chi_{\pi_{t}}^{2}\alpha_{t}^{2}\right]\eta^{2}\sigma_{\star}^{2}}{2T^{2}}\\
 & \overset{(d)}{\geq}\frac{\sum_{t=1}^{T-1}(1-\exp(-(T-t)\eta\bar{\beta}))^{2}\sigma_{\star}^{2}}{3\bar{\beta}T^{2}}\overset{(e)}{\geq}\frac{(1-\exp(-T\eta\bar{\beta}/2))^{2}\sigma_{\star}^{2}}{9\bar{\beta}T}\overset{(f)}{\geq}\frac{\eta\sigma_{\star}^{2}}{144}\cdot T\eta\bar{\beta}\\
 & \overset{(\ref{eq:SGD-cvx-lb-beta})}{=}\frac{\eta\sigma_{\star}^{2}}{144}\min\left\{ \frac{1}{2},T\eta\Lavg\right\} \geq\frac{\eta\sigma_{\star}^{2}}{288}\min\left\{ 1,T\eta\Lavg\right\} =\frac{\eta\sigma_{\star}^{2}}{288}\min\left\{ 1,Kn\eta\Lavg\right\} ,
\end{align*}
where $(c)$ is due to $\E\left[\chi_{\pi_{t}}\chi_{\pi_{t^{\prime}}}\alpha_{t}\alpha_{t^{\prime}}\right]=0,\forall t\neq t^{\prime}\in\left[T-1\right]$,
$(d)$ is by the independence of $\chi_{\pi_{t}}^{2}$ and $\alpha_{t}^{2}$,
$\E\left[\chi_{\pi_{t}}^{2}\right]=\begin{cases}
1 & n\text{ is even}\\
\frac{n-1}{n} & n\text{ is odd}
\end{cases}\geq\frac{2}{3}$ since $n\geq2$, and
\[
\E\left[\alpha_{t}^{2}\right]\geq\left(\E\left[\alpha_{t}\right]\right)^{2}=\left(\sum_{s=t+1}^{T}\prod_{r=t+1}^{s-1}(1-\eta\bar{\beta})\right)^{2}=\frac{(1-(1-\eta\bar{\beta})^{T-t})^{2}}{\eta^{2}\bar{\beta}^{2}}\overset{\eta\bar{\beta}\in\left[0,1\right]}{\geq}\frac{(1-\exp(-(T-t)\eta\bar{\beta}))^{2}}{\eta^{2}\bar{\beta}^{2}},
\]
 $(e)$ holds by
\begin{align*}
\sum_{t=1}^{T-1}(1-\exp(-(T-t)\eta\bar{\beta}))^{2} & \geq\sum_{t=1}^{\left\lfloor \frac{T}{2}\right\rfloor }(1-\exp(-(T-t)\eta\bar{\beta}))^{2}\\
 & \geq\left\lfloor \frac{T}{2}\right\rfloor (1-\exp(-T\eta\bar{\beta}/2))^{2}\overset{T=nK\geq2}{\ge}\frac{T}{3}(1-\exp(-T\eta\bar{\beta}/2))^{2},
\end{align*}
and $(f)$ follows from $1-\exp(-x)\geq\frac{x}{2},\forall x\in\left[0,\frac{1}{2}\right]$
and $0\leq\frac{T\eta\bar{\beta}}{2}\overset{(\ref{eq:SGD-cvx-lb-beta})}{\leq}\frac{1}{2}$.
\end{casenv}
Combining the two cases above, we conclude that
\[
\sup_{f\in\F(d,n,D,L_{1:n},\sigma_{\star})}\E\left[f(\bar{\bx}_{K})-f_{\star}\right]\gtrsim\left(\frac{D^{2}}{\eta nK}+\eta\sigma_{\star}^{2}\right)\min\left\{ 1,Kn\eta\Lavg\right\} .
\]
Consequently,
\begin{eqnarray*}
\sup_{f\in\F(d,n,D,L_{1:n},\sigma_{\star})}\E\left[f(\bar{\bx}_{K})-f_{\star}\right]\gtrsim\frac{D^{2}}{\eta nK}+\eta\sigma_{\star}^{2} & \text{when} & K\eta n\Lavg\geq1.
\end{eqnarray*}
\end{proof}

\section{Theoretical Analysis: Application of $\protect\RR$}

In this section, we first prove Theorem \ref{thm:RR-QG-ub}, which
provides a convergence guarantee for $\RShufflingSGD$ with $\RR$
for functions satisfying the quadratic growth condition and serves
as an application of our new Theorem \ref{thm:RR-cvx-ub}. We then
invoke Theorem \ref{thm:RR-QG-ub} to establish Corollary \ref{cor:RR-QG-complexity},
which gives the sample complexity, measured by the number of individual
gradient evaluations, required to obtain an $\varepsilon$-accurate
solution, i.e., to find a point $\bx\in\R^{d}$ such that $\E\left[f(\bx)-f_{\star}\right]\lesssim\varepsilon$.

\subsection{Proof of Theorem \ref{thm:RR-QG-ub}\label{subsec:RR-QG-ub-Proof}}

\begin{proof}[Proof of Theorem \ref{thm:RR-QG-ub}]
Note that $f$ is real-valued. Hence, under Assumptions \ref{assu:minimizer}
and \ref{assu:convexity}, $\X_{\star}=\arg\min_{\bx\in\R^{d}}f(\bx)$
is a nonempty closed convex set. Therefore, for any $\bx\in\R^{d}$,
the problem $\inf_{\bx_{\star}\in\X_{\star}}\left\Vert \bx-\bx_{\star}\right\Vert $
admits a unique solution, which we denote by $\left[\bx\right]_{\star}\defeq\arg\min_{\bx_{\star}\in\X_{\star}}\left\Vert \bx-\bx_{\star}\right\Vert $
in the following. Moreover, we let $\F_{t}\defeq\sigma(\text{permutations used in the first }t\text{ restart stages})$,
i.e., the $\sigma$-algebra generated by the permutations in the first
$t$ restarts.

By Theorem \ref{thm:RR-cvx-ub}, given $t\in\left[T\right]$ and any
$\F_{t-1}$-measurable $\bx_{\star}\in\X_{\star}$,
\[
\E\left[f(\bx_{t+1})-f_{\star}\mid\F_{t-1}\right]\lesssim\frac{\left\Vert \bx_{t}-\bx_{\star}\right\Vert ^{2}}{(2-\eta\Lmax)\eta nK}+\frac{\min\left\{ 1,\eta n\Lavg\right\} \eta\sigma_{\star}^{2}}{(2-\eta\Lmax)^{2}}.
\]
Note that $\bx_{t}$ is $\F_{t-1}$-measurable, and hence $\left[\bx_{t}\right]_{\star}$
is also $\F_{t-1}$-measurable. We can therefore replace $\bx_{\star}$
with $\left[\bx_{t}\right]_{\star}$ in the above inequality to obtain
\begin{align*}
\E\left[f(\bx_{t+1})-f_{\star}\mid\F_{t-1}\right] & \lesssim\frac{\dist^{2}(\bx_{t},\X_{\star})}{(2-\eta\Lmax)\eta nK}+\frac{\min\left\{ 1,\eta n\Lavg\right\} \eta\sigma_{\star}^{2}}{(2-\eta\Lmax)^{2}}\\
\Rightarrow\E\left[f(\bx_{t+1})-f_{\star}\right] & \lesssim\frac{\E\left[\dist^{2}(\bx_{t},\X_{\star})\right]}{(2-\eta\Lmax)\eta nK}+\frac{\min\left\{ 1,\eta n\Lavg\right\} \eta\sigma_{\star}^{2}}{(2-\eta\Lmax)^{2}}.
\end{align*}
Let $C>0$ denote the hidden numerical constant in the above inequality
(equivalently, the hidden constant in Theorem \ref{thm:RR-cvx-ub}).
When $K=\left\lceil 4CK_{\star}\right\rceil =\left\lceil \frac{4C}{(2-\eta\Lmax)\eta n\mu}\right\rceil $,
we have
\begin{equation}
\E\left[f(\bx_{t+1})-f_{\star}\right]\leq\frac{\mu\E\left[\dist^{2}(\bx_{t},\X_{\star})\right]}{4}+C\frac{\min\left\{ 1,\eta n\Lavg\right\} \eta\sigma_{\star}^{2}}{(2-\eta\Lmax)^{2}}.\label{eq:RR-QG-ub-1}
\end{equation}

Next, by Assumption \ref{assu:quadratic-growth},
\[
\frac{\mu}{2}\dist^{2}(\bx_{t+1},\X_{\star})\leq f(\bx_{t+1})-f_{\star},
\]
which implies that
\[
\E\left[\dist^{2}(\bx_{t+1},\X_{\star})\right]\leq\frac{\E\left[\dist^{2}(\bx_{t},\X_{\star})\right]}{2}+2C\frac{\min\left\{ 1,\eta n\Lavg\right\} \eta\sigma_{\star}^{2}}{\mu(2-\eta\Lmax)^{2}}.
\]
Unrolling the above recursion gives
\begin{equation}
\E\left[\dist^{2}(\bx_{t},\X_{\star})\right]\leq\frac{\dist^{2}(\bx_{1},\X_{\star})}{2^{t-1}}+4C\frac{\min\left\{ 1,\eta n\Lavg\right\} \eta\sigma_{\star}^{2}}{\mu(2-\eta\Lmax)^{2}},\forall t\in\left[T+1\right].\label{eq:RR-QG-ub-2}
\end{equation}

Finally, we plug (\ref{eq:RR-QG-ub-2}) back into (\ref{eq:RR-QG-ub-1})
to obtain that
\[
\E\left[f(\bx_{t+1})-f_{\star}\right]\lesssim\frac{\mu\dist^{2}(\bx_{1},\X_{\star})}{2^{t}}+\frac{\min\left\{ 1,\eta n\Lavg\right\} \eta\sigma_{\star}^{2}}{(2-\eta\Lmax)^{2}},\forall t\in\left[T\right].
\]
In particular, taking $t=T$ in the above inequality yields the desired
result.
\end{proof}

\subsection{Proof of Corollary \ref{cor:RR-QG-complexity}\label{subsec:RR-QG-complexity-Proof}}

\begin{proof}[Proof of Corollary \ref{cor:RR-QG-complexity}]
By Theorem \ref{thm:RR-QG-ub}, we have
\[
\E\left[f(\bx_{T+1})-f_{\star}\right]\lesssim\frac{\mu\dist^{2}(\bx_{1},\X_{\star})}{2^{T}}+\frac{\min\left\{ 1,\eta n\Lavg\right\} \eta\sigma_{\star}^{2}}{(2-\eta\Lmax)^{2}}.
\]

When $T=T_{\star}=\left\lceil \log_{2}\frac{\mu\dist^{2}(\bx_{1},\X_{\star})}{\varepsilon}\right\rceil $,
we have 
\[
\frac{\mu\dist^{2}(\bx_{1},\X_{\star})}{2^{T_{\star}}}\leq\varepsilon.
\]

When $\eta=\eta_{\star}=\min\left\{ \frac{1}{\Lmax},\max\left\{ \frac{\varepsilon}{\sigma_{\star}^{2}},\frac{1}{\sigma_{\star}}\sqrt{\frac{\varepsilon}{n\Lavg}}\right\} \right\} $,
we observe
\[
\frac{\min\left\{ 1,\eta_{\star}n\Lavg\right\} \eta_{\star}\sigma_{\star}^{2}}{(2-\eta_{\star}\Lmax)^{2}}\overset{\eta_{\star}\Lmax\leq1}{\leq}\min\left\{ 1,\eta_{\star}n\Lavg\right\} \eta_{\star}\sigma_{\star}^{2}.
\]
Next, we use $\eta_{\star}\leq\max\left\{ \frac{\varepsilon}{\sigma_{\star}^{2}},\frac{1}{\sigma_{\star}}\sqrt{\frac{\varepsilon}{n\Lavg}}\right\} $
to obtain that
\[
\min\left\{ 1,\eta_{\star}n\Lavg\right\} \leq\min\left\{ 1,\max\left\{ \frac{n\Lavg\varepsilon}{\sigma_{\star}^{2}},\frac{\sqrt{n\Lavg\varepsilon}}{\sigma_{\star}}\right\} \right\} =\min\left\{ 1,\frac{\sqrt{n\Lavg\varepsilon}}{\sigma_{\star}}\right\} 
\]
and
\[
\eta_{\star}\sigma_{\star}^{2}\leq\max\left\{ \varepsilon,\sigma_{\star}\sqrt{\frac{\varepsilon}{n\Lavg}}\right\} .
\]
The above three inequalities together imply that
\[
\frac{\min\left\{ 1,\eta_{\star}n\Lavg\right\} \eta_{\star}\sigma_{\star}^{2}}{(2-\eta_{\star}\Lmax)^{2}}\leq\max\left\{ \varepsilon,\sigma_{\star}\sqrt{\frac{\varepsilon}{n\Lavg}}\right\} \min\left\{ 1,\frac{\sqrt{n\Lavg\varepsilon}}{\sigma_{\star}}\right\} \leq\varepsilon.
\]

Therefore, we know
\[
\E\left[f(\bx_{T_{\star}+1})-f_{\star}\right]\lesssim\varepsilon.
\]
Finally, for sufficiently small $\varepsilon$ such that $K_{\star}\geq1$,
the number of individual gradient evaluations is
\begin{align*}
nKT_{\star} & \asymp nK_{\star}T_{\star}\asymp n\cdot\frac{1}{(2-\eta_{\star}\Lmax)\eta_{\star}n\mu}\cdot\log\frac{\mu\dist^{2}(\bx_{1},\X_{\star})}{\varepsilon}\overset{\eta_{\star}\Lmax\in\left[0,1\right]}{\asymp}\frac{1}{\eta_{\star}\mu}\log\frac{\mu\dist^{2}(\bx_{1},\X_{\star})}{\varepsilon}\\
 & \asymp\left(\frac{\Lmax}{\mu}+\min\left\{ \frac{\sigma_{\star}^{2}}{\mu\varepsilon},\frac{\sigma_{\star}}{\mu}\sqrt{\frac{n\Lavg}{\varepsilon}}\right\} \right)\log\frac{\mu\dist^{2}(\bx_{1},\X_{\star})}{\varepsilon}.
\end{align*}
\end{proof}

\section{Conclusion and Future Work}

In this work, we prove that $\ShufflingSGD$ under $\RR$ dominates
standard $\SGD$ in smooth convex optimization under any reasonable
stepsize after any finite number of epochs. Our Theorem \ref{thm:RR-cvx-ub}
follows from combining two novel convergence results, whose analysis
may each be of independent interest. We also establish new algorithm-dependent
lower bounds for $\RR$ and $\SGD$ in Theorems \ref{thm:RR-cvx-lb}
and \ref{thm:SGD-cvx-lb}. Together, these results provide (nearly)
sharp worst-case characterizations of $\RR$ and standard $\SGD$.
Moreover, we combine our new convergence result with a restart framework
and apply it to functions satisfying the quadratic growth condition,
yielding a new convergence rate in Theorem \ref{thm:RR-QG-ub} and
a better sample complexity bound in Corollary \ref{cor:RR-QG-complexity}
that improves upon the best known results, even those established
under stronger assumptions.

Our work suggests several new directions for future research. First,
our current proof of Theorem \ref{thm:RR-cvx-ub} is split into two
distinct parts. It is therefore worthwhile to investigate whether
Theorem \ref{thm:RR-cvx-ub} can be proved via a unified analysis.
Second, the current upper and lower bounds for both $\RR$ and $\SGD$
differ by factors polynomial in $1/(2-\eta\Lmax)$, and closing this
gap remains an important direction. Finally, our results under the
quadratic growth condition rely on a restart scheme. It would be interesting
to determine whether $\ShufflingSGD$ can achieve comparable finite-time
convergence guarantees without restarting.

\section*{Timeline and AI Disclosure}

The improved stepsize threshold of $2/\Lmax$ for the upper bounds
in Theorem \ref{thm:RR-cvx-ub} and its application to quadratic-growth
functions in Theorem \ref{thm:RR-QG-ub} were obtained in March 2026
without any assistance from AI, shortly after the submission of the
conference version \citep{pmlr-v336-liu26d}, whose results were likewise
proved without AI assistance. The lower bounds for $\RR$ and $\SGD$
in Theorems \ref{thm:RR-cvx-lb} and \ref{thm:SGD-cvx-lb} were obtained
in July 2026. The former was derived with assistance from GPT-5.6
Sol (particularly in establishing Theorem \ref{thm:RR-cvx-lb-II},
a key intermediate result used in the proof of Theorem \ref{thm:RR-cvx-lb}),
whereas the latter was developed without AI assistance. GPT-5.6 Sol
was also used to polish and refine the presentation during the final
preparation of this extended version in September 2026.

\bibliographystyle{plainnat}
\bibliography{ref}

\clearpage

\appendix

\section{Additional Related Work\label{sec:Additional-Related-Work}}

This section provides additional discussion of the related work. We
mainly focus on smooth convex optimization under $\RR$. As for nonsmooth
convex optimization under $\RR/\SS/\IG$, the reader could refer to
\citet{kibardin1979decomposition,doi:10.1137/S1052623499362111,NEURIPS2022_7bc4f74e,pmlr-v267-liu25ct}.
See also \citet{NIPS2016_c74d97b0} for $\RR$ under structured problems.

The first breakthrough on $\RR$ is by \citet{gurbuzbalaban2021random}
for smooth strongly convex optimization, which, however, requires
each $f_{i}$ to be quadratic or to have a Lipschitz Hessian. Since
then, extensive studies have emerged. Among them, a series of works
continues to study the convergence behavior of $\RR$ in smooth strongly
convex optimization for quadratic objectives or under other additional
conditions \citep{8514028,pmlr-v97-haochen19a,pmlr-v125-safran20a,pmlr-v119-rajput20a,NEURIPS2020_cb8acb1d,NEURIPS2021_803ef568,rajput2022permutationbased}.

To the best of our knowledge, the first work that drops the strong
convexity assumption is \citet{pmlr-v97-nagaraj19a}, which provides
a convergence rate of $\frac{D^{2}}{\eta nK}+\eta G^{2}$ under the
requirement $\eta\leq\frac{2}{L}$ for $L$-smooth $G$-Lipschitz
convex $f_{i}$, where we remind the reader that $D$ denotes the
distance between the initial point and the optimal solution, $\eta$
represents the stepsize, and $K$ is the number of epochs. However,
this rate cannot reflect any advantage of $\RR$ over standard $\SGD$.
Subsequently, two works \citep{NEURIPS2020_c8cc6e90,JMLR:v22:20-1238}
further remove the extra Lipschitz assumption and establish the bound
$\frac{D^{2}}{\eta nK}+\eta^{2}nL\sigma_{\star}^{2}$ under the condition
$\eta\lesssim\frac{1}{nL}$, where $\sigma_{\star}^{2}$ is the gradient
variance at the optimal solution. This rate is faster than standard
$\SGD$ under its required regime and remains the best bound so far.
In fact, it is unimprovable for small $\eta\lesssim\frac{1}{nL}$
and large $K\gtrsim\frac{nL^{2}D^{2}}{\sigma_{\star}^{2}}$ due to
the lower bound of $(\frac{L\sigma_{\star}^{2}D^{4}}{nK^{2}})^{\frac{1}{3}}$
by \citet{pmlr-v202-cha23a}. Recently, \citet{pmlr-v235-liu24cg}
and \citet{ICLR2025_fea9f93f} extend the above rate from the average
iterate to the last iterate (up to additional polylogarithmic factors).

\section{Summary of Notation}

For readability, we recall and summarize the notation used in the
paper.
\begin{itemize}
\item $f=\frac{1}{n}\sum_{i=1}^{n}f_{i}$ is the objective, where each $f_{i}:\R^{d}\to\R$
is differentiable and convex.
\item $\bx_{\star}$ denotes the minimizer of $f$. $f_{\star}=f(\bx_{\star})$
is the optimal function value. $\sigma_{\star}^{2}=\frac{1}{n}\sum_{i=1}^{n}\left\Vert \nabla f_{i}(\bx_{\star})\right\Vert ^{2}$
is the variance of the gradient at the optimal solution.
\item $L_{i}>0$ is the smoothness parameter of $f_{i}$. $\Lavg=\frac{1}{n}\sum_{i=1}^{n}L_{i}$
is the average smoothness parameter. $\Lmax=\max_{i\in\left[n\right]}L_{i}$
is the maximum smoothness parameter.
\item $\Breg_{h}(\bx,\by)=h(\bx)-h(\by)-\left\langle \nabla h(\by),\bx-\by\right\rangle $
denotes the Bregman divergence induced by a real-valued differentiable
function $h$ (not necessarily convex). We write $\Breg$ (resp. $\Breg_{i}$)
as shorthand for $\Breg_{f}$ (resp. $\Breg_{f_{i}}$).
\end{itemize}

\section{Missing Proofs of Core Lemmas\label{sec:Core-Lemmas}}

This section contains the missing proofs of the three core lemmas
presented in Subsection \ref{subsec:Core-Lemmas}.

Before providing the proofs, we recall two notions introduced in (\ref{eq:main-virtual-pi})
and (\ref{eq:main-virtual-x}), respectively. Given a permutation
$\pi$ of $\left[n\right]$ and two indices $i,j\in\left[n\right]$
satisfying $j\leq i$, $\pi(i,j)$ is the permutation generated by
exchanging the elements $\pi_{i}$ and $\pi_{j}$ in $\pi$, i.e.,
\begin{equation}
\pi(i,j)=\left(\pi^{1},\mydots,\pi^{j-1},\pi^{i},\pi^{j+1},\mydots,\pi^{i-1},\pi^{j},\pi^{i+1},\mydots,\pi^{n}\right).\label{eq:virtual-pi}
\end{equation}
For any given $k\in\left[K\right]$, $\bx_{k}^{l}(i,j),\forall l\in\left[n+1\right]$
denotes the trajectory of the $k$-th epoch starting from $\bx_{k}^{1}$,
produced by $\ShufflingSGD$, but under the permutation $\pi_{k}(i,j)$,
i.e.,
\begin{eqnarray}
 & \bx_{k}^{l+1}(i,j)=\bx_{k}^{l}(i,j)-\eta_{k}\nabla f_{\pi_{k}^{l}(i,j)}(\bx_{k}^{l}(i,j)),\forall l\in\left[n\right], & \text{where}\enspace\bx_{k}^{1}(i,j)=\bx_{k}^{1}.\label{eq:virtual-x}
\end{eqnarray}

\subsection{Proof of Lemma \ref{lem:core-function}}

\begin{proof}
Given $k\in\left[K\right]$ and $i\in\left[n\right]$, we first have
the decomposition
\[
\ell(\bx_{k}^{i})-\ell_{\pi_{k}^{i}}(\bx_{k}^{i})=\frac{1}{n}\left(\sum_{j=1}^{i-1}\ell_{\pi_{k}^{j}}(\bx_{k}^{i})-\ell_{\pi_{k}^{i}}(\bx_{k}^{i})+\sum_{j=i}^{n}\ell_{\pi_{k}^{j}}(\bx_{k}^{i})-\ell_{\pi_{k}^{i}}(\bx_{k}^{i})\right).
\]
Note that $\E\left[\ell_{\pi_{k}^{j}}(\bx_{k}^{i})\right]=\E\left[\ell_{\pi_{k}^{i}}(\bx_{k}^{i})\right]$
holds for any $j\in\left\{ i,\mydots,n\right\} $ under $\RR$, since
$\ell_{\pi_{k}^{j}}(\bx_{k}^{i})$ and $\ell_{\pi_{k}^{i}}(\bx_{k}^{i})$
are equal in distribution conditioning on $\pi_{k}^{1}$ to $\pi_{k}^{i-1}$
and $\pi_{1}$ to $\pi_{k-1}$. Therefore, we obtain
\begin{equation}
\E\left[\ell(\bx_{k}^{i})-\ell_{\pi_{k}^{i}}(\bx_{k}^{i})\right]=\frac{1}{n}\sum_{j=1}^{i-1}\E\left[\ell_{\pi_{k}^{j}}(\bx_{k}^{i})-\ell_{\pi_{k}^{i}}(\bx_{k}^{i})\right].\label{eq:core-function-1}
\end{equation}
For any fixed $j\in\left[i\right]$, it is known that $\pi_{k}$ equals
$\pi_{k}(i,j)$ in distribution (e.g., Lemma C.3 of \citet{pmlr-v267-liu25ct}),
which implies that $(\pi_{k},\bx_{k})$ also equals $(\pi_{k}(i,j),\bx_{k}(i,j))$
in distribution, since the trajectory of the $k$-th epoch generated
by $\ShufflingSGD$ is deterministically determined by the permutation
and stepsize. This implies that
\[
\E\left[\ell_{\pi_{k}^{j}}(\bx_{k}^{i})\right]=\E\left[\ell_{\pi_{k}^{j}(i,j)}(\bx_{k}^{i}(i,j))\right]=\E\left[\ell_{\pi_{k}^{i}}(\bx_{k}^{i}(i,j))\right],
\]
where the last step is due to $\pi_{k}^{j}(i,j)=\pi_{k}^{i}$ from
its definition (\ref{eq:main-virtual-pi}). Hence, we finally obtain
\[
\E\left[\ell(\bx_{k}^{i})-\ell_{\pi_{k}^{i}}(\bx_{k}^{i})\right]=\frac{1}{n}\sum_{j=1}^{i-1}\E\left[\ell_{\pi_{k}^{i}}(\bx_{k}^{i}(i,j))-\ell_{\pi_{k}^{i}}(\bx_{k}^{i})\right].
\]
\end{proof}

\subsection{Proof of Lemma \ref{lem:core-stability}}

\begin{proof}
By the definition of $\pi(i,j)$ in (\ref{eq:main-virtual-pi}), we
know
\begin{equation}
\pi_{k}^{l}(i,j)=\pi_{k}^{l},\forall l\notin\left\{ i,j\right\} \Rightarrow\nabla f_{\pi_{k}^{l}(i,j)}=\nabla f_{\pi_{k}^{l}},\forall l\notin\left\{ i,j\right\} .\label{eq:core-stability-1}
\end{equation}
Therefore, given $l\notin\left\{ i,j\right\} $, by the definition
of $\bx_{k}(i,j)$ (see (\ref{eq:main-virtual-x})) and the update
rule of $\ShufflingSGD$, we have
\begin{align}
\left\Vert \bx_{k}^{l+1}(i,j)-\bx_{k}^{l+1}\right\Vert ^{2}= & \left\Vert \bx_{k}^{l}(i,j)-\bx_{k}^{l}-\eta_{k}\left(\nabla f_{\pi_{k}^{l}(i,j)}(\bx_{k}^{l}(i,j))-\nabla f_{\pi_{k}^{l}}(\bx_{k}^{l})\right)\right\Vert ^{2}\nonumber \\
\overset{(\ref{eq:core-stability-1})}{=} & \left\Vert \bx_{k}^{l}(i,j)-\bx_{k}^{l}-\eta_{k}\left(\nabla f_{\pi_{k}^{l}}(\bx_{k}^{l}(i,j))-\nabla f_{\pi_{k}^{l}}(\bx_{k}^{l})\right)\right\Vert ^{2}\nonumber \\
= & \left\Vert \bx_{k}^{l}(i,j)-\bx_{k}^{l}\right\Vert ^{2}+\eta_{k}^{2}\left\Vert \nabla f_{\pi_{k}^{l}}(\bx_{k}^{l}(i,j))-\nabla f_{\pi_{k}^{l}}(\bx_{k}^{l})\right\Vert ^{2}\nonumber \\
 & -2\eta_{k}\left\langle \nabla f_{\pi_{k}^{l}}(\bx_{k}^{l}(i,j))-\nabla f_{\pi_{k}^{l}}(\bx_{k}^{l}),\bx_{k}^{l}(i,j)-\bx_{k}^{l}\right\rangle \nonumber \\
\overset{(a)}{\leq} & \left\Vert \bx_{k}^{l}(i,j)-\bx_{k}^{l}\right\Vert ^{2}+\left(\eta_{k}^{2}-\frac{2\eta_{k}}{L_{\pi_{k}^{l}}}\right)\left\Vert \nabla f_{\pi_{k}^{l}}(\bx_{k}^{l}(i,j))-\nabla f_{\pi_{k}^{l}}(\bx_{k}^{l})\right\Vert ^{2}\nonumber \\
\overset{(b)}{\leq} & \left\Vert \bx_{k}^{l}(i,j)-\bx_{k}^{l}\right\Vert ^{2},\label{eq:core-stability-2}
\end{align}
where $(a)$ is due to Lemma \ref{lem:co-coercivity-cvx} and $(b)$
holds by $\eta_{k}\leq\frac{2}{\Lmax}\Rightarrow\eta_{k}^{2}-\frac{2\eta_{k}}{L_{\pi_{k}^{l}}}\leq0$.

Apply (\ref{eq:core-stability-2}) from $l=j+1$ to $l=i-1$ to obtain
\begin{align*}
 & \left\Vert \bx_{k}^{i}(i,j)-\bx_{k}^{i}\right\Vert ^{2}\leq\left\Vert \bx_{k}^{j+1}(i,j)-\bx_{k}^{j+1}\right\Vert ^{2}\\
= & \left\Vert \bx_{k}^{j}(i,j)-\bx_{k}^{j}-\eta_{k}\left(\nabla f_{\pi_{k}^{j}(i,j)}(\bx_{k}^{j}(i,j))-\nabla f_{\pi_{k}^{j}}(\bx_{k}^{j})\right)\right\Vert ^{2}\\
\overset{(c)}{=} & \eta_{k}^{2}\left\Vert \nabla f_{\pi_{k}^{j}(i,j)}(\bx_{k}^{j})-\nabla f_{\pi_{k}^{j}}(\bx_{k}^{j})\right\Vert ^{2}\overset{(\ref{eq:main-virtual-pi})}{=}\eta_{k}^{2}\left\Vert \nabla f_{\pi_{k}^{i}}(\bx_{k}^{j})-\nabla f_{\pi_{k}^{j}}(\bx_{k}^{j})\right\Vert ^{2},
\end{align*}
where $(c)$ holds by $\bx_{k}^{j}(i,j)=\bx_{k}^{j}$, since $\bx_{k}^{l}(i,j)=\bx_{k}^{l},\forall l\in\left[j\right]$
implied by (\ref{eq:core-stability-2}) and $\bx_{k}^{1}(i,j)\overset{(\ref{eq:main-virtual-x})}{=}\bx_{k}^{1}$
together.
\end{proof}

\subsection{Proof of Lemma \ref{lem:core-Bregman}}

\begin{proof}
Expanding the definitions of $\Breg$ and $\Breg_{\pi_{k}^{i}}$,
we know
\begin{align*}
\Breg(\bx_{k}^{i},\bx_{\star})-\Breg_{\pi_{k}^{i}}(\bx_{k}^{i},\bx_{\star})= & f(\bx_{k}^{i})-f(\bx_{\star})-\left\langle \nabla f(\bx_{\star}),\bx_{k}^{i}-\bx_{\star}\right\rangle \\
 & -\left(f_{\pi_{k}^{i}}(\bx_{k}^{i})-f_{\pi_{k}^{i}}(\bx_{\star})-\left\langle \nabla f_{\pi_{k}^{i}}(\bx_{\star}),\bx_{k}^{i}-\bx_{\star}\right\rangle \right)\\
= & f(\bx_{k}^{i})-f_{\pi_{k}^{i}}(\bx_{k}^{i})-f(\bx_{\star})+f_{\pi_{k}^{i}}(\bx_{\star})\\
 & -\left\langle \nabla f(\bx_{\star})-\nabla f_{\pi_{k}^{i}}(\bx_{\star}),\bx_{k}^{i}-\bx_{\star}\right\rangle .
\end{align*}
Since $\E\left[f_{\pi_{k}^{i}}(\bx_{\star})\right]=f(\bx_{\star})$
and $\E\left[\nabla f_{\pi_{k}^{i}}(\bx_{\star})\right]=\nabla f(\bx_{\star})$
under $\RR$, after taking expectations on both sides, we obtain that
\begin{equation}
\E\left[\Breg(\bx_{k}^{i},\bx_{\star})-\Breg_{\pi_{k}^{i}}(\bx_{k}^{i},\bx_{\star})\right]=\E\left[f(\bx_{k}^{i})-f_{\pi_{k}^{i}}(\bx_{k}^{i})-\left\langle \nabla f(\bx_{\star})-\nabla f_{\pi_{k}^{i}}(\bx_{\star}),\bx_{k}^{i}\right\rangle \right].\label{eq:core-Bregman-1}
\end{equation}
Now, we denote by $\ell_{i}(\bx)\defeq f_{i}(\bx)-\left\langle \nabla f_{i}(\bx_{\star}),\bx\right\rangle ,\forall i\in\left[n\right]$
and $\ell(\bx)\defeq\frac{1}{n}\sum_{i=1}^{n}\ell_{i}(\bx)=f(\bx)-\left\langle \nabla f(\bx_{\star}),\bx\right\rangle $.
Then, (\ref{eq:core-Bregman-1}) implies that
\begin{equation}
\E\left[\Breg(\bx_{k}^{i},\bx_{\star})-\Breg_{\pi_{k}^{i}}(\bx_{k}^{i},\bx_{\star})\right]=\E\left[\ell(\bx_{k}^{i})-\ell_{\pi_{k}^{i}}(\bx_{k}^{i})\right]\overset{(a)}{=}\frac{1}{n}\sum_{j=1}^{i-1}\E\left[\ell_{\pi_{k}^{i}}(\bx_{k}^{i}(i,j))-\ell_{\pi_{k}^{i}}(\bx_{k}^{i})\right],\label{eq:core-Bregman-2}
\end{equation}
where $(a)$ holds by Lemma \ref{lem:core-function}.

Note that $\ell_{\pi_{k}^{i}}$ is $L_{\pi_{k}^{i}}$-smooth by its
definition and Assumption \ref{assu:smoothness}, we therefore have,
for any $j\in\left[i-1\right]$,
\begin{align}
\ell_{\pi_{k}^{i}}(\bx_{k}^{i}(i,j))-\ell_{\pi_{k}^{i}}(\bx_{k}^{i}) & \leq\left\langle \nabla\ell_{\pi_{k}^{i}}(\bx_{k}^{i}),\bx_{k}^{i}(i,j)-\bx_{k}^{i}\right\rangle +\frac{L_{\pi_{k}^{i}}}{2}\left\Vert \bx_{k}^{i}(i,j)-\bx_{k}^{i}\right\Vert ^{2}\nonumber \\
 & =\left\langle \nabla f_{\pi_{k}^{i}}(\bx_{k}^{i})-\nabla f_{\pi_{k}^{i}}(\bx_{\star}),\bx_{k}^{i}(i,j)-\bx_{k}^{i}\right\rangle +\frac{L_{\pi_{k}^{i}}}{2}\left\Vert \bx_{k}^{i}(i,j)-\bx_{k}^{i}\right\Vert ^{2}\nonumber \\
 & \overset{(b)}{\leq}\frac{\left\Vert \nabla f_{\pi_{k}^{i}}(\bx_{k}^{i})-\nabla f_{\pi_{k}^{i}}(\bx_{\star})\right\Vert ^{2}}{2L_{\pi_{k}^{i}}}+L_{\pi_{k}^{i}}\left\Vert \bx_{k}^{i}(i,j)-\bx_{k}^{i}\right\Vert ^{2}\nonumber \\
 & \overset{(c)}{\leq}\Breg_{\pi_{k}^{i}}(\bx_{k}^{i},\bx_{\star})+L_{\pi_{k}^{i}}\left\Vert \bx_{k}^{i}(i,j)-\bx_{k}^{i}\right\Vert ^{2}\nonumber \\
 & \overset{(d)}{\leq}\Breg_{\pi_{k}^{i}}(\bx_{k}^{i},\bx_{\star})+\eta_{k}^{2}L_{\pi_{k}^{i}}\left\Vert \nabla f_{\pi_{k}^{i}}(\bx_{k}^{j})-\nabla f_{\pi_{k}^{j}}(\bx_{k}^{j})\right\Vert ^{2},\label{eq:core-Bregman-3}
\end{align}
where $(b)$ is by Cauchy-Schwarz inequality and AM-GM inequality,
$(c)$ is due to Lemma \ref{lem:co-coercivity-cvx}, and $(d)$ holds
by Lemma \ref{lem:core-stability}. Furthermore, we can bound
\begin{align*}
 & \left\Vert \nabla f_{\pi_{k}^{i}}(\bx_{k}^{j})-\nabla f_{\pi_{k}^{j}}(\bx_{k}^{j})\right\Vert ^{2}\\
\leq & 4\left\Vert \nabla f_{\pi_{k}^{i}}(\bx_{k}^{j})-\nabla f_{\pi_{k}^{i}}(\bx_{\star})\right\Vert ^{2}+4\left\Vert \nabla f_{\pi_{k}^{j}}(\bx_{k}^{j})-\nabla f_{\pi_{k}^{j}}(\bx_{\star})\right\Vert ^{2}+4\left\Vert \nabla f_{\pi_{k}^{i}}(\bx_{\star})\right\Vert ^{2}+4\left\Vert \nabla f_{\pi_{k}^{j}}(\bx_{\star})\right\Vert ^{2}\\
\leq & 8L_{\pi_{k}^{i}}\Breg_{\pi_{k}^{i}}(\bx_{k}^{j},\bx_{\star})+8L_{\pi_{k}^{j}}\Breg_{\pi_{k}^{j}}(\bx_{k}^{j},\bx_{\star})+4\left\Vert \nabla f_{\pi_{k}^{i}}(\bx_{\star})\right\Vert ^{2}+4\left\Vert \nabla f_{\pi_{k}^{j}}(\bx_{\star})\right\Vert ^{2},
\end{align*}
where the last step is by Lemma \ref{lem:co-coercivity-cvx} again.
Plug the above inequality back into (\ref{eq:core-Bregman-3}) and
use $\eta_{k}\le\frac{2}{\Lmax}$ to have
\begin{align}
\ell_{\pi_{k}^{i}}(\bx_{k}^{i}(i,j))-\ell_{\pi_{k}^{i}}(\bx_{k}^{i})\leq & \Breg_{\pi_{k}^{i}}(\bx_{k}^{i},\bx_{\star})+32\Breg_{\pi_{k}^{i}}(\bx_{k}^{j},\bx_{\star})+32\Breg_{\pi_{k}^{j}}(\bx_{k}^{j},\bx_{\star})\nonumber \\
 & +4\eta_{k}^{2}L_{\pi_{k}^{i}}\left\Vert \nabla f_{\pi_{k}^{i}}(\bx_{\star})\right\Vert ^{2}+4\eta_{k}^{2}L_{\pi_{k}^{i}}\left\Vert \nabla f_{\pi_{k}^{j}}(\bx_{\star})\right\Vert ^{2}.\label{eq:core-Bregman-4}
\end{align}

Combine (\ref{eq:core-Bregman-2}) and (\ref{eq:core-Bregman-4})
to obtain
\begin{align*}
\E\left[\Breg(\bx_{k}^{i},\bx_{\star})-\Breg_{\pi_{k}^{i}}(\bx_{k}^{i},\bx_{\star})\right]\leq & \frac{1}{n}\sum_{j=1}^{i-1}\E\left[\Breg_{\pi_{k}^{i}}(\bx_{k}^{i},\bx_{\star})+32\Breg_{\pi_{k}^{i}}(\bx_{k}^{j},\bx_{\star})+32\Breg_{\pi_{k}^{j}}(\bx_{k}^{j},\bx_{\star})\right]\\
 & +\frac{4\eta_{k}^{2}}{n}\sum_{j=1}^{i-1}\E\left[L_{\pi_{k}^{i}}\left\Vert \nabla f_{\pi_{k}^{i}}(\bx_{\star})\right\Vert ^{2}+L_{\pi_{k}^{i}}\left\Vert \nabla f_{\pi_{k}^{j}}(\bx_{\star})\right\Vert ^{2}\right].
\end{align*}
When $j<i$, we observe that $\E\left[\Breg_{\pi_{k}^{i}}(\bx_{k}^{j},\bx_{\star})\right]=\E\left[\Breg_{\pi_{k}^{j}}(\bx_{k}^{j},\bx_{\star})\right]$
and the following two inequalities hold
\begin{align*}
\E\left[L_{\pi_{k}^{i}}\left\Vert \nabla f_{\pi_{k}^{i}}(\bx_{\star})\right\Vert ^{2}\right] & =\frac{\sum_{l=1}^{n}L_{l}\left\Vert \nabla f_{l}(\bx_{\star})\right\Vert ^{2}}{n}\leq\Lmax\sigma_{\star}^{2},\\
\E\left[L_{\pi_{k}^{i}}\left\Vert \nabla f_{\pi_{k}^{j}}(\bx_{\star})\right\Vert ^{2}\right] & =\frac{\sum_{l=1}^{n}\frac{n\Lavg-L_{l}}{n-1}\left\Vert \nabla f_{l}(\bx_{\star})\right\Vert ^{2}}{n}\leq2\Lavg\sigma_{\star}^{2}\leq2\Lmax\sigma_{\star}^{2}.
\end{align*}
Put everything together, rearrange terms, and use $i-1\leq n$ to
conclude the desired inequality.
\end{proof}

\section{Missing Proofs of Lemmas for Bound I\label{sec:RR-cvx-ub-I}}

In this section, we provide the missing proofs of the lemmas presented
in Subsection \ref{subsec:RR-cvx-ub-I}, which were used to prove
the first convergence rate in Theorem \ref{thm:RR-cvx-ub-I}.

\subsection{Proof of Lemma \ref{lem:RR-cvx-ub-I-descent}}

\begin{proof}
Given $k\in\left[K\right]$ and $i\in\left[n\right]$, by the update
rule of Algorithm \ref{alg:Shuffling-SGD}, we have
\begin{align}
 & \left\Vert \bx_{k}^{i+1}-\bx_{\star}\right\Vert ^{2}=\left\Vert \bx_{k}^{i}-\bx_{\star}-\eta_{k}\nabla f_{\pi_{k}^{i}}(\bx_{k}^{i})\right\Vert ^{2}\nonumber \\
= & \left\Vert \bx_{k}^{i}-\bx_{\star}\right\Vert ^{2}+\eta_{k}^{2}\left\Vert \nabla f_{\pi_{k}^{i}}(\bx_{k}^{i})\right\Vert ^{2}-2\eta_{k}\left\langle \nabla f_{\pi_{k}^{i}}(\bx_{k}^{i}),\bx_{k}^{i}-\bx_{\star}\right\rangle .\label{eq:RR-cvx-ub-I-descent-1}
\end{align}
Note that $\eta_{k}<\frac{2}{\Lmax}\Rightarrow2-\eta_{k}L_{\pi_{k}^{i}}>0$,
we hence can bound
\begin{align}
\left\Vert \nabla f_{\pi_{k}^{i}}(\bx_{k}^{i})\right\Vert ^{2} & \leq\left(1+\frac{2-\eta_{k}L_{\pi_{k}^{i}}}{2\eta_{k}L_{\pi_{k}^{i}}}\right)\left\Vert \nabla f_{\pi_{k}^{i}}(\bx_{k}^{i})-\nabla f_{\pi_{k}^{i}}(\bx_{\star})\right\Vert ^{2}+\left(1+\frac{2\eta_{k}L_{\pi_{k}^{i}}}{2-\eta_{k}L_{\pi_{k}^{i}}}\right)\left\Vert \nabla f_{\pi_{k}^{i}}(\bx_{\star})\right\Vert ^{2}\nonumber \\
 & =\frac{2+\eta_{k}L_{\pi_{k}^{i}}}{2\eta_{k}L_{\pi_{k}^{i}}}\left\Vert \nabla f_{\pi_{k}^{i}}(\bx_{k}^{i})-\nabla f_{\pi_{k}^{i}}(\bx_{\star})\right\Vert ^{2}+\frac{2+\eta_{k}L_{\pi_{k}^{i}}}{2-\eta_{k}L_{\pi_{k}^{i}}}\left\Vert \nabla f_{\pi_{k}^{i}}(\bx_{\star})\right\Vert ^{2}\nonumber \\
 & \leq\frac{2+\eta_{k}L_{\pi_{k}^{i}}}{2\eta_{k}}\left\langle \nabla f_{\pi_{k}^{i}}(\bx_{k}^{i})-\nabla f_{\pi_{k}^{i}}(\bx_{\star}),\bx_{k}^{i}-\bx_{\star}\right\rangle +\frac{2+\eta_{k}L_{\pi_{k}^{i}}}{2-\eta_{k}L_{\pi_{k}^{i}}}\left\Vert \nabla f_{\pi_{k}^{i}}(\bx_{\star})\right\Vert ^{2},\label{eq:RR-cvx-ub-I-descent-2}
\end{align}
where the last step is due to Lemma \ref{lem:co-coercivity-cvx}.

We plug (\ref{eq:RR-cvx-ub-I-descent-2}) back into (\ref{eq:RR-cvx-ub-I-descent-1})
to obtain that
\begin{align*}
\left\Vert \bx_{k}^{i+1}-\bx_{\star}\right\Vert ^{2}\leq & \left\Vert \bx_{k}^{i}-\bx_{\star}\right\Vert ^{2}-\frac{2-\eta_{k}L_{\pi_{k}^{i}}}{2}\eta_{k}\left\langle \nabla f_{\pi_{k}^{i}}(\bx_{k}^{i})-\nabla f_{\pi_{k}^{i}}(\bx_{\star}),\bx_{k}^{i}-\bx_{\star}\right\rangle \\
 & -2\eta_{k}\left\langle \nabla f_{\pi_{k}^{i}}(\bx_{\star}),\bx_{k}^{i}-\bx_{\star}\right\rangle +\frac{2+\eta_{k}L_{\pi_{k}^{i}}}{2-\eta_{k}L_{\pi_{k}^{i}}}\eta_{k}^{2}\left\Vert \nabla f_{\pi_{k}^{i}}(\bx_{\star})\right\Vert ^{2}\\
= & \left\Vert \bx_{k}^{i}-\bx_{\star}\right\Vert ^{2}-\frac{2-\eta_{k}L_{\pi_{k}^{i}}}{2}\eta_{k}\left(\Breg_{\pi_{k}^{i}}(\bx_{k}^{i},\bx_{\star})+\Breg_{\pi_{k}^{i}}(\bx_{\star},\bx_{k}^{i})\right)\\
 & -2\eta_{k}\left\langle \nabla f_{\pi_{k}^{i}}(\bx_{\star}),\bx_{k}^{i}-\bx_{\star}\right\rangle +\frac{2+\eta_{k}L_{\pi_{k}^{i}}}{2-\eta_{k}L_{\pi_{k}^{i}}}\eta_{k}^{2}\left\Vert \nabla f_{\pi_{k}^{i}}(\bx_{\star})\right\Vert ^{2}.
\end{align*}
After rearranging terms, we have
\begin{align*}
\frac{2-\eta_{k}L_{\pi_{k}^{i}}}{2}\eta_{k}\Breg_{\pi_{k}^{i}}(\bx_{k}^{i},\bx_{\star})\leq & \left\Vert \bx_{k}^{i}-\bx_{\star}\right\Vert ^{2}-\left\Vert \bx_{k}^{i+1}-\bx_{\star}\right\Vert ^{2}+\frac{2+\eta_{k}L_{\pi_{k}^{i}}}{2-\eta_{k}L_{\pi_{k}^{i}}}\eta_{k}^{2}\left\Vert \nabla f_{\pi_{k}^{i}}(\bx_{\star})\right\Vert ^{2}\\
 & -2\eta_{k}\left\langle \nabla f_{\pi_{k}^{i}}(\bx_{\star}),\bx_{k}^{i}-\bx_{\star}\right\rangle -\frac{2-\eta_{k}L_{\pi_{k}^{i}}}{2}\eta_{k}\Breg_{\pi_{k}^{i}}(\bx_{\star},\bx_{k}^{i}).
\end{align*}
Therefore, we conclude the proof after applying
\begin{eqnarray*}
2-\eta_{k}\Lmax\leq2-\eta_{k}L_{\pi_{k}^{i}} & \text{and} & \frac{2+\eta_{k}L_{\pi_{k}^{i}}}{2-\eta_{k}L_{\pi_{k}^{i}}}\leq\frac{2+\eta_{k}\Lmax}{2-\eta_{k}\Lmax},
\end{eqnarray*}
both of which hold by $L_{\pi_{k}^{i}}\leq\Lmax$.
\end{proof}

\subsection{Proof of Lemma \ref{lem:RR-cvx-ub-I-core}}

\begin{proof}
Note that $\nabla f(\bx_{\star})=\bzero$ now, which implies that
\[
-\left\langle \nabla f_{\pi_{k}^{i}}(\bx_{\star}),\bx_{k}^{i}-\bx_{\star}\right\rangle =\left\langle \nabla f(\bx_{\star})-\nabla f_{\pi_{k}^{i}}(\bx_{\star}),\bx_{k}^{i}-\bx_{\star}\right\rangle .
\]
Since $\E\left[\nabla f_{\pi_{k}^{i}}(\bx_{\star})\right]=\nabla f(\bx_{\star})$
under $\RR$, after taking expectations on both sides, we obtain that
\begin{equation}
\E\left[-\left\langle \nabla f_{\pi_{k}^{i}}(\bx_{\star}),\bx_{k}^{i}-\bx_{\star}\right\rangle \right]=\E\left[\left\langle \nabla f(\bx_{\star})-\nabla f_{\pi_{k}^{i}}(\bx_{\star}),\bx_{k}^{i}\right\rangle \right].\label{eq:RR-cvx-ub-I-core-1}
\end{equation}
Now, we denote by $\ell_{i}(\bx)\defeq\left\langle \nabla f_{i}(\bx_{\star}),\bx\right\rangle ,\forall i\in\left[n\right]$
and $\ell(\bx)\defeq\frac{1}{n}\sum_{i=1}^{n}\ell_{i}(\bx)=\left\langle \nabla f(\bx_{\star}),\bx\right\rangle $.
Then, (\ref{eq:RR-cvx-ub-I-core-1}) implies that
\begin{equation}
\E\left[-\left\langle \nabla f_{\pi_{k}^{i}}(\bx_{\star}),\bx_{k}^{i}-\bx_{\star}\right\rangle \right]=\E\left[\ell(\bx_{k}^{i})-\ell_{\pi_{k}^{i}}(\bx_{k}^{i})\right]\overset{(a)}{=}\frac{1}{n}\sum_{j=1}^{i-1}\E\left[\ell_{\pi_{k}^{i}}(\bx_{k}^{i}(i,j))-\ell_{\pi_{k}^{i}}(\bx_{k}^{i})\right],\label{eq:RR-cvx-ub-I-core-2}
\end{equation}
where $(a)$ holds by Lemma \ref{lem:core-function}.

Note that, for any $j<i$, we have
\begin{align*}
\ell_{\pi_{k}^{i}}(\bx_{k}^{i}(i,j))-\ell_{\pi_{k}^{i}}(\bx_{k}^{i}) & =\left\langle \nabla f_{\pi_{k}^{i}}(\bx_{\star}),\bx_{k}^{i}(i,j)-\bx_{k}^{i}\right\rangle \leq\left\Vert \nabla f_{\pi_{k}^{i}}(\bx_{\star})\right\Vert \left\Vert \bx_{k}^{i}(i,j)-\bx_{k}^{i}\right\Vert \\
 & \overset{(b)}{\leq}\eta_{k}\left\Vert \nabla f_{\pi_{k}^{i}}(\bx_{\star})\right\Vert \left\Vert \nabla f_{\pi_{k}^{i}}(\bx_{k}^{j})-\nabla f_{\pi_{k}^{j}}(\bx_{k}^{j})\right\Vert ,
\end{align*}
where $(b)$ follows from Lemma \ref{lem:core-stability}. Furthermore,
we can bound
\[
\left\Vert \nabla f_{\pi_{k}^{i}}(\bx_{k}^{j})-\nabla f_{\pi_{k}^{j}}(\bx_{k}^{j})\right\Vert \leq\left\Vert \nabla f_{\pi_{k}^{i}}(\bx_{k}^{j})-\nabla f_{\pi_{k}^{i}}(\bx_{\star})\right\Vert +\left\Vert \nabla f_{\pi_{k}^{j}}(\bx_{k}^{j})-\nabla f_{\pi_{k}^{j}}(\bx_{\star})\right\Vert +\left\Vert \nabla f_{\pi_{k}^{i}}(\bx_{\star})\right\Vert +\left\Vert \nabla f_{\pi_{k}^{j}}(\bx_{\star})\right\Vert ,
\]
which implies that
\begin{align}
 & \ell_{\pi_{k}^{i}}(\bx_{k}^{i}(i,j))-\ell_{\pi_{k}^{i}}(\bx_{k}^{i})\nonumber \\
\leq & \eta_{k}\left\Vert \nabla f_{\pi_{k}^{i}}(\bx_{\star})\right\Vert \left\Vert \nabla f_{\pi_{k}^{i}}(\bx_{k}^{j})-\nabla f_{\pi_{k}^{i}}(\bx_{\star})\right\Vert +\eta_{k}\left\Vert \nabla f_{\pi_{k}^{i}}(\bx_{\star})\right\Vert \left\Vert \nabla f_{\pi_{k}^{j}}(\bx_{k}^{j})-\nabla f_{\pi_{k}^{j}}(\bx_{\star})\right\Vert \nonumber \\
 & +\eta_{k}\left\Vert \nabla f_{\pi_{k}^{i}}(\bx_{\star})\right\Vert ^{2}+\eta_{k}\left\Vert \nabla f_{\pi_{k}^{i}}(\bx_{\star})\right\Vert \left\Vert \nabla f_{\pi_{k}^{j}}(\bx_{\star})\right\Vert .\label{eq:RR-cvx-ub-I-core-3}
\end{align}
By AM-GM inequality, the following three inequalities hold
\begin{align*}
\eta_{k}\left\Vert \nabla f_{\pi_{k}^{i}}(\bx_{\star})\right\Vert \left\Vert \nabla f_{\pi_{k}^{i}}(\bx_{k}^{j})-\nabla f_{\pi_{k}^{i}}(\bx_{\star})\right\Vert  & \leq\frac{4\eta_{k}^{2}\Lmax}{2-\eta_{k}\Lmax}\left\Vert \nabla f_{\pi_{k}^{i}}(\bx_{\star})\right\Vert ^{2}+\frac{2-\eta_{k}\Lmax}{16\Lmax}\left\Vert \nabla f_{\pi_{k}^{i}}(\bx_{k}^{j})-\nabla f_{\pi_{k}^{i}}(\bx_{\star})\right\Vert ^{2},\\
\eta_{k}\left\Vert \nabla f_{\pi_{k}^{i}}(\bx_{\star})\right\Vert \left\Vert \nabla f_{\pi_{k}^{j}}(\bx_{k}^{j})-\nabla f_{\pi_{k}^{j}}(\bx_{\star})\right\Vert  & \leq\frac{4\eta_{k}^{2}\Lmax}{2-\eta_{k}\Lmax}\left\Vert \nabla f_{\pi_{k}^{i}}(\bx_{\star})\right\Vert ^{2}+\frac{2-\eta_{k}\Lmax}{16\Lmax}\left\Vert \nabla f_{\pi_{k}^{j}}(\bx_{k}^{j})-\nabla f_{\pi_{k}^{j}}(\bx_{\star})\right\Vert ^{2},\\
\eta_{k}\left\Vert \nabla f_{\pi_{k}^{i}}(\bx_{\star})\right\Vert \left\Vert \nabla f_{\pi_{k}^{j}}(\bx_{\star})\right\Vert  & \leq\frac{\eta_{k}}{2}\left\Vert \nabla f_{\pi_{k}^{i}}(\bx_{\star})\right\Vert ^{2}+\frac{\eta_{k}}{2}\left\Vert \nabla f_{\pi_{k}^{j}}(\bx_{\star})\right\Vert ^{2}.
\end{align*}
Plug the above inequalities back into (\ref{eq:RR-cvx-ub-I-core-3})
to have
\begin{align}
\ell_{\pi_{k}^{i}}(\bx_{k}^{i}(i,j))-\ell_{\pi_{k}^{i}}(\bx_{k}^{i})\leq & \frac{2-\eta_{k}\Lmax}{16\Lmax}\left\Vert \nabla f_{\pi_{k}^{i}}(\bx_{k}^{j})-\nabla f_{\pi_{k}^{i}}(\bx_{\star})\right\Vert ^{2}+\frac{2-\eta_{k}\Lmax}{16\Lmax}\left\Vert \nabla f_{\pi_{k}^{j}}(\bx_{k}^{j})-\nabla f_{\pi_{k}^{j}}(\bx_{\star})\right\Vert ^{2}\nonumber \\
 & +\left(\frac{8\eta_{k}\Lmax}{2-\eta_{k}\Lmax}+\frac{3}{2}\right)\eta_{k}\left\Vert \nabla f_{\pi_{k}^{i}}(\bx_{\star})\right\Vert ^{2}+\frac{\eta_{k}}{2}\left\Vert \nabla f_{\pi_{k}^{j}}(\bx_{\star})\right\Vert ^{2}\nonumber \\
\leq & \frac{2-\eta_{k}\Lmax}{8}\Breg_{\pi_{k}^{i}}(\bx_{\star},\bx_{k}^{j})+\frac{2-\eta_{k}\Lmax}{8}\Breg_{\pi_{k}^{j}}(\bx_{\star},\bx_{k}^{j})\nonumber \\
 & +\left(\frac{8\eta_{k}\Lmax}{2-\eta_{k}\Lmax}+\frac{3}{2}\right)\eta_{k}\left\Vert \nabla f_{\pi_{k}^{i}}(\bx_{\star})\right\Vert ^{2}+\frac{\eta_{k}}{2}\left\Vert \nabla f_{\pi_{k}^{j}}(\bx_{\star})\right\Vert ^{2},\label{eq:RR-cvx-ub-I-core-4}
\end{align}
where the last step is by Lemma \ref{lem:co-coercivity-cvx}, $L_{\pi_{k}^{i}}\leq\Lmax$
and $L_{\pi_{k}^{j}}\leq\Lmax$.

Combine (\ref{eq:RR-cvx-ub-I-core-2}) and (\ref{eq:RR-cvx-ub-I-core-4})
to obtain that
\begin{align*}
\E\left[-\left\langle \nabla f_{\pi_{k}^{i}}(\bx_{\star}),\bx_{k}^{i}-\bx_{\star}\right\rangle \right]\leq & \frac{1}{n}\sum_{j=1}^{i-1}\E\left[\frac{2-\eta_{k}\Lmax}{8}\Breg_{\pi_{k}^{i}}(\bx_{\star},\bx_{k}^{j})+\frac{2-\eta_{k}\Lmax}{8}\Breg_{\pi_{k}^{j}}(\bx_{\star},\bx_{k}^{j})\right]\\
 & +\frac{1}{n}\sum_{j=1}^{i-1}\E\left[\left(\frac{8\eta_{k}\Lmax}{2-\eta_{k}\Lmax}+\frac{3}{2}\right)\eta_{k}\left\Vert \nabla f_{\pi_{k}^{i}}(\bx_{\star})\right\Vert ^{2}+\frac{\eta_{k}}{2}\left\Vert \nabla f_{\pi_{k}^{j}}(\bx_{\star})\right\Vert ^{2}\right].
\end{align*}
We observe that $\E\left[\Breg_{\pi_{k}^{i}}(\bx_{\star},\bx_{k}^{j})\right]=\E\left[\Breg_{\pi_{k}^{j}}(\bx_{\star},\bx_{k}^{j})\right],\forall j<i$
and $\E\left[\left\Vert \nabla f_{\pi_{k}^{i}}(\bx_{\star})\right\Vert ^{2}\right]=\E\left[\left\Vert \nabla f_{\pi_{k}^{j}}(\bx_{\star})\right\Vert ^{2}\right]=\sigma_{\star}^{2}$.
Put everything together to conclude the desired inequality.
\end{proof}

\section{Missing Proofs of Lemmas for Bound II\label{sec:Bound-II}}

In this section, we provide the missing proofs of the lemmas presented
in Subsection \ref{subsec:RR-cvx-ub-II}, which were used to prove
the second convergence rate in Theorem \ref{thm:RR-cvx-ub-II}. 

Before presenting the proofs, we recall the notion introduced in (\ref{eq:main-virtual-y-def}).
For any $k\in\left[K\right]$, the virtual sequence $\by_{k}^{i},\forall i\in\left[n+1\right]$
follows the equations
\begin{eqnarray}
\by_{k}^{i+1}=\by_{k}^{i}-\eta_{k}\nabla f_{\pi_{k}^{i}}(\bx_{\star}),\forall i\in\left[n\right], & \by_{k+1}^{1}=\by_{k}^{n+1}, & \text{where}\enspace\by_{1}^{1}=\bx_{\star}.\label{eq:virtual-y-def}
\end{eqnarray}
As shown in (\ref{eq:main-virtual-y-prop}), the virtual sequence
satisfies that
\begin{equation}
\by_{k}^{n+1}=\by_{k}^{1}=\bx_{\star},\forall k\in\left[K\right].\label{eq:virtual-y-prop}
\end{equation}

\subsection{Proof of Lemma \ref{lem:RR-cvx-ub-II-descent}}

\begin{proof}
Given $k\in\left[K\right]$ and $i\in\left[n\right]$, by the update
rule of Algorithm \ref{alg:Shuffling-SGD} and the definition of $\by_{k}^{i+1}$
(see (\ref{eq:virtual-y-def})), we have
\begin{align}
 & \left\Vert \bx_{k}^{i+1}-\by_{k}^{i+1}\right\Vert ^{2}=\left\Vert \bx_{k}^{i}-\by_{k}^{i}-\eta_{k}\left(\nabla f_{\pi_{k}^{i}}(\bx_{k}^{i})-\nabla f_{\pi_{k}^{i}}(\bx_{\star})\right)\right\Vert ^{2}\nonumber \\
= & \left\Vert \bx_{k}^{i}-\by_{k}^{i}\right\Vert ^{2}+\eta_{k}^{2}\left\Vert \nabla f_{\pi_{k}^{i}}(\bx_{k}^{i})-\nabla f_{\pi_{k}^{i}}(\bx_{\star})\right\Vert ^{2}+2\eta_{k}\left\langle \nabla f_{\pi_{k}^{i}}(\bx_{\star})-\nabla f_{\pi_{k}^{i}}(\bx_{k}^{i}),\bx_{k}^{i}-\by_{k}^{i}\right\rangle .\label{eq:RR-cvx-ub-II-descent-1}
\end{align}
Note that, by Lemma \ref{lem:co-coercivity-cvx}, we can find that
\begin{align*}
 & \left\Vert \nabla f_{\pi_{k}^{i}}(\bx_{k}^{i})-\nabla f_{\pi_{k}^{i}}(\bx_{\star})\right\Vert ^{2}\leq L_{\pi_{k}^{i}}\left\langle \nabla f_{\pi_{k}^{i}}(\bx_{k}^{i})-\nabla f_{\pi_{k}^{i}}(\bx_{\star}),\bx_{k}^{i}-\bx_{\star}\right\rangle \\
= & L_{\pi_{k}^{i}}\left\langle \nabla f_{\pi_{k}^{i}}(\bx_{k}^{i})-\nabla f_{\pi_{k}^{i}}(\bx_{\star}),\by_{k}^{i}-\bx_{\star}\right\rangle +L_{\pi_{k}^{i}}\left\langle \nabla f_{\pi_{k}^{i}}(\bx_{k}^{i})-\nabla f_{\pi_{k}^{i}}(\bx_{\star}),\bx_{k}^{i}-\by_{k}^{i}\right\rangle \\
\leq & \frac{2-\eta_{k}L_{\pi_{k}^{i}}}{4}\left\Vert \nabla f_{\pi_{k}^{i}}(\bx_{k}^{i})-\nabla f_{\pi_{k}^{i}}(\bx_{\star})\right\Vert ^{2}+\frac{L_{\pi_{k}^{i}}^{2}}{2-\eta_{k}L_{\pi_{k}^{i}}}\left\Vert \by_{k}^{i}-\bx_{\star}\right\Vert ^{2}+L_{\pi_{k}^{i}}\left\langle \nabla f_{\pi_{k}^{i}}(\bx_{k}^{i})-\nabla f_{\pi_{k}^{i}}(\bx_{\star}),\bx_{k}^{i}-\by_{k}^{i}\right\rangle ,
\end{align*}
where the last step is due to Cauchy-Schwarz inequality, AM-GM inequality,
and $\eta_{k}<\frac{2}{\Lmax}\Rightarrow2-\eta_{k}L_{\pi_{k}^{i}}>0$.
Rearrange terms in the above inequality to obtain that
\begin{equation}
\left\Vert \nabla f_{\pi_{k}^{i}}(\bx_{k}^{i})-\nabla f_{\pi_{k}^{i}}(\bx_{\star})\right\Vert ^{2}\leq\frac{4L_{\pi_{k}^{i}}^{2}}{4-\eta_{k}^{2}L_{\pi_{k}^{i}}^{2}}\left\Vert \by_{k}^{i}-\bx_{\star}\right\Vert ^{2}+\frac{4L_{\pi_{k}^{i}}}{2+\eta_{k}L_{\pi_{k}^{i}}}\left\langle \nabla f_{\pi_{k}^{i}}(\bx_{k}^{i})-\nabla f_{\pi_{k}^{i}}(\bx_{\star}),\bx_{k}^{i}-\by_{k}^{i}\right\rangle .\label{eq:RR-cvx-ub-II-descent-2}
\end{equation}

Next, we plug (\ref{eq:RR-cvx-ub-II-descent-2}) back into (\ref{eq:RR-cvx-ub-II-descent-1})
to have
\begin{align}
\left\Vert \bx_{k}^{i+1}-\by_{k}^{i+1}\right\Vert ^{2} & \leq\left\Vert \bx_{k}^{i}-\by_{k}^{i}\right\Vert ^{2}+2\eta_{k}\frac{2-\eta_{k}L_{\pi_{k}^{i}}}{2+\eta_{k}L_{\pi_{k}^{i}}}\left\langle \nabla f_{\pi_{k}^{i}}(\bx_{\star})-\nabla f_{\pi_{k}^{i}}(\bx_{k}^{i}),\bx_{k}^{i}-\by_{k}^{i}\right\rangle +\frac{4\eta_{k}^{2}L_{\pi_{k}^{i}}^{2}}{4-\eta_{k}^{2}L_{\pi_{k}^{i}}^{2}}\left\Vert \by_{k}^{i}-\bx_{\star}\right\Vert ^{2}\nonumber \\
 & =\left\Vert \bx_{k}^{i}-\by_{k}^{i}\right\Vert ^{2}+2\eta_{k}\frac{2-\eta_{k}L_{\pi_{k}^{i}}}{2+\eta_{k}L_{\pi_{k}^{i}}}\left(\Breg_{\pi_{k}^{i}}(\by_{k}^{i},\bx_{\star})-\Breg_{\pi_{k}^{i}}(\by_{k}^{i},\bx_{k}^{i})-\Breg_{\pi_{k}^{i}}(\bx_{k}^{i},\bx_{\star})\right)+\frac{4\eta_{k}^{2}L_{\pi_{k}^{i}}^{2}}{4-\eta_{k}^{2}L_{\pi_{k}^{i}}^{2}}\left\Vert \by_{k}^{i}-\bx_{\star}\right\Vert ^{2}\nonumber \\
 & \overset{(a)}{\leq}\left\Vert \bx_{k}^{i}-\by_{k}^{i}\right\Vert ^{2}+2\eta_{k}\frac{2-\eta_{k}L_{\pi_{k}^{i}}}{2+\eta_{k}L_{\pi_{k}^{i}}}\left(\Breg_{\pi_{k}^{i}}(\by_{k}^{i},\bx_{\star})-\Breg_{\pi_{k}^{i}}(\bx_{k}^{i},\bx_{\star})\right)+\frac{4\eta_{k}^{2}L_{\pi_{k}^{i}}^{2}}{4-\eta_{k}^{2}L_{\pi_{k}^{i}}^{2}}\left\Vert \by_{k}^{i}-\bx_{\star}\right\Vert ^{2}\nonumber \\
 & \overset{(b)}{\leq}\left\Vert \bx_{k}^{i}-\by_{k}^{i}\right\Vert ^{2}-2\eta_{k}\frac{2-\eta_{k}L_{\pi_{k}^{i}}}{2+\eta_{k}L_{\pi_{k}^{i}}}\Breg_{\pi_{k}^{i}}(\bx_{k}^{i},\bx_{\star})+\frac{4+\eta_{k}^{2}L_{\pi_{k}^{i}}^{2}}{4-\eta_{k}^{2}L_{\pi_{k}^{i}}^{2}}\eta_{k}L_{\pi_{k}^{i}}\left\Vert \by_{k}^{i}-\bx_{\star}\right\Vert ^{2},\label{eq:RR-cvx-ub-II-descent-3}
\end{align}
where we use $\Breg_{\pi_{k}^{i}}(\by_{k}^{i},\bx_{k}^{i})\geq0$
in $(a)$ and $\Breg_{\pi_{k}^{i}}(\by_{k}^{i},\bx_{\star})\leq\frac{L_{\pi_{k}^{i}}}{2}\left\Vert \by_{k}^{i}-\bx_{\star}\right\Vert ^{2}$
in $(b)$.

Finally, we sum up (\ref{eq:RR-cvx-ub-II-descent-3}) from $i=1$
to $n$, use $\bx_{k+1}^{1}=\bx_{k}^{n+1}$ and $\by_{k}^{n+1}=\by_{k}^{1}=\bx_{\star}$
(see (\ref{eq:virtual-y-prop})), and rearrange terms to obtain
\[
2\eta_{k}\sum_{i=1}^{n}\frac{2-\eta_{k}L_{\pi_{k}^{i}}}{2+\eta_{k}L_{\pi_{k}^{i}}}\Breg_{\pi_{k}^{i}}(\bx_{k}^{i},\bx_{\star})\leq\left\Vert \bx_{k}^{1}-\bx_{\star}\right\Vert ^{2}-\left\Vert \bx_{k+1}^{1}-\bx_{\star}\right\Vert ^{2}+\sum_{i=1}^{n}\frac{4+\eta_{k}^{2}L_{\pi_{k}^{i}}^{2}}{4-\eta_{k}^{2}L_{\pi_{k}^{i}}^{2}}\eta_{k}L_{\pi_{k}^{i}}\left\Vert \by_{k}^{i}-\bx_{\star}\right\Vert ^{2}.
\]
Therefore, we conclude the proof after applying
\begin{eqnarray*}
\frac{2-\eta_{k}\Lmax}{2+\eta_{k}\Lmax}\leq\frac{2-\eta_{k}L_{\pi_{k}^{i}}}{2+\eta_{k}L_{\pi_{k}^{i}}} & \text{and} & \frac{4+\eta_{k}^{2}L_{\pi_{k}^{i}}^{2}}{4-\eta_{k}^{2}L_{\pi_{k}^{i}}^{2}}\leq\frac{4+\eta_{k}^{2}\Lmax^{2}}{4-\eta_{k}^{2}\Lmax^{2}},
\end{eqnarray*}
both of which hold by $L_{\pi_{k}^{i}}\leq\Lmax$.
\end{proof}

\subsection{Proof of Lemma \ref{lem:RR-cvx-ub-II-residual}}

\begin{proof}
Note that we have
\[
L_{\pi_{k}^{i}}\left\Vert \by_{k}^{i}-\bx_{\star}\right\Vert ^{2}\overset{(\ref{eq:virtual-y-prop})}{=}L_{\pi_{k}^{i}}\left\Vert \by_{k}^{i}-\by_{k}^{1}\right\Vert ^{2}\overset{(\ref{eq:virtual-y-def})}{=}\eta_{k}^{2}L_{\pi_{k}^{i}}\left\Vert \sum_{j=1}^{i-1}\nabla f_{\pi_{k}^{j}}(\bx_{\star})\right\Vert ^{2}.
\]
Therefore, we can bound 
\[
\E\left[\sum_{i=1}^{n}L_{\pi_{k}^{i}}\left\Vert \by_{k}^{i}-\bx_{\star}\right\Vert ^{2}\right]=\eta_{k}^{2}\E\left[\sum_{i=1}^{n}L_{\pi_{k}^{i}}\left\Vert \sum_{j=1}^{i-1}\nabla f_{\pi_{k}^{j}}(\bx_{\star})\right\Vert ^{2}\right]\leq\frac{\eta_{k}^{2}n^{2}\Lavg\sigma_{\star}^{2}}{3},
\]
where the last step is due to Lemma E.1 of \citet{pmlr-v235-liu24cg}
(the constant here is slightly better, since $\nabla f(\bx_{\star})=\bzero$
in our setting leads to a provable improvement).
\end{proof}

\section{Rate for $\protect\SGD$\label{sec:SGD-cvx-ub}}

In this section, we present a convergence rate for standard $\SGD$
with epoch-dependent stepsizes. The proof relies on Lemma \ref{lem:RR-cvx-ub-I-descent},
a direct extension of Theorem 8 of \citet{NEURIPS2025_97335d45},
as discussed previously.
\begin{thm}
\label{thm:SGD-cvx-ub}Under Assumptions \ref{assu:minimizer}, \ref{assu:convexity},
and \ref{assu:smoothness}, suppose $\eta_{k}<\frac{2}{\Lmax},\forall k\in\left[K\right]$,
then $\SGD$ (Algorithm \ref{alg:Shuffling-SGD}, with $\pi_{k}^{i}$
sampled independently and uniformly from $\left[n\right]$) guarantees
that
\[
\E\left[f(\bar{\bx}_{K})-f_{\star}\right]\leq\frac{2\left\Vert \bx_{1}^{1}-\bx_{\star}\right\Vert ^{2}}{n\sum_{k=1}^{K}(2-\eta_{k}\Lmax)\eta_{k}}+\frac{8\sum_{k=1}^{K}\eta_{k}^{2}/(2-\eta_{k}\Lmax)}{\sum_{k=1}^{K}(2-\eta_{k}\Lmax)\eta_{k}}\sigma_{\star}^{2},
\]
where $\bar{\bx}_{K}$ is defined in (\ref{eq:avg-x}).
\end{thm}
\begin{proof}
We invoke Lemma \ref{lem:RR-cvx-ub-I-descent} and drop the nonpositive
term $-\frac{2-\eta_{k}\Lmax}{2}\eta_{k}\Breg_{\pi_{k}^{i}}(\bx_{\star},\bx_{k}^{i})$
to have, for any $k\in\left[K\right]$ and $i\in\left[n\right]$,
\begin{equation}
\frac{2-\eta_{k}\Lmax}{2}\eta_{k}\Breg_{\pi_{k}^{i}}(\bx_{k}^{i},\bx_{\star})\leq\left\Vert \bx_{k}^{i}-\bx_{\star}\right\Vert ^{2}-\left\Vert \bx_{k}^{i+1}-\bx_{\star}\right\Vert ^{2}+\frac{2+\eta_{k}\Lmax}{2-\eta_{k}\Lmax}\eta_{k}^{2}\left\Vert \nabla f_{\pi_{k}^{i}}(\bx_{\star})\right\Vert ^{2}-2\eta_{k}\left\langle \nabla f_{\pi_{k}^{i}}(\bx_{\star}),\bx_{k}^{i}-\bx_{\star}\right\rangle .\label{eq:SGD-cvx-ub-1}
\end{equation}
Note that $\pi_{k}^{i}$ is now sampled uniformly from $\left[n\right]$
and is independent from the history. Combining with the fact that
$\nabla f(\bx_{\star})=\bzero$, this implies that
\begin{eqnarray*}
\E\left[\Breg_{\pi_{k}^{i}}(\bx_{k}^{i},\bx_{\star})\right]=\E\left[f(\bx_{k}^{i})-f_{\star}\right], & \E\left[\left\Vert \nabla f_{\pi_{k}^{i}}(\bx_{\star})\right\Vert ^{2}\right]=\sigma_{\star}^{2}, & \E\left[\left\langle \nabla f_{\pi_{k}^{i}}(\bx_{\star}),\bx_{k}^{i}-\bx_{\star}\right\rangle \right]=0.
\end{eqnarray*}
Therefore, after taking expectations on both sides of (\ref{eq:SGD-cvx-ub-1}),
multiplying both sides by $2$, and using $\eta_{k}\Lmax\leq2$, we
obtain that
\[
(2-\eta_{k}\Lmax)\eta_{k}\E\left[f(\bx_{k}^{i})-f_{\star}\right]\leq2\E\left[\left\Vert \bx_{k}^{i}-\bx_{\star}\right\Vert ^{2}\right]-2\E\left[\left\Vert \bx_{k}^{i+1}-\bx_{\star}\right\Vert ^{2}\right]+\frac{8\eta_{k}^{2}}{2-\eta_{k}\Lmax}\sigma_{\star}^{2}.
\]
Finally, a direct calculation yields the rate
\[
\E\left[f(\bar{\bx}_{K})-f_{\star}\right]\leq\frac{2\left\Vert \bx_{1}^{1}-\bx_{\star}\right\Vert ^{2}}{n\sum_{k=1}^{K}(2-\eta_{k}\Lmax)\eta_{k}}+\frac{8\sum_{k=1}^{K}\eta_{k}^{2}/(2-\eta_{k}\Lmax)}{\sum_{k=1}^{K}(2-\eta_{k}\Lmax)\eta_{k}}\sigma_{\star}^{2}.
\]
\end{proof}

\end{document}